\documentclass[a4]{article}

\usepackage{amsmath}
\usepackage{amssymb}
\usepackage{amsthm}

\begin{document}

\bibliographystyle{amsalpha}
\newtheorem{Theorem}{Theorem}[section]
\newtheorem{Lemma}{Lemma}[section]
\newtheorem{Remark}{Remark}[section]
\newtheorem{Corollary}{Corollary}[section]
\newtheorem{Proposition}{Proposition}[section]
\newtheorem{Example}{Example}[section]
\newtheorem{Definition}{Definition}[section]
\newtheorem{Problem}{Problem}[section]
\newtheorem{Proof}{Proof}[section]

\title{Second order polynomial Hamiltonian systems with ${\tilde W}(E_6^{(1)}),{\tilde W}(E_7^{(1)})$ and $W(E_8^{(1)})$-symmetry \\}

\author{By\\
Yusuke Sasano}
\maketitle

\begin{abstract} We find and study a six (resp. seven, eight)-parameter family of polynomial Hamiltonian systems of second order, respectively. This system admits the affine Weyl group symmetry of type $E_6^{(1)}$ (resp. $E_7^{(1)}, E_8^{(1)}$) as the group of its B{\"a}cklund transformations. Each system is the first example which gave second-order polynomial Hamiltonian system with ${\tilde W}(E_6^{(1)})$ (resp. ${\tilde W}(E_7^{(1)}), W(E_8^{(1)})$)-symmetry.
 We also show that its space of initial conditions $S$ is obtained by gluing eight (resp. nine, ten) copies of ${\Bbb C}^2$ via the birational and symplectic transformations.


\end{abstract}

\section{Main results of the system with ${\tilde W}(E_6^{(1)})$-symmetry}

In this paper, we find and study a 6-parameter family of polynomial Hamiltonian systems of second order. This system admits extended affine Weyl group symmetry of type $E_6^{(1)}$ (see Figure 1) as the group of its B{\"a}cklund transformations. This system is the first example which gave second-order polynomial Hamiltonian system with ${\tilde W}(E_6^{(1)})$-symmetry.

By eliminating $p$ or $q$, we obtain the second-order ordinary differential equation. However, its form is not normal \rm{(cf. \cite{Cosgrove1,Cosgrove2})}.

We also show that after a series of explicit blowing-ups at nine points including the infinitely near points of the Hirzebruch surface ${\Sigma_2}$ (see Figure 2) and blowing-down along the $(-1)$-curve $H' \cong {\Bbb P}^1$ to a nonsingular point (see Figure 3), we obtain the rational surface $\tilde{S}$ and a birational morphism 
$$
\varphi:\tilde{S} \leftarrow S_9 \rightarrow \cdots \rightarrow S_1 \rightarrow {\Sigma_2}.
$$
Here, the symbol $H'$ denotes the strict transform of $H$, each ${E}_i$ denotes the exceptional divisors, and $-K_{\Sigma_2}=2H, \ H \cong {\Bbb P}^1, \ (H)^2=2$. In order to obtain a minimal compactification of the space of initial conditions, we must blow down along the $(-1)$-curve $H'$.

Its canonical divisor $K_{\tilde{S}}$ of $\tilde{S}$ is given by
\begin{align}
\begin{split}
K_{\tilde{S}}&=- \sum_{i=1}^{3} {E}_i, \quad ({E}_i)^2=-2, \quad {E}_i \cong {\Bbb P}^1, \quad E_1 \cap E_2 \cap E_3 \not=\varnothing, \quad (E_j,E_k)=1 \quad (j \not= k).
\end{split}
\end{align}
This configuration of $(-K_{\tilde{S}})_{red}$ is type ${A_2^{(1)}}^{*}$ (see Figure 1). This type of rational surface $\tilde S$ does not appear in the list of Painlev\'e equations (see \cite{Sakai}).

The space of initial conditions $S$ is obtained by gluing eight copies of ${\Bbb C}^2$
\begin{align}
\begin{split}
S&={\tilde{S}}-(-{K_{\tilde{S}}})_{red}\\
&={\Bbb C}^2 \cup \bigcup_{i=0}^{6} {\Bbb C}^2
\end{split}
\end{align}
via the birational and symplectic transformations $r_j$ (see Theorem 2.1). 

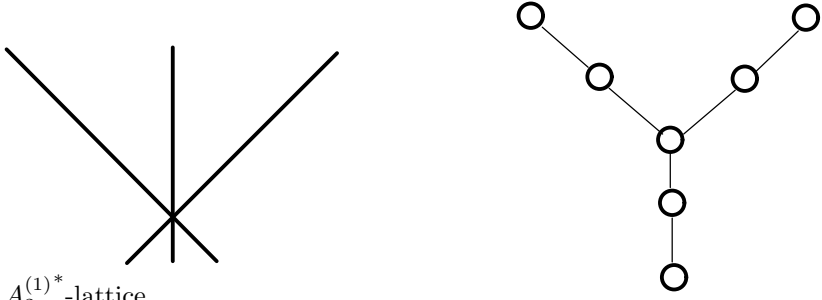
\begin{figure}[h]
\unitlength 0.1in
\begin{picture}( 42.5600, 15.7000)( 11.7000,-17.1000)
%
\special{pn 20}%
\special{pa 1180 380}%
\special{pa 2280 1490}%
\special{fp}%
\special{pa 2050 370}%
\special{pa 2050 1490}%
\special{fp}%
%
\special{pn 20}%
\special{pa 1810 1500}%
\special{pa 2910 400}%
\special{fp}%
\put(11.7000,-17.5000){\makebox(0,0)[lb]{${A_2^{(1)}}^{*}$-lattice}}%
%
\special{pn 20}%
\special{ar 3924 204 64 64  0.0000000 6.2831853}%
%
\special{pn 20}%
\special{ar 4284 524 64 64  0.0000000 6.2831853}%
%
\special{pn 8}%
\special{pa 3984 264}%
\special{pa 4224 474}%
\special{fp}%
%
\special{pn 8}%
\special{pa 4334 584}%
\special{pa 4614 824}%
\special{fp}%
%
\special{pn 20}%
\special{ar 4654 854 64 64  0.0000000 6.2831853}%
%
\special{pn 20}%
\special{ar 5044 534 64 64  0.0000000 6.2831853}%
%
\special{pn 20}%
\special{ar 5364 214 64 64  0.0000000 6.2831853}%
%
\special{pn 20}%
\special{ar 4664 1184 64 64  0.0000000 6.2831853}%
%
\special{pn 20}%
\special{ar 4674 1574 64 64  0.0000000 6.2831853}%
%
\special{pn 8}%
\special{pa 4724 824}%
\special{pa 4994 594}%
\special{fp}%
%
\special{pn 8}%
\special{pa 5094 504}%
\special{pa 5324 284}%
\special{fp}%
%
\special{pn 8}%
\special{pa 4654 924}%
\special{pa 4654 1104}%
\special{fp}%
%
\special{pn 8}%
\special{pa 4664 1264}%
\special{pa 4664 1494}%
\special{fp}%
\put(37.6000,-18.8000){\makebox(0,0)[lb]{Dynkin diagram of type $E_6^{(1)}$}}%
\end{picture}%
\label{fig:E6figure3}
\caption{${A_2^{(1)}}^{*}$-lattice and Dynkin diagram of type $E_6^{(1)}$}
\end{figure}

\begin{center}
\begin{tabular}{|c||c|c|c|c|c|c|c|c|} \hline 
Painlev\'e equations &  PVI & PV & PIV & ${PIII}^{D_6^{(1)}}$ & ${PIII}^{D_7^{(1)}}$ & ${PIII}^{D_8^{(1)}}$ & PII & PI\\ \hline
Type of surface &  $D_4^{(1)}$ & $D_5^{(1)}$ & $E_6^{(1)}$ & $D_6^{(1)}$ & $D_7^{(1)}$ & $D_8^{(1)}$ & $E_7^{(1)}$ & $E_8^{(1)}$ \\ \hline
Symmetry &  $D_4^{(1)}$ & $A_3^{(1)}$ & $A_2^{(1)}$ & $C_2^{(1)}$ & $A_1^{(1)}$ & none & $A_1^{(1)}$ & none \\ \hline
\end{tabular}
\end{center}

This system is the first example whose minimal model $\tilde S$ is the rational surface of type ${A_2^{(1)}}^{*}$.

The author believes that this system can be obtained by holonomic deformation of the 3rd-order linear ordinary differential equation
\begin{align}\label{34}
\begin{split}
&\frac{d^3y}{dx^3}+a_1(x)\frac{d^2y}{dx^2}+a_2(x)\frac{dy}{dx}+a_3(x)y=0 \quad (a_i \in {\Bbb C}(x))
\end{split}
\end{align}
satisfying the Riemann scheme:
\begin{equation}\label{scheme}
\begin{pmatrix}
x=0 & x=\frac{1}{\varepsilon} & x=q & x=\infty\\
0 & 0 & 0 & \alpha_0\\
\alpha_2 & \alpha_4 & 1 & \alpha_0+\alpha_5\\
\alpha_1+\alpha_2 & \alpha_3+\alpha_4 & 3 & \alpha_0+\alpha_5+\alpha_6
\end{pmatrix} \; ,
\end{equation}
where $x=q$ is an apparent singular point. In this case, $\varepsilon=1$.

It is still an open question whether the equation \eqref{34} satisfying \eqref{scheme} tends to the equation \eqref{34} satisfying the Riemann scheme:
\begin{equation}\label{scheme2}
\begin{pmatrix}
x=0 & x=q & x=\infty & x=\infty\\
0 & 0 & 0 & \alpha_4\\
\alpha_2 & 1 & 1 & \alpha_3\\
\alpha_1+\alpha_2 & 3 & t & \alpha_0
\end{pmatrix}
\end{equation}
as $\varepsilon \rightarrow 0$, where $x=\infty$ is an irregular singular point with Poincar\'e rank 1 (cf. \cite{Mazzocco,Joshi}).

\section{Holomorphy}
\begin{Theorem}
Let us consider a polynomial Hamiltonian system with Hamiltonian $I \in {\Bbb C}(t)[q,p]$. We assume that

$(A1)$ $deg(I)=7$ with respect to $q,p$.

$(A2)$ This system becomes again a polynomial Hamiltonian system in each coordinate $r_i \ (i=0,1,2,3,4)${\rm : \rm}
\begin{align*}
&r_0:x_0=1/q, \ y_0=-(qp+\alpha_0)q,\\
&r_1:x_1=-(pq-(\alpha_1+\alpha_2))p, \ y_1=1/p,\\
&r_2:x_2=-(pq-\alpha_2)p, \ y_2=1/p,\\
&r_3:x_3=-(p(q-1)-(\alpha_3+\alpha_4))p, \ y_3=1/p,\\
&r_4:x_4=-(p(q-1)-\alpha_4)p, \ y_4=1/p.
\end{align*}

$(A3)$ In addition to the assumption $(A2)$, the Hamiltonian system in the coordinate $r_0$ becomes again a polynomial Hamiltonian system in the coordinates $r_5,r_6${\rm : \rm}
\begin{align*}
&r_5:x_5=-(x_0y_0-\alpha_5)y_0, \ y_5=1/y_0,\\
&r_6:x_6=-(x_0y_0-(\alpha_5+\alpha_6))y_0, \ y_6=1/y_0.
\end{align*}
Then such a system coincides with the system
\begin{align}\label{11}
\begin{split}
&\frac{dq}{dt}=\frac{\partial I}{\partial p}, \quad \frac{dp}{dt}=-\frac{\partial I}{\partial q},\\
I:=&(q-1)^2q^2p^3-q(q-1)\{(\alpha_1+2\alpha_2+\alpha_3+2\alpha_4)q-\alpha_1-2\alpha_2\}p^2\\
&+[\{-3\alpha_0^2-2\alpha_0(\alpha_1+2\alpha_2+\alpha_3+2\alpha_4)-3\alpha_0 \alpha_5-\alpha_5(\alpha_1+2\alpha_2+\alpha_3+2\alpha_4+\alpha_5)\}q^2\\
&\{-3\alpha_0^2-\alpha_2^2+2\alpha_0(\alpha_1+2\alpha_2+\alpha_3+2\alpha_4)+3\alpha_0 \alpha_5+2\alpha_2 \alpha_5+\alpha_1(\alpha_5-\alpha_2)\\
&+(\alpha_4+\alpha_5)(\alpha_3+\alpha_4+\alpha_5)\}q+\alpha_2(\alpha_1+\alpha_2)]p+\alpha_0(\alpha_0+\alpha_5)(\alpha_0+\alpha_5+\alpha_6)q,
\end{split}
\end{align}
where the constant parameters $\alpha_i$ satisfy the relation:
\begin{equation}
3\alpha_0+\alpha_1+2\alpha_2+\alpha_3+2\alpha_4+2\alpha_5+\alpha_6=0.
\end{equation}
\end{Theorem}
Since each transformation $r_i$ is symplectic, the system \eqref{11} is transformed into a Hamiltonian system, whose Hamiltonian may have poles. It is remarkable that the transformed system becomes again a polynomial system for any $i=0,1,\ldots,6$.

The holomorphy conditions $(A2),(A3)$ are new. Theorem 2.1 can be checked by a direct calculation.

\begin{Proposition}
The Hamiltonian $I$ is its first integral.
\end{Proposition}

\begin{Remark}
For the Hamiltonian system in each coordinate system $(x_i,y_i) \ (i=0,1,\ldots,6)$ given by $(A2)$ and $(A3)$ in Theorem 2.1, by eliminating $x_i$ or $y_i$, we obtain the second-order ordinary differential equation. However, its form is not normal \rm{(cf. \cite{Cosgrove1,Cosgrove2})}.
\end{Remark}

\section{Symmetry}

\begin{Theorem}
The system \eqref{11} admits extended affine Weyl group symmetry of type $E_6^{(1)}$ as the group of its B{\"a}cklund transformations whose generators $s_i, \ i=0,1,\ldots,6, \ {\pi}_j, \ j=1,2,3$ are explicitly given as follows{\rm : \rm}with the notation $(*):=(q,p,t;\alpha_0,\alpha_1,\ldots,\alpha_6)$,
\begin{align*}
        s_{0}: (*) &\rightarrow \left(q+\frac{\alpha_0}{p},p,t;-\alpha_0,\alpha_1,\alpha_0+\alpha_2,\alpha_3,\alpha_4+\alpha_0,\alpha_5+\alpha_0,\alpha_6 \right),\\
        s_{1}: (*) &\rightarrow (q,p,t;\alpha_0,-\alpha_1,\alpha_2+\alpha_1,\alpha_3,\alpha_4,\alpha_5,\alpha_6), \\
        s_{2}: (*) &\rightarrow \left(q,p-\frac{\alpha_2}{q},t;\alpha_0+\alpha_2,\alpha_1+\alpha_2,-\alpha_2,\alpha_3,\alpha_4,\alpha_5,\alpha_6 \right), \\
        s_{3}: (*) &\rightarrow (q,p,t;\alpha_0,\alpha_1,\alpha_2,-\alpha_3,\alpha_4+\alpha_3,\alpha_5,\alpha_6), \\
        s_{4}: (*) &\rightarrow \left(q,p-\frac{\alpha_4}{q-1},t;\alpha_0+\alpha_4,\alpha_1,\alpha_2,\alpha_3+\alpha_4,-\alpha_4,\alpha_5,\alpha_6 \right),\\
        s_{5}: (*) &\rightarrow (q,p,t;\alpha_0+\alpha_5,\alpha_1,\alpha_2,\alpha_3,\alpha_4,-\alpha_5,\alpha_6+\alpha_5),\\
        s_{6}: (*) &\rightarrow (q,p,t;\alpha_0,\alpha_1,\alpha_2,\alpha_3,\alpha_4,\alpha_5+\alpha_6,-\alpha_6),\\
        \pi_{1}: (*) &\rightarrow (1-q,-p,1-t;\alpha_0,\alpha_3,\alpha_4,\alpha_1,\alpha_2,\alpha_5,\alpha_6),\\
        \pi_{2}: (*) &\rightarrow \left(\frac{1}{q},-(qp+\alpha_0)q,-t;\alpha_0,\alpha_6,\alpha_5,\alpha_3,\alpha_4,\alpha_2,\alpha_1 \right),\\
        \pi_{3}: (*) &\rightarrow \left(\frac{q}{q-1},-(q-1)((q-1)p+\alpha_0),2-t;\alpha_0,\alpha_1,\alpha_2,\alpha_6,\alpha_5,\alpha_4,\alpha_3 \right),
\end{align*}
where $\pi_j, \ j=1,2,3$ are Dynkin diagram automorphisms of type $E_6^{(1)}$.
\end{Theorem}
Theorem 3.1 can be checked by a direct calculation.

\section{Space of initial conditions}

\begin{Theorem}\label{3.1}
After a series of explicit blowing-ups at nine points including the infinitely near points of $\Sigma_2$ and blowing-down along the $(-1)$-curve ${D^{(0)}}' \cong {\Bbb P}^1$, we obtain the rational surface $\tilde{S}$ of the system \eqref{11} and a birational morphism $\varphi:\tilde{S} \cdots \rightarrow \Sigma_2$. Its canonical divisor $K_{\tilde{S}}$ of $\tilde{S}$ is given by
\begin{align}
\begin{split}
K_{\tilde{S}}&=-D_{0}^{(1)}-D_{1}^{(1)}-D_{\infty}^{(1)}, \quad (D_{\nu}^{(1)})^2=-2, \ D_{\nu}^{(1)} \cong {\Bbb P}^1,
\end{split}
\end{align}
where the symbol ${D^{(0)}}'$ denotes the strict transform of $D^{(0)}$, $D_{\nu}^{(1)}$ denote the exceptional divisors and $-K_{\Sigma_2}=2D^{(0)}, \ D^{(0)} \cong {\Bbb P}^1, \ (D^{(0)})^2=2$.
\end{Theorem}

\begin{Theorem}
The space of initial conditions $S$ of the system \eqref{11} is obtained by gluing eight copies of ${\Bbb C}^2$:
\begin{align}
\begin{split}
S&={\tilde{S}}-(-{K_{\tilde{S}}})_{red}\\
&={\Bbb C}^2 \cup \bigcup_{i=0}^{6} U_j,\\
&{\Bbb C}^2 \ni (q,p), \quad U_j \cong {\Bbb C}^2 \ni (x_j,y_j) \ (j=0,1,\ldots,6)
\end{split}
\end{align}
via the birational and symplectic transformations $r_j$ \rm{(see Theorem 2.1)}.
\end{Theorem}

{\bf Proof of Theorems 4.1 and 4.2.}

\begin{figure}[h]
\unitlength 0.1in
\begin{picture}( 24.8000, 16.3000)( 19.9000,-20.0000)
%
\special{pn 8}%
\special{pa 1990 590}%
\special{pa 4320 590}%
\special{fp}%
%
\special{pn 8}%
\special{pa 2010 1810}%
\special{pa 4310 1810}%
\special{fp}%
%
\special{pn 8}%
\special{pa 2410 410}%
\special{pa 2410 2000}%
\special{fp}%
%
\special{pn 8}%
\special{pa 3990 410}%
\special{pa 3990 1990}%
\special{fp}%
\put(44.7000,-13.0000){\makebox(0,0)[lb]{${\Sigma}_{2}$}}%
%
\special{pn 20}%
\special{pa 2410 590}%
\special{pa 2710 590}%
\special{fp}%
\special{sh 1}%
\special{pa 2710 590}%
\special{pa 2644 570}%
\special{pa 2658 590}%
\special{pa 2644 610}%
\special{pa 2710 590}%
\special{fp}%
%
\special{pn 20}%
\special{pa 2410 590}%
\special{pa 2410 870}%
\special{fp}%
\special{sh 1}%
\special{pa 2410 870}%
\special{pa 2430 804}%
\special{pa 2410 818}%
\special{pa 2390 804}%
\special{pa 2410 870}%
\special{fp}%
\put(24.5000,-5.4000){\makebox(0,0)[lb]{$q$}}%
\put(21.4000,-8.3000){\makebox(0,0)[lb]{$p$}}%
\put(44.7000,-18.7000){\makebox(0,0)[lb]{$D^{(0)} \cong {\Bbb P}^1$}}%
%
\special{pn 20}%
\special{pa 2410 1810}%
\special{pa 2410 1540}%
\special{fp}%
\special{sh 1}%
\special{pa 2410 1540}%
\special{pa 2390 1608}%
\special{pa 2410 1594}%
\special{pa 2430 1608}%
\special{pa 2410 1540}%
\special{fp}%
%
\special{pn 20}%
\special{pa 2400 1810}%
\special{pa 2680 1810}%
\special{fp}%
\special{sh 1}%
\special{pa 2680 1810}%
\special{pa 2614 1790}%
\special{pa 2628 1810}%
\special{pa 2614 1830}%
\special{pa 2680 1810}%
\special{fp}%
%
\special{pn 20}%
\special{pa 3990 590}%
\special{pa 3990 870}%
\special{fp}%
\special{sh 1}%
\special{pa 3990 870}%
\special{pa 4010 804}%
\special{pa 3990 818}%
\special{pa 3970 804}%
\special{pa 3990 870}%
\special{fp}%
%
\special{pn 20}%
\special{pa 3980 590}%
\special{pa 3710 590}%
\special{fp}%
\special{sh 1}%
\special{pa 3710 590}%
\special{pa 3778 610}%
\special{pa 3764 590}%
\special{pa 3778 570}%
\special{pa 3710 590}%
\special{fp}%
%
\special{pn 20}%
\special{pa 3980 1810}%
\special{pa 3980 1530}%
\special{fp}%
\special{sh 1}%
\special{pa 3980 1530}%
\special{pa 3960 1598}%
\special{pa 3980 1584}%
\special{pa 4000 1598}%
\special{pa 3980 1530}%
\special{fp}%
%
\special{pn 20}%
\special{pa 3970 1800}%
\special{pa 3720 1800}%
\special{fp}%
\special{sh 1}%
\special{pa 3720 1800}%
\special{pa 3788 1820}%
\special{pa 3774 1800}%
\special{pa 3788 1780}%
\special{pa 3720 1800}%
\special{fp}%
\put(36.3000,-5.4000){\makebox(0,0)[lb]{$z_1$}}%
\put(40.8000,-8.2000){\makebox(0,0)[lb]{$w_1$}}%
\put(20.9000,-17.3000){\makebox(0,0)[lb]{$w_2$}}%
\put(24.6000,-20.0000){\makebox(0,0)[lb]{$z_2$}}%
\put(40.8000,-17.3000){\makebox(0,0)[lb]{$w_3$}}%
\put(36.7000,-20.0000){\makebox(0,0)[lb]{$z_3$}}%
\end{picture}%
\label{fig:E6figure1}
\caption{Hirzebruch surface ${\Sigma_2}$}
\end{figure}
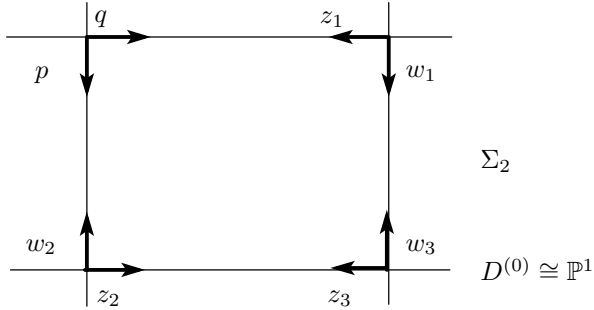
At first, we take the Hirzebruch surface ${\Sigma_2}$ which is obtained by gluing four copies of ${\Bbb C}^2$ via the following identification.
\begin{align}
\begin{split}
&U_j \cong {\Bbb C}^2 \ni (z_j,w_j) \ (j=0,1,2,3)\\
&z_0=q, \ w_0=p, \quad z_1=\frac{1}{q}, \ w_1=-(qp+\alpha_0)q,\\
&z_2=z_0, \ w_2=\frac{1}{w_0}, \quad z_3=z_1, \ w_3=\frac{1}{w_1}.
\end{split}
\end{align}
We define a divisor $D^{(0)}$ on ${\Sigma_2}$:
\begin{equation}
D^{(0)}=\{(z_2,w_2) \in U_2|w_2=0\} \cup \{(z_3,w_3) \in U_3|w_3=0\} \cong {\Bbb P}^1.
\end{equation}
The self-intersection number of $D^{(0)}$ is given by
\begin{equation}
(D^{(0)})^2=2.
\end{equation}

By a direct calculation, we see that the system \eqref{11} has three accessible singular points $a_{\nu}^{(0)} \in D^{(0)} \quad (\nu=0,1,\infty)$:
\begin{align}
\begin{split}
&a_{\nu}^{(0)}=\{(z_2,w_2)=(\nu,0)\} \in U_2 \cap D^{(0)} \ (\nu=0,1),\\
&a_{\infty}^{(0)}=\{(z_3,w_3)=(0,0)\} \in U_3 \cap D^{(0)}.
\end{split}
\end{align}
We perform blowing-ups in ${\Sigma_2}$ at $a_{\nu}^{(0)}$, and let $D_{\nu}^{(1)}$ be the exceptional curves of the blowing-ups at $a_{\nu}^{(0)}$ for $\nu=0,1,\infty$. We can take three coordinate systems $(u_{\nu},v_{\nu})$ around the points at infinity of the exceptional curves $D_{\nu}^{(1)} \quad (\nu=0,1,\infty)$, where
\begin{align}
\begin{split}
&(u_{\nu},v_{\nu})=\left(\frac{z_2-\nu}{w_2},w_2 \right) \ (\nu=0,1),\\
&(u_{\infty},v_{\infty})=\left(\frac{z_3}{w_3},w_3 \right).
\end{split}
\end{align}

\begin{figure}
\unitlength 0.1in
\begin{picture}( 23.6000, 86.5200)( 21.1000,-87.7000)
\put(44.5000,-4.6100){\makebox(0,0)[lb]{$D^{(0)} \cong {\Bbb P}^1, \quad (D^{(0)})^2=2$}}%
%
\special{pn 8}%
\special{pa 2110 346}%
\special{pa 4140 346}%
\special{fp}%
%
\special{pn 20}%
\special{sh 0.600}%
\special{ar 2420 346 20 22  0.0000000 6.2831853}%
%
\special{pn 20}%
\special{sh 0.600}%
\special{ar 3080 346 20 22  0.0000000 6.2831853}%
%
\special{pn 20}%
\special{sh 0.600}%
\special{ar 3760 346 20 22  0.0000000 6.2831853}%
\put(22.9000,-2.8800){\makebox(0,0)[lb]{$a_0^{(0)}$}}%
\put(29.6000,-2.8800){\makebox(0,0)[lb]{$a_1^{(0)}$}}%
\put(36.3000,-2.8800){\makebox(0,0)[lb]{$a_{\infty}^{(0)}$}}%
%
\special{pn 20}%
\special{pa 3086 1138}%
\special{pa 3086 576}%
\special{fp}%
\special{sh 1}%
\special{pa 3086 576}%
\special{pa 3066 644}%
\special{pa 3086 630}%
\special{pa 3106 644}%
\special{pa 3086 576}%
\special{fp}%
\put(31.7500,-10.3700){\makebox(0,0)[lb]{Blow up at $a_0^{(0)},a_1^{(0)}$ and $a_{\infty}^{(0)}$}}%
%
\special{pn 8}%
\special{pa 2120 2736}%
\special{pa 4150 2736}%
\special{fp}%
%
\special{pn 8}%
\special{pa 2430 3068}%
\special{pa 2430 1224}%
\special{dt 0.045}%
%
\special{pn 20}%
\special{sh 0.600}%
\special{ar 2430 2276 20 22  0.0000000 6.2831853}%
%
\special{pn 20}%
\special{sh 0.600}%
\special{ar 2430 1670 20 22  0.0000000 6.2831853}%
\put(21.1000,-18.0000){\makebox(0,0)[lb]{$a_1^{(1)}$}}%
\put(21.1000,-23.4700){\makebox(0,0)[lb]{$a_2^{(1)}$}}%
%
\special{pn 8}%
\special{pa 3080 3054}%
\special{pa 3080 1210}%
\special{dt 0.045}%
%
\special{pn 20}%
\special{sh 0.600}%
\special{ar 3080 2262 20 22  0.0000000 6.2831853}%
%
\special{pn 20}%
\special{sh 0.600}%
\special{ar 3080 1656 20 22  0.0000000 6.2831853}%
\put(27.6000,-17.8600){\makebox(0,0)[lb]{$a_3^{(1)}$}}%
\put(27.6000,-23.3300){\makebox(0,0)[lb]{$a_4^{(1)}$}}%
%
\special{pn 8}%
\special{pa 3780 3068}%
\special{pa 3780 1224}%
\special{dt 0.045}%
%
\special{pn 20}%
\special{sh 0.600}%
\special{ar 3780 2276 20 22  0.0000000 6.2831853}%
%
\special{pn 20}%
\special{sh 0.600}%
\special{ar 3780 1670 20 22  0.0000000 6.2831853}%
\put(34.6000,-18.0000){\makebox(0,0)[lb]{$a_5^{(1)}$}}%
\put(34.6000,-23.4700){\makebox(0,0)[lb]{$a_6^{(1)}$}}%
%
\special{pn 20}%
\special{pa 3086 3730}%
\special{pa 3086 3168}%
\special{fp}%
\special{sh 1}%
\special{pa 3086 3168}%
\special{pa 3066 3236}%
\special{pa 3086 3222}%
\special{pa 3106 3236}%
\special{pa 3086 3168}%
\special{fp}%
\put(31.7500,-36.2900){\makebox(0,0)[lb]{Blow up at $a_j^{(1)}$}}%
%
\special{pn 8}%
\special{pa 2110 5574}%
\special{pa 4140 5574}%
\special{fp}%
%
\special{pn 8}%
\special{pa 2420 5904}%
\special{pa 2420 4062}%
\special{fp}%
%
\special{pn 8}%
\special{pa 2320 5256}%
\special{pa 2670 4926}%
\special{dt 0.045}%
%
\special{pn 8}%
\special{pa 2320 4724}%
\special{pa 2620 4334}%
\special{dt 0.045}%
%
\special{pn 8}%
\special{pa 3070 5918}%
\special{pa 3070 4076}%
\special{fp}%
%
\special{pn 8}%
\special{pa 2970 5270}%
\special{pa 3320 4940}%
\special{dt 0.045}%
%
\special{pn 8}%
\special{pa 2970 4738}%
\special{pa 3270 4350}%
\special{dt 0.045}%
%
\special{pn 8}%
\special{pa 3760 5934}%
\special{pa 3760 4090}%
\special{fp}%
%
\special{pn 8}%
\special{pa 3660 5286}%
\special{pa 4010 4954}%
\special{dt 0.045}%
%
\special{pn 8}%
\special{pa 3660 4752}%
\special{pa 3960 4364}%
\special{dt 0.045}%
\put(44.6000,-28.3700){\makebox(0,0)[lb]{$(D^{(0)})^2=-1$}}%
\put(21.5000,-40.6100){\makebox(0,0)[lb]{$(D_{0}^{(1)})^2=-3$}}%
\put(28.1000,-61.7800){\makebox(0,0)[lb]{$(D_{1}^{(1)})^2=-3$}}%
\put(35.1000,-40.7500){\makebox(0,0)[lb]{$(D_{\infty}^{(1)})^2=-3$}}%
%
\special{pn 20}%
\special{pa 3090 6294}%
\special{pa 3090 6798}%
\special{fp}%
\special{sh 1}%
\special{pa 3090 6798}%
\special{pa 3110 6730}%
\special{pa 3090 6744}%
\special{pa 3070 6730}%
\special{pa 3090 6798}%
\special{fp}%
\put(32.0000,-66.5300){\makebox(0,0)[lb]{Blow down along the $(-1)$-curve $D^{(0)}$ to a nonsingular point}}%
%
\special{pn 20}%
\special{pa 2210 7158}%
\special{pa 3310 8756}%
\special{fp}%
\special{pa 3080 7142}%
\special{pa 3080 8756}%
\special{fp}%
%
\special{pn 20}%
\special{pa 2840 8770}%
\special{pa 3940 7186}%
\special{fp}%
%
\special{pn 8}%
\special{pa 2340 7604}%
\special{pa 2620 6998}%
\special{dt 0.045}%
%
\special{pn 8}%
\special{pa 2590 7906}%
\special{pa 2780 7474}%
\special{dt 0.045}%
%
\special{pn 8}%
\special{pa 2970 7502}%
\special{pa 3260 7214}%
\special{dt 0.045}%
%
\special{pn 8}%
\special{pa 2980 7892}%
\special{pa 3280 7604}%
\special{dt 0.045}%
%
\special{pn 8}%
\special{pa 3600 7230}%
\special{pa 3910 7604}%
\special{dt 0.045}%
%
\special{pn 8}%
\special{pa 3510 7546}%
\special{pa 3760 8006}%
\special{dt 0.045}%
\put(44.7000,-79.4900){\makebox(0,0)[lb]{${A_2^{(1)}}^{*}$-lattice}}%
\put(38.5000,-85.9700){\makebox(0,0)[lb]{Each bold line denotes $(-2)$-curve}}%
\end{picture}%
\label{fig:E6figure2}
\caption{Resolution of accessible singular points}
\end{figure}

Note that $\{(u_{\nu},v_{\nu})|v_{\nu}=0\} \subset D_{\nu}^{(1)}$ for $\nu=0,1,\infty$. By a direct calculation, we see that the system \eqref{E811} has six accessible singular points $a_{\nu}^{(1)}$ for $\nu=1,2,3,4,5,6$ in $D_{\nu}^{(1)} \cong {\Bbb P}^1 \ (\nu=0,1,\infty)$.
\begin{align}
\begin{split}
&a_{1}^{(1)}=\{(u_{0},v_{0})=(\alpha_2+\alpha_1,0)\} \in D_{0}^{(1)}, \quad a_{2}^{(1)}=\{(u_{0},v_{0})=(\alpha_2,0)\} \in D_{0}^{(1)},\\
&a_{3}^{(1)}=\{(u_{1},v_{1})=(\alpha_3+\alpha_4,0)\} \in D_{1}^{(1)}, \quad a_{4}^{(1)}=\{(u_{1},v_{1})=(\alpha_4,0)\} \in D_{1}^{(1)},\\
&a_{5}^{(1)}=\{(u_{\infty},v_{\infty})=(\alpha_5,0)\} \in D_{\infty}^{(1)}, \quad a_{6}^{(1)}=\{(u_{\infty},v_{\infty})=(\alpha_5+\alpha_6,0)\} \in D_{\infty}^{(1)}.
\end{split}
\end{align}
Let us perform blowing-ups at $a_{j}^{(1)}$, and denote $D_{j}^{(2)}$ for the exceptional curves, respectively. We take six coordinate systems $(W_j,V_j)$ around the points at infinity of $D_{j}^{(2)}$ for $j=1,2,3,4,5,6$, where
\begin{align}
\begin{split}
&(W_{1},V_{1})=\left(\frac{u_0-(\alpha_1+\alpha_2)}{v_0},v_0 \right),\\
&(W_{2},V_{2})=\left(\frac{u_0-\alpha_2}{v_0},v_0 \right),\\
&(W_{3},V_{3})=\left(\frac{u_1-(\alpha_3+\alpha_4)}{v_1},v_1 \right),\\
&(W_{4},V_{4})=\left(\frac{u_1-\alpha_4}{v_1},v_1 \right),\\
&(W_{5},V_{5})=\left(\frac{u_{\infty}-\alpha_5}{v_{\infty}},v_{\infty} \right),\\
&(W_{6},V_{6})=\left(\frac{u_{\infty}-(\alpha_5+\alpha_6)}{v_{\infty}},v_{\infty} \right).
\end{split}
\end{align}
For the strict transform of $D^{(0)}$, $D_{\nu}^{(1)}$ and $D_{j}^{(2)}$ by the blowing-ups, we also denote by same symbol, respectively.  Here, the self-intersection number of $D^{(0)}, D_{\nu}^{(1)}$ and $D_{j}^{(2)}$ is given by
\begin{equation}
(D^{(0)})^2=-1. \quad (D_{\nu}^{(1)})^2=-3.
\end{equation}
In order to obtain a minimal compactification of the space of initial conditions, we must blow down along the curve $D^{(0)} \cong {\Bbb P}^1$ to a nonsingular point.  For the strict transform of $D_{\nu}^{(1)}$ and $D_{j}^{(2)}$ by the blowing-down, we also denote by same symbol, respectively. Let ${\tilde S} \cdots \rightarrow {\Sigma_2}$ be the composition of above nine blowing-ups and one blowing-down. Then, we see that the canonical divisor class $K_{{\tilde S}}$ of ${\tilde S}$ is given by
\begin{equation}
K_{{\tilde S}}:=-D_{0}^{(1)}-D_{1}^{(1)}-D_{\infty}^{(1)},
\end{equation}
where the self-intersection number of $D_{\nu}^{(1)} \cong {\Bbb P}^1$ is given by
\begin{equation}
(D_{\nu}^{(1)})^2=-2,
\end{equation}
and
\begin{equation}
D_{0}^{(1)} \cap D_{1}^{(1)} \cap D_{\infty}^{(1)} \not=\varnothing, \quad (D_{j}^{(1)},D_{k}^{(1)})=1 \quad (j \not= k).
\end{equation}
The configuration of the divisor $(-K_{{\tilde S}})_{red}$ on $\tilde S$ is of type ${A_2^{(1)}}^{*}$ (see Figure 3). And we see that ${\tilde S}-(-K_{{\tilde S}})_{red}$ is covered by eight Zariski open sets
\begin{align}
\begin{split}
& \rm{Spec} \ {\Bbb C}[W_{j},V_{j}] \quad (j=1,2,3,4,5,6),\\
& \rm{Spec} \ {\Bbb C}[z_0,w_0],\\
& \rm{Spec} \ {\Bbb C}[z_1,w_1].
\end{split}
\end{align}
The relations between $(W_{j},V_{j})$ and $(x_j,y_j)$ are given by
\begin{equation}
(-W_{j},V_{j})=(x_j,y_j) \quad (j=1,2,3,4,5,6).
\end{equation}
We see that the pole divisor of the symplectic 2-form $dp \wedge dq$ coincides with $(-K_{{\tilde S}})_{red}$. Thus, we have completed the proof of Theorems 4.1 and 4.2. \qed

\section{PVI case}
\begin{Theorem}
Let us consider a polynomial Hamiltonian system with Hamiltonian $G \in {\Bbb C}(t)[q,p]$. We assume that

$(A1)$ $deg(G)=7$ with respect to $q,p$.

$(A2)$ This system becomes again a polynomial Hamiltonian system in each coordinate $rr_i \ (i=0,1,2,3,4)${\rm : \rm}
\begin{align}\label{holoPVI}
\begin{split}
&rr_0:x_0=q+\frac{\alpha_0-\alpha_4}{p}+\frac{t}{p^2},\ y_0=p,\\
&rr_1:x_1=-(qp-(\alpha_1+\alpha_2+\alpha_4))p, \ y_1=\frac{1}{p}, \\
&rr_2:x_2=-(qp-(\alpha_2+\alpha_4))p, \ y_2=\frac{1}{p},\\
&rr_3:x_3=q+\frac{\alpha_3-\alpha_4}{p}+\frac{1}{p^2},\ y_3=p,\\
&rr_4:x_4=-(qp-\alpha_4)p,\ y_4=\frac{1}{p}.
\end{split}
\end{align}
Then such a system coincides with the system
\begin{align}\label{SPVI}
\begin{split}
&\frac{dq}{dt}=\frac{\partial G}{\partial p}, \quad \frac{dp}{dt}=-\frac{\partial G}{\partial q},\\
&G:=-\frac{q^3p^4}{t(t-1)}-\frac{(\alpha_0+\alpha_3-2\alpha_4-1)q^2p^3}{t(t-1)}\\
&-\frac{(t+1)q^2p^2}{t(t-1)}-\frac{(\alpha_1\alpha_2+\alpha_2^2+2\alpha_4-2\alpha_0\alpha_4-2\alpha_3\alpha_4+\alpha_4^2)qp^2}{t(t-1)}\\
&-\frac{\{(\alpha_3-\alpha_4)t+\alpha_0-\alpha_4-1\}qp}{t(t-1)}-\frac{q}{t-1}+\frac{\alpha_4(\alpha_2+\alpha_4)(\alpha_1+\alpha_2+\alpha_4)p}{t(t-1)},
\end{split}
\end{align}
where the constant parameters $\alpha_i$ satisfy the relation:
\begin{equation}
\alpha_0+\alpha_1+2\alpha_2+\alpha_3+\alpha_4=1.
\end{equation}
\end{Theorem}
This Hamiltonian $G$ is equivalent to well-known Hamiltonian $H_{VI}$ of the Painlev\'e VI system by the birational and symplectic transformation $\varphi$
$$
\varphi:Q=-(qp-\alpha_4)p, \quad P=\frac{1}{p},
$$
where $H_{VI}$ is explicitly given by
\begin{align}
\begin{split}
&H_{VI}(q,p,t;\alpha_0,\alpha_1,\alpha_2,\alpha_3,\alpha_4)\\
&=\frac{1}{t(t-1)}[p^2(q-t)(q-1)q-\{(\alpha_0-1)(q-1)q+\alpha_3(q-t)q\\
&+\alpha_4(q-t)(q-1)\}p+\alpha_2(\alpha_1+\alpha_2)(q-t)] \quad (\alpha_0+\alpha_1+2\alpha_2+\alpha_3+\alpha_4=1).
\end{split}
\end{align}
Theorem 5.1 can be checked by a direct calculation.

\begin{Theorem}
The system \eqref{SPVI} is invariant under the following transformations, whose generators $w_i, \ i=0,1,2,3,4$, are given by
\begin{align}\label{D4}
\begin{split}
w_0(q,p,t;\alpha_0,\alpha_1,\dots,\alpha_4) \rightarrow &\left((1-t)\left(q+\frac{\alpha_0-\alpha_4}{p}+\frac{t}{p^2}\right),\frac{p}{1-t},\frac{t}{t-1};\alpha_4,\alpha_1,\alpha_2,\alpha_3,\alpha_0 \right),\\
w_1(q,p,t;\alpha_0,\alpha_1,\dots,\alpha_4) \rightarrow &(q,p,t;\alpha_0,-\alpha_1,\alpha_2+\alpha_1,\alpha_3,\alpha_4),\\
w_2(q,p,t;\alpha_0,\alpha_1,\dots,\alpha_4) \rightarrow &(q,p,t;\alpha_0+\alpha_2,\alpha_1+\alpha_2,-\alpha_2,\alpha_3+\alpha_2,\alpha_4+\alpha_2),\\
w_3(q,p,t;\alpha_0,\alpha_1,\dots,\alpha_4) \rightarrow &\left(-\left(q+\frac{\alpha_3-\alpha_4}{p}+\frac{1}{p^2}\right),-p,1-t;\alpha_0,\alpha_1,\alpha_2,\alpha_4,\alpha_3 \right),\\
w_4(q,p,t;\alpha_0,\alpha_1,\dots,\alpha_4) \rightarrow &\left(q,p-\frac{\alpha_4}{q},t;\alpha_0,\alpha_1,\alpha_2+\alpha_4,\alpha_3,-\alpha_4 \right).
\end{split}
\end{align}
\end{Theorem}
Theorem 5.2 can be checked by a direct calculation.

\section{Main results of the system with ${\tilde W}(E_7^{(1)})$-symmetry}
In this section, by using a relation between holomorphy property and Lax equation, we try to make a second-order polynomial Hamiltonian system with symmetry of the affine Weyl group of type $E_7^{(1)}$. However, for a while, we have not succeeded.

By changing our idea, in the process of construction of the system with ${\tilde W}(E_6^{(1)})$-symmetry, we find the following realtions between symmetry
\begin{align*}
        s_{1}: (q,p,t;\alpha_1,\alpha_2,\alpha_3) &\rightarrow \left(q,p-\frac{\alpha_1}{q},t;-\alpha_1,\alpha_2+\alpha_1,\alpha_3 \right), \\
        s_{2}: (q,p,t;\alpha_1,\alpha_2,\alpha_3) &\rightarrow (q,p,t;\alpha_1+\alpha_2,-\alpha_2,\alpha_3+\alpha_2), \\
        s_{3}: (q,p,t;\alpha_1,\alpha_2,\alpha_3) &\rightarrow (q,p,t;\alpha_1,\alpha_2+\alpha_3,-\alpha_3),
\end{align*}
and holomorphy conditions
\begin{align*}
&r_1:x_1=-(pq-\alpha_1)p, \ y_1=1/p,\\
&r_2:x_2=-(pq-(\alpha_1+\alpha_2))p, \ y_2=1/p,\\
&r_3:x_3=-(pq-(\alpha_1+\alpha_2+\alpha_3))p, \ y_3=1/p.
\end{align*}
By using this key property, we try to make a representation of the affine Weyl group symmetry of type $E_7^{(1)}$ and associated holomorphy conditoins.

In this paper, we find and study a 7-parameter family of polynomial Hamiltonian systems of second order. This system admits extended affine Weyl group symmetry of type $E_7^{(1)}$ as the group of its B{\"a}cklund transformations (see Figure 4). This system is the first example which gave second-order polynomial Hamiltonian systems with ${\tilde W}(E_7^{(1)})$-symmetry.

By eliminating $p$ or $q$, we obtain the second-order ordinary differential equation. However, its form is not normal \rm{(cf. \cite{Cosgrove1,Cosgrove2})}.

We also show that after a series of explicit blowing-ups at ten points including the infinitely near points of the Hirzebruch surface ${\Sigma_2}$ (see Figure 5) and two times blowing-downs along the $(-1)$-curve to a nonsingular point (see Figure 6), respectively, we obtain the rational surface $\tilde{S}$ and a birational morphism 
$$
\varphi:\tilde{S}=S_{12} \leftarrow S_{11} \leftarrow S_{10} \rightarrow \cdots \rightarrow S_1 \rightarrow {\Sigma_2}.
$$
Here, $-K_{\Sigma_2}=2H, \ H \cong {\Bbb P}^1, \ (H)^2=2$. In order to obtain a minimal compactification of the space of initial conditions, we must blow down along the $(-1)$-curves.

Its canonical divisor $K_{\tilde{S}}$ of $\tilde{S}$ is given by
\begin{align}
\begin{split}
K_{\tilde{S}}&=- \sum_{i=1}^{2} {E}_i, \quad ({E}_i)^2=-2, \quad {E}_i \cong {\Bbb P}^1, \quad E_1 \cap E_2  \not=\varnothing, \quad (E_1,E_2)=1.
\end{split}
\end{align}
This configuration of $(-K_{\tilde{S}})_{red}$ is type $A_2$ (see Figure 4). This type of rational surface $\tilde S$ does not appear in the list of Painlev\'e equations (see \cite{Sakai}).

The space of initial conditions $S$ is obtained by gluing nine copies of ${\Bbb C}^2$
\begin{align}
\begin{split}
S&={\tilde{S}}-(-{K_{\tilde{S}}})_{red}\\
&={\Bbb C}^2 \cup \bigcup_{i=0}^{7} {\Bbb C}^2
\end{split}
\end{align}
via the birational and symplectic transformations $r_j$ (see Theorem 7.1).

This system is the first example whose minimal model $\tilde S$ is the rational surface of type $A_2$.

The author believes that this system can be obtained by holonomic deformation of the 4rd-order linear ordinary differential equation
\begin{align}\label{34}
\begin{split}
&\frac{d^4y}{dx^4}+a_1(x)\frac{d^3y}{dx^3}+a_2(x)\frac{d^2y}{dx^2}+a_3(x)\frac{dy}{dx}+a_4(x)y=0 \quad (a_i \in {\Bbb C}(x))
\end{split}
\end{align}
satisfying the Riemann scheme:
\begin{equation}\label{scheme}
\begin{pmatrix}
x=0 & x=1 & x=q_1 & x=q_2 & x=q_3 & x=\infty\\
0 & 0 & 0 & 0 & 0 & \alpha_0\\
\alpha_1 & \alpha_4 & 1 & 1 & 1 & \alpha_0\\
\alpha_1+\alpha_2 & \alpha_4+\alpha_5 & 2 & 2 & 2 & \alpha_0+\alpha_7\\
\alpha_1+\alpha_2+\alpha_3 & \alpha_4+\alpha_5+\alpha_6 & 4 & 4 & 4 & \alpha_0+\alpha_7\\
\end{pmatrix} \; ,
\end{equation}
where each $x=q_i$ is an apparent singular point.

The author conjectures that three apparent singular points $x=q_i$ satisfy $q_i \in {\Bbb C}(t)(q)$ or $q_i=q$ $(i=1,2,3)$.

\section{Holomorphy}
\begin{Theorem}
Let us consider a polynomial Hamiltonian system with Hamiltonian $I \in {\Bbb C}(t)[q,p]$. We assume that

$(A1)$ $deg(I)=10$ with respect to $q,p$.

$(A2)$ This system becomes again a polynomial Hamiltonian system in each coordinate $r_i \ (i=0,1,\ldots,6)${\rm : \rm}
\begin{align*}
&r_0:x_0=1/q, \ y_0=-(qp+\alpha_0)q,\\
&r_1:x_1=-(pq-\alpha_1)p, \ y_1=1/p,\\
&r_2:x_2=-(pq-(\alpha_1+\alpha_2))p, \ y_2=1/p,\\
&r_3:x_3=-(pq-(\alpha_1+\alpha_2+\alpha_3))p, \ y_3=1/p,\\
&r_4:x_4=-(p(q-1)-\alpha_4)p, \ y_4=1/p,\\
&r_5:x_5=-(p(q-1)-(\alpha_4+\alpha_5))p, \ y_5=1/p,\\
&r_6:x_6=-(p(q-1)-(\alpha_4+\alpha_5+\alpha_6))p, \ y_6=1/p
\end{align*}

$(A3)$ In addition to the assumption $(A2)$, the Hamiltonian system in the coordinate $r_0$ becomes again a polynomial Hamiltonian system in the coordinate $r_7${\rm : \rm}
\begin{align*}
&r_7:x_7=-(x_0y_0-\alpha_7)y_0, \ y_7=1/y_0.
\end{align*}
Then such a system coincides with the system
\begin{align}\label{E711}
\begin{split}
&\frac{dq}{dt}=\frac{\partial I}{\partial p}, \quad \frac{dp}{dt}=-\frac{\partial I}{\partial q},\\
I:=&(q-1)^3q^3p^4-(q-1)^2q^2\{(3\alpha_1+2\alpha_2+\alpha_3+3\alpha_4+2\alpha_5+\alpha_6)q-3\alpha_1-2\alpha_2-\alpha_3\}p^3\\
&+(q-1)q\{(6\alpha_0^2+6\alpha_0 \alpha_7+\alpha_7^2)q^2+(-6\alpha_0^2-3\alpha_1^2-4\alpha_1 \alpha_2-\alpha_2^2-2\alpha_1 \alpha_3-\alpha_2 \alpha_3+3\alpha_4^2+4\alpha_4 \alpha_5\\
&+\alpha_5^2+2\alpha_4 \alpha_6+\alpha_5 \alpha_6-6\alpha_0 \alpha_7-\alpha_7^2)q+3\alpha_1^2+\alpha_2^2+\alpha_2 \alpha_3+\alpha_1(4\alpha_2+2\alpha_3)\}p^2\\
&f p+\alpha_1(\alpha_1+\alpha_2)(\alpha_1+\alpha_2+\alpha_3)p+\alpha_0^2\{9\alpha_0^2+\alpha_7^2+\alpha_0(6\alpha_1+4\alpha_2+2\alpha_3+6\alpha_4+4\alpha_5+2\alpha_6+6\alpha_7)\}q^2\\
&+\alpha_0[15\alpha_0^3+3\alpha_1^2 \alpha_7+\alpha_2^2 \alpha_7+3\alpha_4^2 \alpha_7+\alpha_5^2 \alpha_7+\alpha_6 \alpha_7^2+\alpha_7^3\\
&+\alpha_0^2(25\alpha_1+14\alpha_2+3\alpha_3+15\alpha_4+10\alpha_5+5\alpha_6+12\alpha_7)+\alpha_5(\alpha_6 \alpha_7+2\alpha_7^2)\\
&+\alpha_2(3\alpha_4 \alpha_7+2\alpha_5 \alpha_7+\alpha_6 \alpha_7+2\alpha_7^2)+\alpha_4(4\alpha_5 \alpha_7+2\alpha_6 \alpha_7+3\alpha_7^2)\\
&+\alpha_1(3\alpha_2 \alpha_7+6\alpha_4 \alpha_7+4\alpha_5 \alpha_7+2\alpha_6 \alpha_7+4\alpha_7^2)+\alpha_0\{9\alpha_1^2+3\alpha_2^2+3\alpha_4^2+\alpha_5^2+3\alpha_6 \alpha_7+4\alpha_7^2)\\
&+\alpha_5(\alpha_6+6\alpha_7)+\alpha_2(\alpha_3+6\alpha_4+4\alpha_5+2\alpha_6+8\alpha_7)+\alpha_4(4\alpha_5+2\alpha_6+9\alpha_7)\\
&+\alpha_1(10\alpha_2+2\alpha_3+12\alpha_4+8\alpha_5+4\alpha_6+16\alpha_7)\}]q,
\end{split}
\end{align}
where $f$ is explicitly given by
\begin{align}
\begin{split}
f=&\alpha_0\{16\alpha_0^2+2\alpha_7^2+\alpha_0(9\alpha_1+6\alpha_2+3\alpha_3+9\alpha_4+6\alpha_5+3\alpha_6+12\alpha_7)\}q^3\\
&+[\alpha_0^2(16\alpha_1+8\alpha_2+9\alpha_4+6\alpha_5+3\alpha_6)\alpha_0^2+3\alpha_1^2 \alpha_7+\alpha_2^2 \alpha_7+3\alpha_4^2 \alpha_7+\alpha_5^2 \alpha_7+\alpha_6 \alpha_7^2+\alpha_7^3\\
&+\alpha_5(\alpha_6 \alpha_7+2\alpha_7^2)+\alpha_2(3\alpha_4 \alpha_7+2\alpha_5 \alpha_7+\alpha_6 \alpha_7+2\alpha_7^2)\\
&+\alpha_4(4\alpha_5 \alpha_7+2\alpha_6 \alpha_7+3\alpha_7^2)+\alpha_1(3\alpha_2 \alpha_7+6\alpha_4 \alpha_7+4\alpha_5 \alpha_7+2\alpha_6 \alpha_7+4\alpha_7^2)\\
&+\alpha_0\{6\alpha_1^2+2\alpha_2^2+6\alpha_4^2+2\alpha_5^2+3\alpha_6 \alpha_7+2\alpha_7^2+\alpha_5(2\alpha_6+6\alpha_7)+\alpha_2(6\alpha_4+4\alpha_5+2\alpha_6+8\alpha_7)\\
&+\alpha_4(8\alpha_5+4\alpha_6+9\alpha_7)+\alpha_1(6\alpha_2+12\alpha_4+8\alpha_5+4\alpha_6+16\alpha_7)\}]q^2\\
&+[-16\alpha_0^3-\alpha_1^3-\alpha_4^3+\alpha_0^2(-25\alpha_1-14\alpha_2-3\alpha_3-18\alpha_4-12\alpha_5-6\alpha_6-12\alpha_7)\\
&+\alpha_1^2(-2\alpha_2-\alpha_3-3\alpha_7)+\alpha_4^2(-2\alpha_5-\alpha_6-3\alpha_7)-\alpha_2^2 \alpha_7-\alpha_5^2 \alpha_7-\alpha_6 \alpha_7^2-\alpha_7^3\\
&+\alpha_0\{-6\alpha_1^2-2\alpha_2^2-6\alpha_4^2-2\alpha_5^2+\alpha_1(-6\alpha_2-12\alpha_4-8\alpha_5-4\alpha_6-16\alpha_7)+\alpha_4(-8\alpha_5-4\alpha_6-9\alpha_7)\\
&+\alpha_2(-6\alpha_4-4\alpha_5-2\alpha_6-8\alpha_7)+\alpha_5(-2\alpha_6-6\alpha_7)-3\alpha_6 \alpha_7-4\alpha_7^2\}\\
&+\alpha_1(-\alpha_2^2+\alpha_2(-\alpha_3-3\alpha_7)-6\alpha_4 \alpha_7-4\alpha_5 \alpha_7-2\alpha_6 \alpha_7-4\alpha_7^2)+\alpha_5(-\alpha_6 \alpha_7-2\alpha_7^2)\\
&+\alpha_4(-\alpha_5^2+\alpha_5(-\alpha_6-4\alpha_7)-2\alpha_6 \alpha_7-3\alpha_7^2)+\alpha_2(-3\alpha_4 \alpha_7-2\alpha_5 \alpha_7-\alpha_6 \alpha_7-2\alpha_7^2)]q.
\end{split}
\end{align}
Here, the constant parameters $\alpha_i$ satisfy the relation:
\begin{equation}
4\alpha_0+3\alpha_1+2\alpha_2+\alpha_3+3\alpha_4+2\alpha_5+\alpha_6+2\alpha_7=0.
\end{equation}
\end{Theorem}
Since each transformation $r_i$ is symplectic, the system \eqref{E711} is transformed into a Hamiltonian system, whose Hamiltonian may have poles. It is remarkable that the transformed system becomes again a polynomial system for any $i=0,1,\ldots,7$.

The holomorphy conditions $(A2),(A3)$ are new. Theorem 7.1 can be checked by a direct calculation.

\begin{Proposition}
The Hamiltonian $I$ is its first integral.
\end{Proposition}

\begin{Remark}
For the Hamiltonian system in each coordinate system $(x_i,y_i) \ (i=0,1,\ldots,7)$ given by $(A2)$ and $(A3)$ in Theorem 7.1, by eliminating $x_i$ or $y_i$, we obtain the second-order ordinary differential equation. However, its form is not normal \rm{(cf. \cite{Cosgrove1,Cosgrove2})}.
\end{Remark}

\section{Symmetry}

\begin{figure}[h]
\unitlength 0.1in
\begin{picture}( 48.9600, 13.1000)( 16.6000,-17.1000)
%
\special{pn 20}%
\special{pa 1810 1500}%
\special{pa 2910 400}%
\special{fp}%
\put(20.5000,-18.0000){\makebox(0,0)[lb]{$A_2$-lattice}}%
%
\special{pn 20}%
\special{ar 5216 1044 64 64  0.0000000 6.2831853}%
%
\special{pn 20}%
\special{ar 5226 1374 64 64  0.0000000 6.2831853}%
%
\special{pn 8}%
\special{pa 5216 1114}%
\special{pa 5216 1294}%
\special{fp}%
\put(43.2000,-18.8000){\makebox(0,0)[lb]{Dynkin diagram of type $E_7^{(1)}$}}%
%
\special{pn 20}%
\special{pa 1660 400}%
\special{pa 2780 1510}%
\special{fp}%
%
\special{pn 8}%
\special{pa 5294 1030}%
\special{pa 5574 1030}%
\special{fp}%
%
\special{pn 20}%
\special{ar 5644 1030 64 64  0.0000000 6.2831853}%
%
\special{pn 8}%
\special{pa 5714 1034}%
\special{pa 5994 1034}%
\special{fp}%
%
\special{pn 20}%
\special{ar 6064 1034 64 64  0.0000000 6.2831853}%
%
\special{pn 8}%
\special{pa 6144 1024}%
\special{pa 6424 1024}%
\special{fp}%
%
\special{pn 20}%
\special{ar 6494 1024 64 64  0.0000000 6.2831853}%
%
\special{pn 8}%
\special{pa 4014 1034}%
\special{pa 4294 1034}%
\special{fp}%
%
\special{pn 20}%
\special{ar 4364 1034 64 64  0.0000000 6.2831853}%
%
\special{pn 8}%
\special{pa 4434 1036}%
\special{pa 4714 1036}%
\special{fp}%
%
\special{pn 20}%
\special{ar 4784 1036 64 64  0.0000000 6.2831853}%
%
\special{pn 8}%
\special{pa 4864 1026}%
\special{pa 5144 1026}%
\special{fp}%
%
\special{pn 20}%
\special{ar 3934 1040 64 64  0.0000000 6.2831853}%
\end{picture}%
\label{fig:E7figure1}
\caption{$A_2$-lattice and Dynkin diagram of type $E_7^{(1)}$}
\end{figure}

\begin{Theorem}
The system \eqref{E711} admits extended affine Weyl group symmetry of type $E_7^{(1)}$ as the group of its B{\"a}cklund transformations whose generators $s_i, \ i=0,1,\ldots,7$ and $\pi$ are explicitly given as follows{\rm : \rm}with the notation $(*):=(q,p,t;\alpha_0,\alpha_1,\ldots,\alpha_7)$,
\begin{align*}
        s_{0}: (*) &\rightarrow \left(q+\frac{\alpha_0}{p},p,t;-\alpha_0,\alpha_1+\alpha_0,\alpha_2,\alpha_3,\alpha_4+\alpha_0,\alpha_5,\alpha_6,\alpha_7+\alpha_0 \right),\\
        s_{1}: (*) &\rightarrow \left(q,p-\frac{\alpha_1}{q},t;\alpha_0+\alpha_1,-\alpha_1,\alpha_2+\alpha_1,\alpha_3,\alpha_4,\alpha_5,\alpha_6,\alpha_7 \right), \\
        s_{2}: (*) &\rightarrow (q,p,t;\alpha_0,\alpha_1+\alpha_2,-\alpha_2,\alpha_3+\alpha_2,\alpha_4,\alpha_5,\alpha_6,\alpha_7), \\
        s_{3}: (*) &\rightarrow (q,p,t;\alpha_0,\alpha_1,\alpha_2+\alpha_3,-\alpha_3,\alpha_4,\alpha_5,\alpha_6,\alpha_7), \\
        s_{4}: (*) &\rightarrow \left(q,p-\frac{\alpha_4}{q-1},t;\alpha_0+\alpha_4,\alpha_1,\alpha_2,\alpha_3,-\alpha_4,\alpha_5+\alpha_4,\alpha_6,\alpha_7 \right),\\
        s_{5}: (*) &\rightarrow (q,p,t;\alpha_0,\alpha_1,\alpha_2,\alpha_3,\alpha_4+\alpha_5,-\alpha_5,\alpha_6+\alpha_5,\alpha_7),\\
        s_{6}: (*) &\rightarrow (q,p,t;\alpha_0,\alpha_1,\alpha_2,\alpha_3,\alpha_4,\alpha_5+\alpha_6,-\alpha_6,\alpha_7),\\
        s_{7}: (*) &\rightarrow (q,p,t;\alpha_0+\alpha_7,\alpha_1,\alpha_2,\alpha_3,\alpha_4,\alpha_5,\alpha_6,-\alpha_7),\\
        \pi: (*) &\rightarrow (1-q,-p,t+1;\alpha_0,\alpha_4,\alpha_5,\alpha_6,\alpha_1,\alpha_2,\alpha_3,\alpha_7),
\end{align*}
where $\pi$ is the Dynkin diagram automorphism of type $E_7^{(1)}$.
\end{Theorem}
The list should be read as
\begin{align}
\begin{split}
&s_{0}(\alpha_0)=-\alpha_0, \ s_{0}(\alpha_1)=\alpha_1+\alpha_0, \ s_{0}(\alpha_2)=\alpha_2, \ s_{0}(\alpha_3)=\alpha_3, \ s_{0}(\alpha_4)=\alpha_4+\alpha_0,\\
&s_{0}(\alpha_5)=\alpha_5, \ s_{0}(\alpha_6)=\alpha_6, \ s_{0}(\alpha_7)=\alpha_7+\alpha_0,\\
&s_{0}(q)=q+\frac{\alpha_0}{p}, \ s_{0}(p)=p, \ s_{0}(t)=t.
\end{split}
\end{align}
Theorem 8.1 can be checked by a direct calculation.

\section{Space of initial conditions}

\begin{Theorem}\label{3.1}
After a series of explicit blowing-ups at ten points including the infinitely near points of ${\Sigma_2}$ and successive blowing-down along the $(-1)$-curve ${D^{(0)}}' \cong {\Bbb P}^1$ and $D_{\infty}^{(1)} \cong {\Bbb P}^1$, we obtain the rational surface $\tilde{S}$ of the system \eqref{E711} and a birational morphism $\varphi:\tilde{S} \cdots \rightarrow {\Sigma_2}$. Its canonical divisor $K_{\tilde{S}}$ of $\tilde{S}$ is given by
\begin{align}
\begin{split}
K_{\tilde{S}}&=-D_{0}^{(1)}-D_{1}^{(1)}, \quad (D_{\nu}^{(1)})^2=-2, \ D_{\nu}^{(1)} \cong {\Bbb P}^1,
\end{split}
\end{align}
where the symbol ${D^{(0)}}'$ denotes the strict transform of $D^{(0)}$, $D_{\nu}^{(1)}$ denote the exceptional divisors and $-K_{\Sigma_2}=2D^{(0)}, \ D^{(0)} \cong {\Bbb P}^1, \ (D^{(0)})^2=2$.
\end{Theorem}

\begin{Theorem}
The space of initial conditions $S$ of the system \eqref{E711} is obtained by gluing nine copies of ${\Bbb C}^2$:
\begin{align}
\begin{split}
S&={\tilde{S}}-(-{K_{\tilde{S}}})_{red}\\
&={\Bbb C}^2 \cup \bigcup_{i=0}^{7} U_j,\\
&{\Bbb C}^2 \ni (q,p), \quad U_j \cong {\Bbb C}^2 \ni (x_j,y_j) \ (j=0,1,\ldots,7)
\end{split}
\end{align}
via the birational and symplectic transformations $r_j$ \rm{(see Theorem 8.1)}.
\end{Theorem}

{\bf Proof of Theorems 9.1 and 9.2.}

At first, we take the Hirzebruch surface ${\Sigma_2}$. By a direct calculation, we see that the system \eqref{E711} has three accessible singular points $a_{\nu}^{(0)} \in D^{(0)} \quad (\nu=0,1,\infty)$:
\begin{align}
\begin{split}
&a_{\nu}^{(0)}=\{(z_2,w_2)=(\nu,0)\} \in U_2 \cap D^{(0)} \ (\nu=0,1),\\
&a_{\infty}^{(0)}=\{(z_3,w_3)=(0,0)\} \in U_3 \cap D^{(0)}.
\end{split}
\end{align}
We perform blowing-ups in ${\Sigma_2}$ at $a_{\nu}^{(0)}$, and let $D_{\nu}^{(1)}$ be the exceptional curves of the blowing-ups at $a_{\nu}^{(0)}$ for $\nu=0,1,\infty$. We can take three coordinate systems $(u_{\nu},v_{\nu})$ around the points at infinity of the exceptional curves $D_{\nu}^{(1)} \quad (\nu=0,1,\infty)$, where
\begin{align}
\begin{split}
&(u_{\nu},v_{\nu})=\left(\frac{z_2-\nu}{w_2},w_2 \right) \ (\nu=0,1),\\
&(u_{\infty},v_{\infty})=\left(\frac{z_3}{w_3},w_3 \right).
\end{split}
\end{align}

\begin{figure}
\unitlength 0.1in
\begin{picture}( 28.2000, 84.9800)( 16.4000,-86.1600)
\put(44.5000,-4.1400){\makebox(0,0)[lb]{$D^{(0)} \cong {\Bbb P}^1, \quad (D^{(0)})^2=2$}}%
%
\special{pn 8}%
\special{pa 2110 332}%
\special{pa 4140 332}%
\special{fp}%
%
\special{pn 20}%
\special{sh 0.600}%
\special{ar 2420 332 20 16  0.0000000 6.2831853}%
%
\special{pn 20}%
\special{sh 0.600}%
\special{ar 3080 332 20 16  0.0000000 6.2831853}%
%
\special{pn 20}%
\special{sh 0.600}%
\special{ar 3760 332 20 16  0.0000000 6.2831853}%
\put(22.9000,-2.8800){\makebox(0,0)[lb]{$a_0^{(0)}$}}%
\put(29.6000,-2.8800){\makebox(0,0)[lb]{$a_1^{(0)}$}}%
\put(36.3000,-2.8800){\makebox(0,0)[lb]{$a_{\infty}^{(0)}$}}%
%
\special{pn 20}%
\special{pa 3086 918}%
\special{pa 3086 502}%
\special{fp}%
\special{sh 1}%
\special{pa 3086 502}%
\special{pa 3066 570}%
\special{pa 3086 556}%
\special{pa 3106 570}%
\special{pa 3086 502}%
\special{fp}%
\put(31.7500,-8.4100){\makebox(0,0)[lb]{Blow up at $a_0^{(0)},a_1^{(0)}$ and $a_{\infty}^{(0)}$}}%
%
\special{pn 8}%
\special{pa 2120 2094}%
\special{pa 4150 2094}%
\special{fp}%
%
\special{pn 8}%
\special{pa 2430 2338}%
\special{pa 2430 978}%
\special{dt 0.045}%
%
\special{pn 20}%
\special{sh 0.600}%
\special{ar 2430 1754 20 16  0.0000000 6.2831853}%
%
\special{pn 20}%
\special{sh 0.600}%
\special{ar 2430 1308 20 18  0.0000000 6.2831853}%
\put(21.1000,-14.0400){\makebox(0,0)[lb]{$a_1^{(1)}$}}%
\put(21.1000,-18.0800){\makebox(0,0)[lb]{$a_2^{(1)}$}}%
%
\special{pn 8}%
\special{pa 3080 2328}%
\special{pa 3080 968}%
\special{dt 0.045}%
%
\special{pn 20}%
\special{sh 0.600}%
\special{ar 3080 1742 20 18  0.0000000 6.2831853}%
%
\special{pn 20}%
\special{sh 0.600}%
\special{ar 3080 1298 20 18  0.0000000 6.2831853}%
\put(27.6000,-13.9300){\makebox(0,0)[lb]{$a_4^{(1)}$}}%
\put(27.6000,-17.9600){\makebox(0,0)[lb]{$a_5^{(1)}$}}%
%
\special{pn 8}%
\special{pa 3780 2338}%
\special{pa 3780 978}%
\special{dt 0.045}%
%
\special{pn 20}%
\special{sh 0.600}%
\special{ar 3780 1308 20 18  0.0000000 6.2831853}%
\put(34.6000,-14.0400){\makebox(0,0)[lb]{$a_7^{(1)}$}}%
%
\special{pn 20}%
\special{pa 3086 2826}%
\special{pa 3086 2412}%
\special{fp}%
\special{sh 1}%
\special{pa 3086 2412}%
\special{pa 3066 2480}%
\special{pa 3086 2466}%
\special{pa 3106 2480}%
\special{pa 3086 2412}%
\special{fp}%
\put(31.7500,-27.5300){\makebox(0,0)[lb]{Blow up at $a_j^{(1)}$}}%
%
\special{pn 8}%
\special{pa 2110 4188}%
\special{pa 4140 4188}%
\special{fp}%
%
\special{pn 8}%
\special{pa 2420 4430}%
\special{pa 2420 3072}%
\special{fp}%
%
\special{pn 8}%
\special{pa 2320 3954}%
\special{pa 2670 3708}%
\special{dt 0.045}%
%
\special{pn 8}%
\special{pa 2330 3708}%
\special{pa 2630 3422}%
\special{dt 0.045}%
%
\special{pn 8}%
\special{pa 3070 4442}%
\special{pa 3070 3082}%
\special{fp}%
%
\special{pn 8}%
\special{pa 2970 3960}%
\special{pa 3320 3720}%
\special{dt 0.045}%
%
\special{pn 8}%
\special{pa 2970 3672}%
\special{pa 3270 3384}%
\special{dt 0.045}%
%
\special{pn 8}%
\special{pa 3760 4454}%
\special{pa 3760 3092}%
\special{fp}%
%
\special{pn 8}%
\special{pa 3660 3974}%
\special{pa 4010 3728}%
\special{dt 0.045}%
\put(44.6000,-21.7000){\makebox(0,0)[lb]{$(D^{(0)})^2=-1$}}%
\put(21.5000,-30.7200){\makebox(0,0)[lb]{$(D_{0}^{(1)})^2=-4$}}%
\put(28.1000,-46.3400){\makebox(0,0)[lb]{$(D_{1}^{(1)})^2=-4$}}%
\put(35.1000,-30.8200){\makebox(0,0)[lb]{$(D_{\infty}^{(1)})^2=-2$}}%
%
\special{pn 20}%
\special{pa 3090 4716}%
\special{pa 3090 5090}%
\special{fp}%
\special{sh 1}%
\special{pa 3090 5090}%
\special{pa 3110 5022}%
\special{pa 3090 5036}%
\special{pa 3070 5022}%
\special{pa 3090 5090}%
\special{fp}%
\put(32.0000,-49.8400){\makebox(0,0)[lb]{Blow down along the $(-1)$-curve $D^{(0)}$ to a nonsingular point}}%
%
\special{pn 20}%
\special{pa 2210 5354}%
\special{pa 3310 6534}%
\special{fp}%
\special{pa 3080 5346}%
\special{pa 3080 6534}%
\special{fp}%
%
\special{pn 8}%
\special{pa 2840 6546}%
\special{pa 3940 5378}%
\special{dt 0.045}%
%
\special{pn 8}%
\special{pa 2340 5682}%
\special{pa 2620 5238}%
\special{dt 0.045}%
%
\special{pn 8}%
\special{pa 2590 5910}%
\special{pa 2780 5588}%
\special{dt 0.045}%
%
\special{pn 8}%
\special{pa 2970 5610}%
\special{pa 3260 5396}%
\special{dt 0.045}%
%
\special{pn 8}%
\special{pa 2980 5896}%
\special{pa 3280 5682}%
\special{dt 0.045}%
%
\special{pn 8}%
\special{pa 3510 5642}%
\special{pa 3760 5982}%
\special{dt 0.045}%
\put(36.2000,-85.2000){\makebox(0,0)[lb]{Each bold line denotes $(-2)$-curve}}%
%
\special{pn 8}%
\special{pa 2330 3404}%
\special{pa 2630 3116}%
\special{dt 0.045}%
%
\special{pn 8}%
\special{pa 2970 3404}%
\special{pa 3270 3116}%
\special{dt 0.045}%
%
\special{pn 8}%
\special{pa 2480 5790}%
\special{pa 2670 5472}%
\special{dt 0.045}%
%
\special{pn 8}%
\special{pa 2970 5778}%
\special{pa 3260 5568}%
\special{dt 0.045}%
%
\special{pn 20}%
\special{sh 0.600}%
\special{ar 2430 1998 20 16  0.0000000 6.2831853}%
%
\special{pn 20}%
\special{sh 0.600}%
\special{ar 3070 1998 20 16  0.0000000 6.2831853}%
\put(21.1000,-20.4400){\makebox(0,0)[lb]{$a_3^{(1)}$}}%
\put(27.6000,-20.3200){\makebox(0,0)[lb]{$a_6^{(1)}$}}%
%
\special{pn 20}%
\special{pa 3086 6822}%
\special{pa 3086 7194}%
\special{fp}%
\special{sh 1}%
\special{pa 3086 7194}%
\special{pa 3106 7128}%
\special{pa 3086 7142}%
\special{pa 3066 7128}%
\special{pa 3086 7194}%
\special{fp}%
\put(31.9500,-70.8800){\makebox(0,0)[lb]{Blow down along the $(-1)$-curve $D_{\infty}^{(1)}$ to a nonsingular point}}%
%
\special{pn 20}%
\special{pa 2190 7436}%
\special{pa 3290 8616}%
\special{fp}%
\special{pa 3060 7424}%
\special{pa 3060 8616}%
\special{fp}%
%
\special{pn 8}%
\special{pa 2320 7764}%
\special{pa 2600 7320}%
\special{dt 0.045}%
%
\special{pn 8}%
\special{pa 2570 7990}%
\special{pa 2760 7670}%
\special{dt 0.045}%
%
\special{pn 8}%
\special{pa 2950 7692}%
\special{pa 3240 7478}%
\special{dt 0.045}%
%
\special{pn 8}%
\special{pa 2960 7980}%
\special{pa 3260 7764}%
\special{dt 0.045}%
%
\special{pn 8}%
\special{pa 2460 7872}%
\special{pa 2650 7554}%
\special{dt 0.045}%
%
\special{pn 8}%
\special{pa 2950 7862}%
\special{pa 3240 7648}%
\special{dt 0.045}%
\put(37.2000,-53.5000){\makebox(0,0)[lb]{$(D_{\infty}^{(1)})^2=-1$}}%
%
\special{pn 8}%
\special{pa 2510 8376}%
\special{pa 3480 8376}%
\special{dt 0.045}%
\put(16.4000,-53.0000){\makebox(0,0)[lb]{$(D_{0}^{(1)})^2=-3$}}%
\put(27.2000,-53.1000){\makebox(0,0)[lb]{$(D_{1}^{(1)})^2=-3$}}%
\end{picture}%
\label{fig:E7figure2}
\caption{Resolution of accessible singular points}
\end{figure}

Note that $\{(u_{\nu},v_{\nu})|v_{\nu}=0\} \subset D_{\nu}^{(1)}$ for $\nu=0,1,\infty$. By a direct calculation, we see that the system \eqref{11} has seven accessible singular points $a_{\nu}^{(1)}$ for $\nu=1,2,3,4,5,6,7$ in $D_{\nu}^{(1)} \cong {\Bbb P}^1 \ (\nu=0,1,\infty)$.
\begin{align}
\begin{split}
&a_{1}^{(1)}=\{(u_{0},v_{0})=(\alpha_1,0)\} \in D_{0}^{(1)}, \quad a_{2}^{(1)}=\{(u_{0},v_{0})=(\alpha_1+\alpha_2,0)\} \in D_{0}^{(1)},\\
&a_{3}^{(1)}=\{(u_{0},v_{0})=(\alpha_1+\alpha_2+\alpha_3,0)\} \in D_{0}^{(1)}, \quad a_{4}^{(1)}=\{(u_{1},v_{1})=(\alpha_4,0)\} \in D_{1}^{(1)},\\
&a_{5}^{(1)}=\{(u_{1},v_{1})=(\alpha_4+\alpha_5,0)\} \in D_{1}^{(1)}, \quad a_{6}^{(1)}=\{(u_{1},v_{1})=(\alpha_4+\alpha_5+\alpha_6,0)\} \in D_{1}^{(1)},\\
&a_{7}^{(1)}=\{(u_{\infty},v_{\infty})=(\alpha_7,0)\} \in D_{\infty}^{(1)}.
\end{split}
\end{align}
Let us perform blowing-ups at $a_{j}^{(1)}$, and denote $D_{j}^{(2)}$ for the exceptional curves, respectively. We take seven coordinate systems $(W_j,V_j)$ around the points at infinity of $D_{j}^{(2)}$ for $j=1,2,3,4,5,6,7$, where
\begin{align}
\begin{split}
&(W_{1},V_{1})=\left(\frac{u_0-\alpha_1}{v_0},v_0 \right),\\
&(W_{2},V_{2})=\left(\frac{u_0-(\alpha_1+\alpha_2)}{v_0},v_0 \right),\\
&(W_{3},V_{3})=\left(\frac{u_0-(\alpha_1+\alpha_2+\alpha_3)}{v_0},v_0 \right),\\
&(W_{4},V_{4})=\left(\frac{u_1-\alpha_4}{v_1},v_1 \right),\\
&(W_{5},V_{5})=\left(\frac{u_1-(\alpha_4+\alpha_5)}{v_1},v_1 \right),\\
&(W_{6},V_{6})=\left(\frac{u_1-(\alpha_4+\alpha_5+\alpha_6)}{v_1},v_1 \right),\\
&(W_{7},V_{7})=\left(\frac{u_{\infty}-\alpha_7}{v_{\infty}},v_{\infty} \right).
\end{split}
\end{align}
For the strict transform of $D^{(0)}$, $D_{\nu}^{(1)}$ and $D_{j}^{(2)}$ by the blowing-ups, we also denote by same symbol, respectively.  Here, the self-intersection number of $D^{(0)}, D_{\nu}^{(1)}$ is given by
\begin{equation}
(D^{(0)})^2=-1. \quad (D_{0}^{(1)})^2=(D_{1}^{(1)})^2=-4, \quad (D_{\infty}^{(1)})^2=-2.
\end{equation}
In order to obtain a minimal compactification of the space of initial conditions, we must blow down along the $(-1)$-curve $D^{(0)} \cong {\Bbb P}^1$ to a nonsingular point. For the strict transform of $D_{\nu}^{(1)}$ and $D_{j}^{(2)}$ by the blowing-down, we also denote by same symbol, respectively. Here, the self-intersection number of $D_{\nu}^{(1)}$ is given by
\begin{equation}
(D_{0}^{(1)})^2=(D_{1}^{(1)})^2=-3, \quad (D_{\infty}^{(1)})^2=-1.
\end{equation}
We must blow down again along the $(-1)$-curve $D_{\infty}^{(1)} \cong {\Bbb P}^1$ to a nonsingular point. For the strict transform of $D_{\nu}^{(1)}$ and $D_{j}^{(2)}$ by the blowing-down, we also denote by same symbol, respectively. Here, the self-intersection number of $D_{\nu}^{(1)}$ is given by
\begin{equation}
(D_{0}^{(1)})^2=(D_{1}^{(1)})^2=-2.
\end{equation}

Let ${\tilde S} \cdots \rightarrow {\Sigma_2}$ be the composition of above ten times blowing-ups and two  times blowing-downs. Then, we see that the canonical divisor class $K_{{\tilde S}}$ of ${\tilde S}$ is given by
\begin{equation}
K_{{\tilde S}}:=-D_{0}^{(1)}-D_{1}^{(1)},
\end{equation}
where the self-intersection number of $D_{\nu}^{(1)} \cong {\Bbb P}^1$ is given by
\begin{equation}
(D_{\nu}^{(1)})^2=-2,
\end{equation}
and
\begin{equation}
D_{0}^{(1)} \cap D_{1}^{(1)} \not=\varnothing, \quad (D_{0}^{(1)},D_{1}^{(1)})=1.
\end{equation}
The configuration of the divisor $(-K_{{\tilde S}})_{red}$ on $\tilde S$ is of type $A_2$. And we see that ${\tilde S}-(-K_{{\tilde S}})_{red}$ is covered by nine Zariski open sets
\begin{align}
\begin{split}
& \rm{Spec} \ {\Bbb C}[W_{j},V_{j}] \quad (j=1,2,3,4,5,6,7),\\
& \rm{Spec} \ {\Bbb C}[z_0,w_0],\\
& \rm{Spec} \ {\Bbb C}[z_1,w_1].
\end{split}
\end{align}
The relations between $(W_{j},V_{j})$ and $(x_j,y_j)$ are given by
\begin{equation}
(-W_{j},V_{j})=(x_j,y_j) \quad (j=1,2,3,4,5,6,7).
\end{equation}
We see that the pole divisor of the symplectic 2-form $dp \wedge dq$ coincides with $(-K_{{\tilde S}})_{red}$. Thus, we have completed the proof of Theorems 9.1 and 9.2. \qed

\section{Main results of the system with $W(E_8^{(1)})$-symmetry}

By using the key property, we try to make a second-order polynomial Hamiltonian system with symmetry of the affine Weyl group of type $E_8^{(1)}$.

At first, we make a representation of affine Weyl group of type  $E_8^{(1)}$. Next, we make holomorphy conditions $r_i \ (i=0,1,\ldots,8)$ associated with it.

\begin{Problem}
Can we make a polynomial Hamiltonian system with Hamiltonian $I \in {\Bbb C}(t)[q,p]$ satisfying the following condition $(A)$?:

$(A)$:This system becomes again a polynomial Hamiltonian system in each coordinate $r_i \ (i=0,1,\ldots,8)$
\end{Problem}
Before we solve this problem, we construct the space of initial conditions characterized by holomorphy conditions $r_i \ (i=0,1,\ldots,8)$. After a series of explicit blowing-ups at eleven points including the infinitely near points of the Hirzebruch surface ${\Sigma_2}$ (see Figure 2) and three times blowing-downs along the $(-1)$-curve to a nonsingular point (see Figure 7), respectively, we obtain the  smooth rational surface $\tilde{S}$ and a birational morphism 
$$
\varphi:\tilde{S}=S_{14} \leftarrow S_{13} \leftarrow S_{12} \leftarrow S_{11} \rightarrow \cdots \rightarrow S_1 \rightarrow {\Sigma_2}.
$$
Here, $-K_{\Sigma_2}=2H, \ H \cong {\Bbb P}^1, \ (H)^2=2$. In order to obtain a minimal compactification of the space of initial conditions, we must blow down along the $(-1)$-curves. Its canonical divisor $K_{\tilde{S}}$ of $\tilde{S}$ is given by
\begin{align}
\begin{split}
K_{\tilde{S}}&=-E, \quad E \cong {\Bbb P}^1, \quad (E)^2=-3.
\end{split}
\end{align}
In the case of Painlev\'e equations, each component of the anti-canonical divisor $-K_{\tilde{S}}$ is $(-2)$-curve. However, in this case, it is $(-3)$-curve. In this vein, we take a poor view of holomorphy conditions $r_i \ (i=0,1,\ldots,8)$.

However, we can obtain a 8-parameter family of polynomial Hamiltonian systems determined by $r_i \ (i=0,1,\ldots,8)$. Suprisingly, this system admits the affine Weyl group symmetry of type $E_8^{(1)}$ as the group of its B{\"a}cklund transformations (see Figure 6). This system is the first example which gave second-order polynomial Hamiltonian systems with $W(E_8^{(1)})$-symmetry.

By eliminating $p$ or $q$, we obtain the second-order ordinary differential equation. However, its form is not normal \rm{(cf. \cite{Cosgrove1,Cosgrove2})}.

The space of initial conditions $S$ is obtained by gluing ten copies of ${\Bbb C}^2$
\begin{align}
\begin{split}
S&={\tilde{S}}-(-{K_{\tilde{S}}})_{red}\\
&={\Bbb C}^2 \cup \bigcup_{i=0}^{8} {\Bbb C}^2
\end{split}
\end{align}
via the birational and symplectic transformations $r_j$ (see Theorem 11.1).

\begin{figure}
\unitlength 0.1in
\begin{picture}( 30.8300,  6.8000)( 24.1000,-16.4000)
%
\special{pn 20}%
\special{ar 3326 1044 64 64  0.0000000 6.2831853}%
%
\special{pn 20}%
\special{ar 3336 1374 64 64  0.0000000 6.2831853}%
%
\special{pn 8}%
\special{pa 3326 1114}%
\special{pa 3326 1294}%
\special{fp}%
\put(24.5000,-18.1000){\makebox(0,0)[lb]{Dynkin diagram of type $E_8^{(1)}$}}%
%
\special{pn 8}%
\special{pa 3404 1030}%
\special{pa 3684 1030}%
\special{fp}%
%
\special{pn 20}%
\special{ar 3754 1030 64 64  0.0000000 6.2831853}%
%
\special{pn 8}%
\special{pa 3824 1034}%
\special{pa 4104 1034}%
\special{fp}%
%
\special{pn 20}%
\special{ar 4174 1034 64 64  0.0000000 6.2831853}%
%
\special{pn 8}%
\special{pa 4254 1024}%
\special{pa 4534 1024}%
\special{fp}%
%
\special{pn 20}%
\special{ar 4604 1024 64 64  0.0000000 6.2831853}%
%
\special{pn 20}%
\special{ar 2474 1034 64 64  0.0000000 6.2831853}%
%
\special{pn 8}%
\special{pa 2544 1036}%
\special{pa 2824 1036}%
\special{fp}%
%
\special{pn 20}%
\special{ar 2894 1036 64 64  0.0000000 6.2831853}%
%
\special{pn 8}%
\special{pa 2974 1026}%
\special{pa 3254 1026}%
\special{fp}%
%
\special{pn 8}%
\special{pa 4660 1024}%
\special{pa 4940 1024}%
\special{fp}%
%
\special{pn 20}%
\special{ar 5010 1024 64 64  0.0000000 6.2831853}%
%
\special{pn 8}%
\special{pa 5080 1024}%
\special{pa 5360 1024}%
\special{fp}%
%
\special{pn 20}%
\special{ar 5430 1024 64 64  0.0000000 6.2831853}%
\end{picture}%
\label{fig:E8figure1}
\caption{}
\end{figure}

The author believes that this system can be obtained by holonomic deformation of the 6rd-order linear ordinary differential equation
\begin{align}\label{E834}
\begin{split}
&\frac{d^6y}{dx^6}+a_1(x)\frac{d^5y}{dx^5}+a_2(x)\frac{d^4y}{dx^4}+a_3(x)\frac{d^3y}{dx^3}+a_4(x)\frac{d^2y}{dx^2}+a_5(x)\frac{dy}{dx}+a_6(x)y=0 \quad (a_i \in {\Bbb C}(x))
\end{split}
\end{align}
satisfying the Riemann scheme:
\begin{equation}\label{E8scheme}
\begin{pmatrix}
x=0 & x=1 & x=q_i & x=\infty\\
0 & 0 & 0 & \alpha_0\\
\alpha_1 & 0 & 1 & \alpha_0\\
\alpha_1+\alpha_2 & \alpha_6 & 2 & \alpha_0\\
\alpha_1+\alpha_2+\alpha_3 & \alpha_6 & 3 & \alpha_0+\alpha_8\\
\alpha_1+\alpha_2+\alpha_3+\alpha_4 & \alpha_6+\alpha_7 & 4 & \alpha_0+\alpha_8\\
\alpha_1+\alpha_2+\alpha_3+\alpha_4+\alpha_5 & \alpha_6+\alpha_7 & 6 & \alpha_0+\alpha_8\\
\end{pmatrix} \; ,
\end{equation}
where each $x=q_i \ (i=1,2,\ldots,10)$ is an apparent singular point.

The author conjectures that ten apparent singular points $x=q_i$ satisfy $q_i \in {\Bbb C}(t)(q)$ or $q_i=q$.

\section{Holomorphy}
\begin{Theorem}
Let us consider a polynomial Hamiltonian system with Hamiltonian $I \in {\Bbb C}(t)[q,p]$. We assume that

$(C1)$ $deg(I)=15$ with respect to $q,p$.

$(C2)$ This system becomes again a polynomial Hamiltonian system in each coordinate $r_i \ (i=0,1,\ldots,7)${\rm : \rm}
\begin{align*}
&r_0:x_0=1/q, \ y_0=-(qp+\alpha_0)q,\\
&r_1:x_1=-(pq-\alpha_1)p, \ y_1=1/p,\\
&r_2:x_2=-(pq-(\alpha_1+\alpha_2))p, \ y_2=1/p,\\
&r_3:x_3=-(pq-(\alpha_1+\alpha_2+\alpha_3))p, \ y_3=1/p,\\
&r_4:x_4=-(pq-(\alpha_1+\alpha_2+\alpha_3+\alpha_4))p, \ y_4=1/p,\\
&r_5:x_5=-(pq-(\alpha_1+\alpha_2+\alpha_3+\alpha_4+\alpha_5))p, \ y_5=1/p,\\
&r_6:x_6=-(p(q-1)-\alpha_6)p, \ y_6=1/p,\\
&r_7:x_7=-(p(q-1)-(\alpha_6+\alpha_7))p, \ y_7=1/p.
\end{align*}

$(C3)$ In addition to the assumption $(C2)$, the Hamiltonian system in the coordinate $r_0$ becomes again a polynomial Hamiltonian system in the coordinate $r_7${\rm : \rm}
\begin{align*}
&r_8:x_8=-(x_0y_0-\alpha_8)y_0, \ y_8=1/y_0.
\end{align*}
Then such a system coincides with the system
\begin{equation}\label{E811}
\frac{dq}{dt}=\frac{\partial I}{\partial p}, \quad \frac{dp}{dt}=-\frac{\partial I}{\partial q},
\end{equation}
$
I:=p^6 (-1 + q)^4 q^5 + 
 p^5 (-1 + q)^3 q^4 (-6 \alpha_0 - 4 \alpha_6 - 2 \alpha_7 - 
    3 \alpha_8 + 3 q (2 \alpha_0 + \alpha_8)) + 
 p^4 (-1 + q)^2 q^3 (-24 \alpha_0 \alpha_1 - 10 \alpha_1^2 - 
    18 \alpha_0 \alpha_2 - 15 \alpha_1 \alpha_2 - 6 \alpha_2^2 - 
    12 \alpha_0 \alpha_3 - 10 \alpha_1 \alpha_3 - 
    8 \alpha_2 \alpha_3 - 3 \alpha_3^2 - 6 \alpha_0 \alpha_4 - 
    5 \alpha_1 \alpha_4 - 4 \alpha_2 \alpha_4 - 
    3 \alpha_3 \alpha_4 - \alpha_4^2 - 16 \alpha_1 \alpha_6 - 
    12 \alpha_2 \alpha_6 - 8 \alpha_3 \alpha_6 - 
    4 \alpha_4 \alpha_6 - 8 \alpha_1 \alpha_7 - 
    6 \alpha_2 \alpha_7 - 4 \alpha_3 \alpha_7 - 
    2 \alpha_4 \alpha_7 - 12 \alpha_1 \alpha_8 - 
    9 \alpha_2 \alpha_8 - 6 \alpha_3 \alpha_8 - 
    3 \alpha_4 \alpha_8 + 
    q (-15 \alpha_0^2 + 24 \alpha_0 \alpha_1 + 10 \alpha_1^2 + 
       18 \alpha_0 \alpha_2 + 15 \alpha_1 \alpha_2 + 
       6 \alpha_2^2 + 12 \alpha_0 \alpha_3 + 
       10 \alpha_1 \alpha_3 + 8 \alpha_2 \alpha_3 + 
       3 \alpha_3^2 + 6 \alpha_0 \alpha_4 + 
       5 \alpha_1 \alpha_4 + 4 \alpha_2 \alpha_4 + 
       3 \alpha_3 \alpha_4 + \alpha_4^2 + 16 \alpha_1 \alpha_6 + 
       12 \alpha_2 \alpha_6 + 8 \alpha_3 \alpha_6 + 
       4 \alpha_4 \alpha_6 + 6 \alpha_6^2 + 
       8 \alpha_1 \alpha_7 + 6 \alpha_2 \alpha_7 + 
       4 \alpha_3 \alpha_7 + 2 \alpha_4 \alpha_7 + 
       6 \alpha_6 \alpha_7 + \alpha_7^2 - 15 \alpha_0 \alpha_8 + 
       12 \alpha_1 \alpha_8 + 9 \alpha_2 \alpha_8 + 
       6 \alpha_3 \alpha_8 + 3 \alpha_4 \alpha_8 - 
       3 \alpha_8^2) + 
    3 q^2 (5 \alpha_0^2 + 5 \alpha_0 \alpha_8 + \alpha_8^2)) + 
 q \alpha_0 (\alpha_0 + \alpha_8) (-44 \alpha_0^4 - 
    120 \alpha_0^3 \alpha_1 - 122 \alpha_0^2 \alpha_1^2 - 
    64 \alpha_0 \alpha_1^3 - 15 \alpha_1^4 - 
    90 \alpha_0^3 \alpha_2 - 183 \alpha_0^2 \alpha_1 \alpha_2 - 
    144 \alpha_0 \alpha_1^2 \alpha_2 - 45 \alpha_1^3 \alpha_2 - 
    66 \alpha_0^2 \alpha_2^2 - 
    102 \alpha_0 \alpha_1 \alpha_2^2 - 
    48 \alpha_1^2 \alpha_2^2 - 22 \alpha_0 \alpha_2^3 - 
    21 \alpha_1 \alpha_2^3 - 3 \alpha_2^4 - 
    60 \alpha_0^3 \alpha_3 - 122 \alpha_0^2 \alpha_1 \alpha_3 - 
    96 \alpha_0 \alpha_1^2 \alpha_3 - 30 \alpha_1^3 \alpha_3 - 
    88 \alpha_0^2 \alpha_2 \alpha_3 - 
    136 \alpha_0 \alpha_1 \alpha_2 \alpha_3 - 
    64 \alpha_1^2 \alpha_2 \alpha_3 - 
    44 \alpha_0 \alpha_2^2 \alpha_3 - 
    42 \alpha_1 \alpha_2^2 \alpha_3 - 8 \alpha_2^3 \alpha_3 - 
    27 \alpha_0^2 \alpha_3^2 - 40 \alpha_0 \alpha_1 \alpha_3^2 - 
    19 \alpha_1^2 \alpha_3^2 - 26 \alpha_0 \alpha_2 \alpha_3^2 - 
    25 \alpha_1 \alpha_2 \alpha_3^2 - 7 \alpha_2^2 \alpha_3^2 - 
    4 \alpha_0 \alpha_3^3 - 4 \alpha_1 \alpha_3^3 - 
    2 \alpha_2 \alpha_3^3 - 30 \alpha_0^3 \alpha_4 - 
    61 \alpha_0^2 \alpha_1 \alpha_4 - 
    48 \alpha_0 \alpha_1^2 \alpha_4 - 15 \alpha_1^3 \alpha_4 - 
    44 \alpha_0^2 \alpha_2 \alpha_4 - 
    68 \alpha_0 \alpha_1 \alpha_2 \alpha_4 - 
    32 \alpha_1^2 \alpha_2 \alpha_4 - 
    22 \alpha_0 \alpha_2^2 \alpha_4 - 
    21 \alpha_1 \alpha_2^2 \alpha_4 - 4 \alpha_2^3 \alpha_4 - 
    27 \alpha_0^2 \alpha_3 \alpha_4 - 
    40 \alpha_0 \alpha_1 \alpha_3 \alpha_4 - 
    19 \alpha_1^2 \alpha_3 \alpha_4 - 
    26 \alpha_0 \alpha_2 \alpha_3 \alpha_4 - 
    25 \alpha_1 \alpha_2 \alpha_3 \alpha_4 - 
    7 \alpha_2^2 \alpha_3 \alpha_4 - 
    6 \alpha_0 \alpha_3^2 \alpha_4 - 
    6 \alpha_1 \alpha_3^2 \alpha_4 - 
    3 \alpha_2 \alpha_3^2 \alpha_4 - 5 \alpha_0^2 \alpha_4^2 - 
    6 \alpha_0 \alpha_1 \alpha_4^2 - 3 \alpha_1^2 \alpha_4^2 - 
    4 \alpha_0 \alpha_2 \alpha_4^2 - 
    4 \alpha_1 \alpha_2 \alpha_4^2 - \alpha_2^2 \alpha_4^2 - 
    2 \alpha_0 \alpha_3 \alpha_4^2 - 
    2 \alpha_1 \alpha_3 \alpha_4^2 - \alpha_2 \alpha_3 \
\alpha_4^2 - 128 \alpha_0^3 \alpha_6 - 
    272 \alpha_0^2 \alpha_1 \alpha_6 - 
    200 \alpha_0 \alpha_1^2 \alpha_6 - 56 \alpha_1^3 \alpha_6 - 
    204 \alpha_0^2 \alpha_2 \alpha_6 - 
    300 \alpha_0 \alpha_1 \alpha_2 \alpha_6 - 
    126 \alpha_1^2 \alpha_2 \alpha_6 - 
    108 \alpha_0 \alpha_2^2 \alpha_6 - 
    90 \alpha_1 \alpha_2^2 \alpha_6 - 20 \alpha_2^3 \alpha_6 - 
    136 \alpha_0^2 \alpha_3 \alpha_6 - 
    200 \alpha_0 \alpha_1 \alpha_3 \alpha_6 - 
    84 \alpha_1^2 \alpha_3 \alpha_6 - 
    144 \alpha_0 \alpha_2 \alpha_3 \alpha_6 - 
    120 \alpha_1 \alpha_2 \alpha_3 \alpha_6 - 
    40 \alpha_2^2 \alpha_3 \alpha_6 - 
    44 \alpha_0 \alpha_3^2 \alpha_6 - 
    36 \alpha_1 \alpha_3^2 \alpha_6 - 
    24 \alpha_2 \alpha_3^2 \alpha_6 - 4 \alpha_3^3 \alpha_6 - 
    68 \alpha_0^2 \alpha_4 \alpha_6 - 
    100 \alpha_0 \alpha_1 \alpha_4 \alpha_6 - 
    42 \alpha_1^2 \alpha_4 \alpha_6 - 
    72 \alpha_0 \alpha_2 \alpha_4 \alpha_6 - 
    60 \alpha_1 \alpha_2 \alpha_4 \alpha_6 - 
    20 \alpha_2^2 \alpha_4 \alpha_6 - 
    44 \alpha_0 \alpha_3 \alpha_4 \alpha_6 - 
    36 \alpha_1 \alpha_3 \alpha_4 \alpha_6 - 
    24 \alpha_2 \alpha_3 \alpha_4 \alpha_6 - 
    6 \alpha_3^2 \alpha_4 \alpha_6 - 
    8 \alpha_0 \alpha_4^2 \alpha_6 - 
    6 \alpha_1 \alpha_4^2 \alpha_6 - 
    4 \alpha_2 \alpha_4^2 \alpha_6 - 
    2 \alpha_3 \alpha_4^2 \alpha_6 - 137 \alpha_0^2 \alpha_6^2 - 
    200 \alpha_0 \alpha_1 \alpha_6^2 - 
    78 \alpha_1^2 \alpha_6^2 - 
    150 \alpha_0 \alpha_2 \alpha_6^2 - 
    117 \alpha_1 \alpha_2 \alpha_6^2 - 
    42 \alpha_2^2 \alpha_6^2 - 
    100 \alpha_0 \alpha_3 \alpha_6^2 - 
    78 \alpha_1 \alpha_3 \alpha_6^2 - 
    56 \alpha_2 \alpha_3 \alpha_6^2 - 17 \alpha_3^2 \alpha_6^2 - 
    50 \alpha_0 \alpha_4 \alpha_6^2 - 
    39 \alpha_1 \alpha_4 \alpha_6^2 - 
    28 \alpha_2 \alpha_4 \alpha_6^2 - 
    17 \alpha_3 \alpha_4 \alpha_6^2 - 3 \alpha_4^2 \alpha_6^2 - 
    64 \alpha_0 \alpha_6^3 - 48 \alpha_1 \alpha_6^3 - 
    36 \alpha_2 \alpha_6^3 - 24 \alpha_3 \alpha_6^3 - 
    12 \alpha_4 \alpha_6^3 - 11 \alpha_6^4 - 
    64 \alpha_0^3 \alpha_7 - 136 \alpha_0^2 \alpha_1 \alpha_7 - 
    100 \alpha_0 \alpha_1^2 \alpha_7 - 28 \alpha_1^3 \alpha_7 - 
    102 \alpha_0^2 \alpha_2 \alpha_7 - 
    150 \alpha_0 \alpha_1 \alpha_2 \alpha_7 - 
    63 \alpha_1^2 \alpha_2 \alpha_7 - 
    54 \alpha_0 \alpha_2^2 \alpha_7 - 
    45 \alpha_1 \alpha_2^2 \alpha_7 - 10 \alpha_2^3 \alpha_7 - 
    68 \alpha_0^2 \alpha_3 \alpha_7 - 
    100 \alpha_0 \alpha_1 \alpha_3 \alpha_7 - 
    42 \alpha_1^2 \alpha_3 \alpha_7 - 
    72 \alpha_0 \alpha_2 \alpha_3 \alpha_7 - 
    60 \alpha_1 \alpha_2 \alpha_3 \alpha_7 - 
    20 \alpha_2^2 \alpha_3 \alpha_7 - 
    22 \alpha_0 \alpha_3^2 \alpha_7 - 
    18 \alpha_1 \alpha_3^2 \alpha_7 - 
    12 \alpha_2 \alpha_3^2 \alpha_7 - 2 \alpha_3^3 \alpha_7 - 
    34 \alpha_0^2 \alpha_4 \alpha_7 - 
    50 \alpha_0 \alpha_1 \alpha_4 \alpha_7 - 
    21 \alpha_1^2 \alpha_4 \alpha_7 - 
    36 \alpha_0 \alpha_2 \alpha_4 \alpha_7 - 
    30 \alpha_1 \alpha_2 \alpha_4 \alpha_7 - 
    10 \alpha_2^2 \alpha_4 \alpha_7 - 
    22 \alpha_0 \alpha_3 \alpha_4 \alpha_7 - 
    18 \alpha_1 \alpha_3 \alpha_4 \alpha_7 - 
    12 \alpha_2 \alpha_3 \alpha_4 \alpha_7 - 
    3 \alpha_3^2 \alpha_4 \alpha_7 - 
    4 \alpha_0 \alpha_4^2 \alpha_7 - 
    3 \alpha_1 \alpha_4^2 \alpha_7 - 
    2 \alpha_2 \alpha_4^2 \alpha_7 - \alpha_3 \alpha_4^2 \
\alpha_7 - 137 \alpha_0^2 \alpha_6 \alpha_7 - 
    200 \alpha_0 \alpha_1 \alpha_6 \alpha_7 - 
    78 \alpha_1^2 \alpha_6 \alpha_7 - 
    150 \alpha_0 \alpha_2 \alpha_6 \alpha_7 - 
    117 \alpha_1 \alpha_2 \alpha_6 \alpha_7 - 
    42 \alpha_2^2 \alpha_6 \alpha_7 - 
    100 \alpha_0 \alpha_3 \alpha_6 \alpha_7 - 
    78 \alpha_1 \alpha_3 \alpha_6 \alpha_7 - 
    56 \alpha_2 \alpha_3 \alpha_6 \alpha_7 - 
    17 \alpha_3^2 \alpha_6 \alpha_7 - 
    50 \alpha_0 \alpha_4 \alpha_6 \alpha_7 - 
    39 \alpha_1 \alpha_4 \alpha_6 \alpha_7 - 
    28 \alpha_2 \alpha_4 \alpha_6 \alpha_7 - 
    17 \alpha_3 \alpha_4 \alpha_6 \alpha_7 - 
    3 \alpha_4^2 \alpha_6 \alpha_7 - 
    96 \alpha_0 \alpha_6^2 \alpha_7 - 
    72 \alpha_1 \alpha_6^2 \alpha_7 - 
    54 \alpha_2 \alpha_6^2 \alpha_7 - 
    36 \alpha_3 \alpha_6^2 \alpha_7 - 
    18 \alpha_4 \alpha_6^2 \alpha_7 - 22 \alpha_6^3 \alpha_7 - 
    37 \alpha_0^2 \alpha_7^2 - 56 \alpha_0 \alpha_1 \alpha_7^2 - 
    22 \alpha_1^2 \alpha_7^2 - 42 \alpha_0 \alpha_2 \alpha_7^2 - 
    33 \alpha_1 \alpha_2 \alpha_7^2 - 12 \alpha_2^2 \alpha_7^2 - 
    28 \alpha_0 \alpha_3 \alpha_7^2 - 
    22 \alpha_1 \alpha_3 \alpha_7^2 - 
    16 \alpha_2 \alpha_3 \alpha_7^2 - 5 \alpha_3^2 \alpha_7^2 - 
    14 \alpha_0 \alpha_4 \alpha_7^2 - 
    11 \alpha_1 \alpha_4 \alpha_7^2 - 
    8 \alpha_2 \alpha_4 \alpha_7^2 - 
    5 \alpha_3 \alpha_4 \alpha_7^2 - \alpha_4^2 \alpha_7^2 - 
    52 \alpha_0 \alpha_6 \alpha_7^2 - 
    40 \alpha_1 \alpha_6 \alpha_7^2 - 
    30 \alpha_2 \alpha_6 \alpha_7^2 - 
    20 \alpha_3 \alpha_6 \alpha_7^2 - 
    10 \alpha_4 \alpha_6 \alpha_7^2 - 18 \alpha_6^2 \alpha_7^2 - 
    10 \alpha_0 \alpha_7^3 - 8 \alpha_1 \alpha_7^3 - 
    6 \alpha_2 \alpha_7^3 - 4 \alpha_3 \alpha_7^3 - 
    2 \alpha_4 \alpha_7^3 - 7 \alpha_6 \alpha_7^3 - \alpha_7^4 - 
    88 \alpha_0^3 \alpha_8 - 180 \alpha_0^2 \alpha_1 \alpha_8 - 
    122 \alpha_0 \alpha_1^2 \alpha_8 - 32 \alpha_1^3 \alpha_8 - 
    135 \alpha_0^2 \alpha_2 \alpha_8 - 
    183 \alpha_0 \alpha_1 \alpha_2 \alpha_8 - 
    72 \alpha_1^2 \alpha_2 \alpha_8 - 
    66 \alpha_0 \alpha_2^2 \alpha_8 - 
    51 \alpha_1 \alpha_2^2 \alpha_8 - 11 \alpha_2^3 \alpha_8 - 
    90 \alpha_0^2 \alpha_3 \alpha_8 - 
    122 \alpha_0 \alpha_1 \alpha_3 \alpha_8 - 
    48 \alpha_1^2 \alpha_3 \alpha_8 - 
    88 \alpha_0 \alpha_2 \alpha_3 \alpha_8 - 
    68 \alpha_1 \alpha_2 \alpha_3 \alpha_8 - 
    22 \alpha_2^2 \alpha_3 \alpha_8 - 
    27 \alpha_0 \alpha_3^2 \alpha_8 - 
    20 \alpha_1 \alpha_3^2 \alpha_8 - 
    13 \alpha_2 \alpha_3^2 \alpha_8 - 2 \alpha_3^3 \alpha_8 - 
    45 \alpha_0^2 \alpha_4 \alpha_8 - 
    61 \alpha_0 \alpha_1 \alpha_4 \alpha_8 - 
    24 \alpha_1^2 \alpha_4 \alpha_8 - 
    44 \alpha_0 \alpha_2 \alpha_4 \alpha_8 - 
    34 \alpha_1 \alpha_2 \alpha_4 \alpha_8 - 
    11 \alpha_2^2 \alpha_4 \alpha_8 - 
    27 \alpha_0 \alpha_3 \alpha_4 \alpha_8 - 
    20 \alpha_1 \alpha_3 \alpha_4 \alpha_8 - 
    13 \alpha_2 \alpha_3 \alpha_4 \alpha_8 - 
    3 \alpha_3^2 \alpha_4 \alpha_8 - 
    5 \alpha_0 \alpha_4^2 \alpha_8 - 
    3 \alpha_1 \alpha_4^2 \alpha_8 - 
    2 \alpha_2 \alpha_4^2 \alpha_8 - \alpha_3 \alpha_4^2 \
\alpha_8 - 192 \alpha_0^2 \alpha_6 \alpha_8 - 
    272 \alpha_0 \alpha_1 \alpha_6 \alpha_8 - 
    100 \alpha_1^2 \alpha_6 \alpha_8 - 
    204 \alpha_0 \alpha_2 \alpha_6 \alpha_8 - 
    150 \alpha_1 \alpha_2 \alpha_6 \alpha_8 - 
    54 \alpha_2^2 \alpha_6 \alpha_8 - 
    136 \alpha_0 \alpha_3 \alpha_6 \alpha_8 - 
    100 \alpha_1 \alpha_3 \alpha_6 \alpha_8 - 
    72 \alpha_2 \alpha_3 \alpha_6 \alpha_8 - 
    22 \alpha_3^2 \alpha_6 \alpha_8 - 
    68 \alpha_0 \alpha_4 \alpha_6 \alpha_8 - 
    50 \alpha_1 \alpha_4 \alpha_6 \alpha_8 - 
    36 \alpha_2 \alpha_4 \alpha_6 \alpha_8 - 
    22 \alpha_3 \alpha_4 \alpha_6 \alpha_8 - 
    4 \alpha_4^2 \alpha_6 \alpha_8 - 
    137 \alpha_0 \alpha_6^2 \alpha_8 - 
    100 \alpha_1 \alpha_6^2 \alpha_8 - 
    75 \alpha_2 \alpha_6^2 \alpha_8 - 
    50 \alpha_3 \alpha_6^2 \alpha_8 - 
    25 \alpha_4 \alpha_6^2 \alpha_8 - 32 \alpha_6^3 \alpha_8 - 
    96 \alpha_0^2 \alpha_7 \alpha_8 - 
    136 \alpha_0 \alpha_1 \alpha_7 \alpha_8 - 
    50 \alpha_1^2 \alpha_7 \alpha_8 - 
    102 \alpha_0 \alpha_2 \alpha_7 \alpha_8 - 
    75 \alpha_1 \alpha_2 \alpha_7 \alpha_8 - 
    27 \alpha_2^2 \alpha_7 \alpha_8 - 
    68 \alpha_0 \alpha_3 \alpha_7 \alpha_8 - 
    50 \alpha_1 \alpha_3 \alpha_7 \alpha_8 - 
    36 \alpha_2 \alpha_3 \alpha_7 \alpha_8 - 
    11 \alpha_3^2 \alpha_7 \alpha_8 - 
    34 \alpha_0 \alpha_4 \alpha_7 \alpha_8 - 
    25 \alpha_1 \alpha_4 \alpha_7 \alpha_8 - 
    18 \alpha_2 \alpha_4 \alpha_7 \alpha_8 - 
    11 \alpha_3 \alpha_4 \alpha_7 \alpha_8 - 
    2 \alpha_4^2 \alpha_7 \alpha_8 - 
    137 \alpha_0 \alpha_6 \alpha_7 \alpha_8 - 
    100 \alpha_1 \alpha_6 \alpha_7 \alpha_8 - 
    75 \alpha_2 \alpha_6 \alpha_7 \alpha_8 - 
    50 \alpha_3 \alpha_6 \alpha_7 \alpha_8 - 
    25 \alpha_4 \alpha_6 \alpha_7 \alpha_8 - 
    48 \alpha_6^2 \alpha_7 \alpha_8 - 
    37 \alpha_0 \alpha_7^2 \alpha_8 - 
    28 \alpha_1 \alpha_7^2 \alpha_8 - 
    21 \alpha_2 \alpha_7^2 \alpha_8 - 
    14 \alpha_3 \alpha_7^2 \alpha_8 - 
    7 \alpha_4 \alpha_7^2 \alpha_8 - 
    26 \alpha_6 \alpha_7^2 \alpha_8 - 5 \alpha_7^3 \alpha_8 - 
    63 \alpha_0^2 \alpha_8^2 - 84 \alpha_0 \alpha_1 \alpha_8^2 - 
    28 \alpha_1^2 \alpha_8^2 - 63 \alpha_0 \alpha_2 \alpha_8^2 - 
    42 \alpha_1 \alpha_2 \alpha_8^2 - 15 \alpha_2^2 \alpha_8^2 - 
    42 \alpha_0 \alpha_3 \alpha_8^2 - 
    28 \alpha_1 \alpha_3 \alpha_8^2 - 
    20 \alpha_2 \alpha_3 \alpha_8^2 - 6 \alpha_3^2 \alpha_8^2 - 
    21 \alpha_0 \alpha_4 \alpha_8^2 - 
    14 \alpha_1 \alpha_4 \alpha_8^2 - 
    10 \alpha_2 \alpha_4 \alpha_8^2 - 
    6 \alpha_3 \alpha_4 \alpha_8^2 - \alpha_4^2 \alpha_8^2 - 
    92 \alpha_0 \alpha_6 \alpha_8^2 - 
    64 \alpha_1 \alpha_6 \alpha_8^2 - 
    48 \alpha_2 \alpha_6 \alpha_8^2 - 
    32 \alpha_3 \alpha_6 \alpha_8^2 - 
    16 \alpha_4 \alpha_6 \alpha_8^2 - 33 \alpha_6^2 \alpha_8^2 - 
    46 \alpha_0 \alpha_7 \alpha_8^2 - 
    32 \alpha_1 \alpha_7 \alpha_8^2 - 
    24 \alpha_2 \alpha_7 \alpha_8^2 - 
    16 \alpha_3 \alpha_7 \alpha_8^2 - 
    8 \alpha_4 \alpha_7 \alpha_8^2 - 
    33 \alpha_6 \alpha_7 \alpha_8^2 - 9 \alpha_7^2 \alpha_8^2 - 
    19 \alpha_0 \alpha_8^3 - 12 \alpha_1 \alpha_8^3 - 
    9 \alpha_2 \alpha_8^3 - 6 \alpha_3 \alpha_8^3 - 
    3 \alpha_4 \alpha_8^3 - 14 \alpha_6 \alpha_8^3 - 
    7 \alpha_7 \alpha_8^3 - 2 \alpha_8^4 + 
    q^2 \alpha_0^2 (\alpha_0 + \alpha_8)^2 + 
    q \alpha_0 (\alpha_0 + \alpha_8) (11 \alpha_0^2 + 
       24 \alpha_0 \alpha_1 + 10 \alpha_1^2 + 
       18 \alpha_0 \alpha_2 + 15 \alpha_1 \alpha_2 + 
       6 \alpha_2^2 + 12 \alpha_0 \alpha_3 + 
       10 \alpha_1 \alpha_3 + 8 \alpha_2 \alpha_3 + 
       3 \alpha_3^2 + 6 \alpha_0 \alpha_4 + 
       5 \alpha_1 \alpha_4 + 4 \alpha_2 \alpha_4 + 
       3 \alpha_3 \alpha_4 + \alpha_4^2 + 16 \alpha_0 \alpha_6 + 
       16 \alpha_1 \alpha_6 + 12 \alpha_2 \alpha_6 + 
       8 \alpha_3 \alpha_6 + 4 \alpha_4 \alpha_6 + 
       6 \alpha_6^2 + 8 \alpha_0 \alpha_7 + 
       8 \alpha_1 \alpha_7 + 6 \alpha_2 \alpha_7 + 
       4 \alpha_3 \alpha_7 + 2 \alpha_4 \alpha_7 + 
       6 \alpha_6 \alpha_7 + \alpha_7^2 + 11 \alpha_0 \alpha_8 + 
       12 \alpha_1 \alpha_8 + 9 \alpha_2 \alpha_8 + 
       6 \alpha_3 \alpha_8 + 3 \alpha_4 \alpha_8 + 
       8 \alpha_6 \alpha_8 + 4 \alpha_7 \alpha_8 + 
       3 \alpha_8^2)) + 
 p^3 (-1 + q) q^2 (-36 \alpha_0 \alpha_1^2 - 20 \alpha_1^3 - 
    54 \alpha_0 \alpha_1 \alpha_2 - 45 \alpha_1^2 \alpha_2 - 
    18 \alpha_0 \alpha_2^2 - 33 \alpha_1 \alpha_2^2 - 
    8 \alpha_2^3 - 36 \alpha_0 \alpha_1 \alpha_3 - 
    30 \alpha_1^2 \alpha_3 - 24 \alpha_0 \alpha_2 \alpha_3 - 
    44 \alpha_1 \alpha_2 \alpha_3 - 16 \alpha_2^2 \alpha_3 - 
    6 \alpha_0 \alpha_3^2 - 14 \alpha_1 \alpha_3^2 - 
    10 \alpha_2 \alpha_3^2 - 2 \alpha_3^3 - 
    18 \alpha_0 \alpha_1 \alpha_4 - 15 \alpha_1^2 \alpha_4 - 
    12 \alpha_0 \alpha_2 \alpha_4 - 
    22 \alpha_1 \alpha_2 \alpha_4 - 8 \alpha_2^2 \alpha_4 - 
    6 \alpha_0 \alpha_3 \alpha_4 - 
    14 \alpha_1 \alpha_3 \alpha_4 - 
    10 \alpha_2 \alpha_3 \alpha_4 - 3 \alpha_3^2 \alpha_4 - 
    3 \alpha_1 \alpha_4^2 - 
    2 \alpha_2 \alpha_4^2 - \alpha_3 \alpha_4^2 - 
    24 \alpha_1^2 \alpha_6 - 36 \alpha_1 \alpha_2 \alpha_6 - 
    12 \alpha_2^2 \alpha_6 - 24 \alpha_1 \alpha_3 \alpha_6 - 
    16 \alpha_2 \alpha_3 \alpha_6 - 4 \alpha_3^2 \alpha_6 - 
    12 \alpha_1 \alpha_4 \alpha_6 - 
    8 \alpha_2 \alpha_4 \alpha_6 - 
    4 \alpha_3 \alpha_4 \alpha_6 - 12 \alpha_1^2 \alpha_7 - 
    18 \alpha_1 \alpha_2 \alpha_7 - 6 \alpha_2^2 \alpha_7 - 
    12 \alpha_1 \alpha_3 \alpha_7 - 
    8 \alpha_2 \alpha_3 \alpha_7 - 2 \alpha_3^2 \alpha_7 - 
    6 \alpha_1 \alpha_4 \alpha_7 - 
    4 \alpha_2 \alpha_4 \alpha_7 - 
    2 \alpha_3 \alpha_4 \alpha_7 - 18 \alpha_1^2 \alpha_8 - 
    27 \alpha_1 \alpha_2 \alpha_8 - 9 \alpha_2^2 \alpha_8 - 
    18 \alpha_1 \alpha_3 \alpha_8 - 
    12 \alpha_2 \alpha_3 \alpha_8 - 3 \alpha_3^2 \alpha_8 - 
    9 \alpha_1 \alpha_4 \alpha_8 - 
    6 \alpha_2 \alpha_4 \alpha_8 - 
    3 \alpha_3 \alpha_4 \alpha_8 + 
    q^3 (2 \alpha_0 + \alpha_8) (10 \alpha_0^2 + 
       10 \alpha_0 \alpha_8 + \alpha_8^2) + 
    q (-20 \alpha_0^3 - 96 \alpha_0^2 \alpha_1 - 
       4 \alpha_0 \alpha_1^2 + 20 \alpha_1^3 - 
       72 \alpha_0^2 \alpha_2 - 6 \alpha_0 \alpha_1 \alpha_2 + 
       45 \alpha_1^2 \alpha_2 - 6 \alpha_0 \alpha_2^2 + 
       33 \alpha_1 \alpha_2^2 + 8 \alpha_2^3 - 
       48 \alpha_0^2 \alpha_3 - 4 \alpha_0 \alpha_1 \alpha_3 + 
       30 \alpha_1^2 \alpha_3 - 8 \alpha_0 \alpha_2 \alpha_3 + 
       44 \alpha_1 \alpha_2 \alpha_3 + 16 \alpha_2^2 \alpha_3 - 
       6 \alpha_0 \alpha_3^2 + 14 \alpha_1 \alpha_3^2 + 
       10 \alpha_2 \alpha_3^2 + 2 \alpha_3^3 - 
       24 \alpha_0^2 \alpha_4 - 2 \alpha_0 \alpha_1 \alpha_4 + 
       15 \alpha_1^2 \alpha_4 - 4 \alpha_0 \alpha_2 \alpha_4 + 
       22 \alpha_1 \alpha_2 \alpha_4 + 8 \alpha_2^2 \alpha_4 - 
       6 \alpha_0 \alpha_3 \alpha_4 + 
       14 \alpha_1 \alpha_3 \alpha_4 + 
       10 \alpha_2 \alpha_3 \alpha_4 + 3 \alpha_3^2 \alpha_4 - 
       4 \alpha_0 \alpha_4^2 + 3 \alpha_1 \alpha_4^2 + 
       2 \alpha_2 \alpha_4^2 + \alpha_3 \alpha_4^2 - 
       40 \alpha_0^2 \alpha_6 - 64 \alpha_0 \alpha_1 \alpha_6 + 
       24 \alpha_1^2 \alpha_6 - 48 \alpha_0 \alpha_2 \alpha_6 + 
       36 \alpha_1 \alpha_2 \alpha_6 + 12 \alpha_2^2 \alpha_6 - 
       32 \alpha_0 \alpha_3 \alpha_6 + 
       24 \alpha_1 \alpha_3 \alpha_6 + 
       16 \alpha_2 \alpha_3 \alpha_6 + 4 \alpha_3^2 \alpha_6 - 
       16 \alpha_0 \alpha_4 \alpha_6 + 
       12 \alpha_1 \alpha_4 \alpha_6 + 
       8 \alpha_2 \alpha_4 \alpha_6 + 
       4 \alpha_3 \alpha_4 \alpha_6 - 24 \alpha_0 \alpha_6^2 - 
       4 \alpha_6^3 - 20 \alpha_0^2 \alpha_7 - 
       32 \alpha_0 \alpha_1 \alpha_7 + 12 \alpha_1^2 \alpha_7 - 
       24 \alpha_0 \alpha_2 \alpha_7 + 
       18 \alpha_1 \alpha_2 \alpha_7 + 6 \alpha_2^2 \alpha_7 - 
       16 \alpha_0 \alpha_3 \alpha_7 + 
       12 \alpha_1 \alpha_3 \alpha_7 + 
       8 \alpha_2 \alpha_3 \alpha_7 + 2 \alpha_3^2 \alpha_7 - 
       8 \alpha_0 \alpha_4 \alpha_7 + 
       6 \alpha_1 \alpha_4 \alpha_7 + 
       4 \alpha_2 \alpha_4 \alpha_7 + 
       2 \alpha_3 \alpha_4 \alpha_7 - 
       24 \alpha_0 \alpha_6 \alpha_7 - 6 \alpha_6^2 \alpha_7 - 
       4 \alpha_0 \alpha_7^2 - 2 \alpha_6 \alpha_7^2 - 
       30 \alpha_0^2 \alpha_8 - 96 \alpha_0 \alpha_1 \alpha_8 - 
       2 \alpha_1^2 \alpha_8 - 72 \alpha_0 \alpha_2 \alpha_8 - 
       3 \alpha_1 \alpha_2 \alpha_8 - 3 \alpha_2^2 \alpha_8 - 
       48 \alpha_0 \alpha_3 \alpha_8 - 
       2 \alpha_1 \alpha_3 \alpha_8 - 
       4 \alpha_2 \alpha_3 \alpha_8 - 3 \alpha_3^2 \alpha_8 - 
       24 \alpha_0 \alpha_4 \alpha_8 - \alpha_1 \alpha_4 \
\alpha_8 - 2 \alpha_2 \alpha_4 \alpha_8 - 
       3 \alpha_3 \alpha_4 \alpha_8 - 2 \alpha_4^2 \alpha_8 - 
       40 \alpha_0 \alpha_6 \alpha_8 - 
       32 \alpha_1 \alpha_6 \alpha_8 - 
       24 \alpha_2 \alpha_6 \alpha_8 - 
       16 \alpha_3 \alpha_6 \alpha_8 - 
       8 \alpha_4 \alpha_6 \alpha_8 - 12 \alpha_6^2 \alpha_8 - 
       20 \alpha_0 \alpha_7 \alpha_8 - 
       16 \alpha_1 \alpha_7 \alpha_8 - 
       12 \alpha_2 \alpha_7 \alpha_8 - 
       8 \alpha_3 \alpha_7 \alpha_8 - 
       4 \alpha_4 \alpha_7 \alpha_8 - 
       12 \alpha_6 \alpha_7 \alpha_8 - 2 \alpha_7^2 \alpha_8 - 
       18 \alpha_0 \alpha_8^2 - 24 \alpha_1 \alpha_8^2 - 
       18 \alpha_2 \alpha_8^2 - 12 \alpha_3 \alpha_8^2 - 
       6 \alpha_4 \alpha_8^2 - 12 \alpha_6 \alpha_8^2 - 
       6 \alpha_7 \alpha_8^2 - 4 \alpha_8^3) + 
    q^2 (96 \alpha_0^2 \alpha_1 + 40 \alpha_0 \alpha_1^2 + 
       72 \alpha_0^2 \alpha_2 + 60 \alpha_0 \alpha_1 \alpha_2 + 
       24 \alpha_0 \alpha_2^2 + 48 \alpha_0^2 \alpha_3 + 
       40 \alpha_0 \alpha_1 \alpha_3 + 
       32 \alpha_0 \alpha_2 \alpha_3 + 12 \alpha_0 \alpha_3^2 + 
       24 \alpha_0^2 \alpha_4 + 20 \alpha_0 \alpha_1 \alpha_4 + 
       16 \alpha_0 \alpha_2 \alpha_4 + 
       12 \alpha_0 \alpha_3 \alpha_4 + 4 \alpha_0 \alpha_4^2 + 
       40 \alpha_0^2 \alpha_6 + 64 \alpha_0 \alpha_1 \alpha_6 + 
       48 \alpha_0 \alpha_2 \alpha_6 + 
       32 \alpha_0 \alpha_3 \alpha_6 + 
       16 \alpha_0 \alpha_4 \alpha_6 + 24 \alpha_0 \alpha_6^2 + 
       20 \alpha_0^2 \alpha_7 + 32 \alpha_0 \alpha_1 \alpha_7 + 
       24 \alpha_0 \alpha_2 \alpha_7 + 
       16 \alpha_0 \alpha_3 \alpha_7 + 
       8 \alpha_0 \alpha_4 \alpha_7 + 
       24 \alpha_0 \alpha_6 \alpha_7 + 4 \alpha_0 \alpha_7^2 + 
       96 \alpha_0 \alpha_1 \alpha_8 + 20 \alpha_1^2 \alpha_8 + 
       72 \alpha_0 \alpha_2 \alpha_8 + 
       30 \alpha_1 \alpha_2 \alpha_8 + 12 \alpha_2^2 \alpha_8 + 
       48 \alpha_0 \alpha_3 \alpha_8 + 
       20 \alpha_1 \alpha_3 \alpha_8 + 
       16 \alpha_2 \alpha_3 \alpha_8 + 6 \alpha_3^2 \alpha_8 + 
       24 \alpha_0 \alpha_4 \alpha_8 + 
       10 \alpha_1 \alpha_4 \alpha_8 + 
       8 \alpha_2 \alpha_4 \alpha_8 + 
       6 \alpha_3 \alpha_4 \alpha_8 + 2 \alpha_4^2 \alpha_8 + 
       40 \alpha_0 \alpha_6 \alpha_8 + 
       32 \alpha_1 \alpha_6 \alpha_8 + 
       24 \alpha_2 \alpha_6 \alpha_8 + 
       16 \alpha_3 \alpha_6 \alpha_8 + 
       8 \alpha_4 \alpha_6 \alpha_8 + 12 \alpha_6^2 \alpha_8 + 
       20 \alpha_0 \alpha_7 \alpha_8 + 
       16 \alpha_1 \alpha_7 \alpha_8 + 
       12 \alpha_2 \alpha_7 \alpha_8 + 
       8 \alpha_3 \alpha_7 \alpha_8 + 
       4 \alpha_4 \alpha_7 \alpha_8 + 
       12 \alpha_6 \alpha_7 \alpha_8 + 2 \alpha_7^2 \alpha_8 + 
       6 \alpha_0 \alpha_8^2 + 24 \alpha_1 \alpha_8^2 + 
       18 \alpha_2 \alpha_8^2 + 12 \alpha_3 \alpha_8^2 + 
       6 \alpha_4 \alpha_8^2 + 12 \alpha_6 \alpha_8^2 + 
       6 \alpha_7 \alpha_8^2 + 3 \alpha_8^3)) + 
 p^2 q (-24 \alpha_0 \alpha_1^3 - 15 \alpha_1^4 - 
    54 \alpha_0 \alpha_1^2 \alpha_2 - 45 \alpha_1^3 \alpha_2 - 
    36 \alpha_0 \alpha_1 \alpha_2^2 - 48 \alpha_1^2 \alpha_2^2 - 
    6 \alpha_0 \alpha_2^3 - 21 \alpha_1 \alpha_2^3 - 
    3 \alpha_2^4 - 36 \alpha_0 \alpha_1^2 \alpha_3 - 
    30 \alpha_1^3 \alpha_3 - 
    48 \alpha_0 \alpha_1 \alpha_2 \alpha_3 - 
    64 \alpha_1^2 \alpha_2 \alpha_3 - 
    12 \alpha_0 \alpha_2^2 \alpha_3 - 
    42 \alpha_1 \alpha_2^2 \alpha_3 - 8 \alpha_2^3 \alpha_3 - 
    12 \alpha_0 \alpha_1 \alpha_3^2 - 19 \alpha_1^2 \alpha_3^2 - 
    6 \alpha_0 \alpha_2 \alpha_3^2 - 
    25 \alpha_1 \alpha_2 \alpha_3^2 - 7 \alpha_2^2 \alpha_3^2 - 
    4 \alpha_1 \alpha_3^3 - 2 \alpha_2 \alpha_3^3 - 
    18 \alpha_0 \alpha_1^2 \alpha_4 - 15 \alpha_1^3 \alpha_4 - 
    24 \alpha_0 \alpha_1 \alpha_2 \alpha_4 - 
    32 \alpha_1^2 \alpha_2 \alpha_4 - 
    6 \alpha_0 \alpha_2^2 \alpha_4 - 
    21 \alpha_1 \alpha_2^2 \alpha_4 - 4 \alpha_2^3 \alpha_4 - 
    12 \alpha_0 \alpha_1 \alpha_3 \alpha_4 - 
    19 \alpha_1^2 \alpha_3 \alpha_4 - 
    6 \alpha_0 \alpha_2 \alpha_3 \alpha_4 - 
    25 \alpha_1 \alpha_2 \alpha_3 \alpha_4 - 
    7 \alpha_2^2 \alpha_3 \alpha_4 - 
    6 \alpha_1 \alpha_3^2 \alpha_4 - 
    3 \alpha_2 \alpha_3^2 \alpha_4 - 3 \alpha_1^2 \alpha_4^2 - 
    4 \alpha_1 \alpha_2 \alpha_4^2 - \alpha_2^2 \alpha_4^2 - 
    2 \alpha_1 \alpha_3 \alpha_4^2 - \alpha_2 \alpha_3 \
\alpha_4^2 - 16 \alpha_1^3 \alpha_6 - 
    36 \alpha_1^2 \alpha_2 \alpha_6 - 
    24 \alpha_1 \alpha_2^2 \alpha_6 - 4 \alpha_2^3 \alpha_6 - 
    24 \alpha_1^2 \alpha_3 \alpha_6 - 
    32 \alpha_1 \alpha_2 \alpha_3 \alpha_6 - 
    8 \alpha_2^2 \alpha_3 \alpha_6 - 
    8 \alpha_1 \alpha_3^2 \alpha_6 - 
    4 \alpha_2 \alpha_3^2 \alpha_6 - 
    12 \alpha_1^2 \alpha_4 \alpha_6 - 
    16 \alpha_1 \alpha_2 \alpha_4 \alpha_6 - 
    4 \alpha_2^2 \alpha_4 \alpha_6 - 
    8 \alpha_1 \alpha_3 \alpha_4 \alpha_6 - 
    4 \alpha_2 \alpha_3 \alpha_4 \alpha_6 - 
    8 \alpha_1^3 \alpha_7 - 18 \alpha_1^2 \alpha_2 \alpha_7 - 
    12 \alpha_1 \alpha_2^2 \alpha_7 - 2 \alpha_2^3 \alpha_7 - 
    12 \alpha_1^2 \alpha_3 \alpha_7 - 
    16 \alpha_1 \alpha_2 \alpha_3 \alpha_7 - 
    4 \alpha_2^2 \alpha_3 \alpha_7 - 
    4 \alpha_1 \alpha_3^2 \alpha_7 - 
    2 \alpha_2 \alpha_3^2 \alpha_7 - 
    6 \alpha_1^2 \alpha_4 \alpha_7 - 
    8 \alpha_1 \alpha_2 \alpha_4 \alpha_7 - 
    2 \alpha_2^2 \alpha_4 \alpha_7 - 
    4 \alpha_1 \alpha_3 \alpha_4 \alpha_7 - 
    2 \alpha_2 \alpha_3 \alpha_4 \alpha_7 - 
    12 \alpha_1^3 \alpha_8 - 27 \alpha_1^2 \alpha_2 \alpha_8 - 
    18 \alpha_1 \alpha_2^2 \alpha_8 - 3 \alpha_2^3 \alpha_8 - 
    18 \alpha_1^2 \alpha_3 \alpha_8 - 
    24 \alpha_1 \alpha_2 \alpha_3 \alpha_8 - 
    6 \alpha_2^2 \alpha_3 \alpha_8 - 
    6 \alpha_1 \alpha_3^2 \alpha_8 - 
    3 \alpha_2 \alpha_3^2 \alpha_8 - 
    9 \alpha_1^2 \alpha_4 \alpha_8 - 
    12 \alpha_1 \alpha_2 \alpha_4 \alpha_8 - 
    3 \alpha_2^2 \alpha_4 \alpha_8 - 
    6 \alpha_1 \alpha_3 \alpha_4 \alpha_8 - 
    3 \alpha_2 \alpha_3 \alpha_4 \alpha_8 + 
    3 q^4 \alpha_0 (\alpha_0 + \alpha_8) (5 \alpha_0^2 + 
       5 \alpha_0 \alpha_8 + \alpha_8^2) + 
    q^2 (-105 \alpha_0^4 - 288 \alpha_0^3 \alpha_1 - 
       120 \alpha_0^2 \alpha_1^2 - 24 \alpha_0 \alpha_1^3 - 
       15 \alpha_1^4 - 216 \alpha_0^3 \alpha_2 - 
       180 \alpha_0^2 \alpha_1 \alpha_2 - 
       54 \alpha_0 \alpha_1^2 \alpha_2 - 
       45 \alpha_1^3 \alpha_2 - 72 \alpha_0^2 \alpha_2^2 - 
       36 \alpha_0 \alpha_1 \alpha_2^2 - 
       48 \alpha_1^2 \alpha_2^2 - 6 \alpha_0 \alpha_2^3 - 
       21 \alpha_1 \alpha_2^3 - 3 \alpha_2^4 - 
       144 \alpha_0^3 \alpha_3 - 
       120 \alpha_0^2 \alpha_1 \alpha_3 - 
       36 \alpha_0 \alpha_1^2 \alpha_3 - 
       30 \alpha_1^3 \alpha_3 - 
       96 \alpha_0^2 \alpha_2 \alpha_3 - 
       48 \alpha_0 \alpha_1 \alpha_2 \alpha_3 - 
       64 \alpha_1^2 \alpha_2 \alpha_3 - 
       12 \alpha_0 \alpha_2^2 \alpha_3 - 
       42 \alpha_1 \alpha_2^2 \alpha_3 - 8 \alpha_2^3 \alpha_3 - 
       36 \alpha_0^2 \alpha_3^2 - 
       12 \alpha_0 \alpha_1 \alpha_3^2 - 
       19 \alpha_1^2 \alpha_3^2 - 
       6 \alpha_0 \alpha_2 \alpha_3^2 - 
       25 \alpha_1 \alpha_2 \alpha_3^2 - 
       7 \alpha_2^2 \alpha_3^2 - 4 \alpha_1 \alpha_3^3 - 
       2 \alpha_2 \alpha_3^3 - 72 \alpha_0^3 \alpha_4 - 
       60 \alpha_0^2 \alpha_1 \alpha_4 - 
       18 \alpha_0 \alpha_1^2 \alpha_4 - 
       15 \alpha_1^3 \alpha_4 - 
       48 \alpha_0^2 \alpha_2 \alpha_4 - 
       24 \alpha_0 \alpha_1 \alpha_2 \alpha_4 - 
       32 \alpha_1^2 \alpha_2 \alpha_4 - 
       6 \alpha_0 \alpha_2^2 \alpha_4 - 
       21 \alpha_1 \alpha_2^2 \alpha_4 - 4 \alpha_2^3 \alpha_4 - 
       36 \alpha_0^2 \alpha_3 \alpha_4 - 
       12 \alpha_0 \alpha_1 \alpha_3 \alpha_4 - 
       19 \alpha_1^2 \alpha_3 \alpha_4 - 
       6 \alpha_0 \alpha_2 \alpha_3 \alpha_4 - 
       25 \alpha_1 \alpha_2 \alpha_3 \alpha_4 - 
       7 \alpha_2^2 \alpha_3 \alpha_4 - 
       6 \alpha_1 \alpha_3^2 \alpha_4 - 
       3 \alpha_2 \alpha_3^2 \alpha_4 - 
       12 \alpha_0^2 \alpha_4^2 - 3 \alpha_1^2 \alpha_4^2 - 
       4 \alpha_1 \alpha_2 \alpha_4^2 - \alpha_2^2 \alpha_4^2 - 
       2 \alpha_1 \alpha_3 \alpha_4^2 - \alpha_2 \alpha_3 \
\alpha_4^2 - 240 \alpha_0^3 \alpha_6 - 
       384 \alpha_0^2 \alpha_1 \alpha_6 - 
       152 \alpha_0 \alpha_1^2 \alpha_6 - 
       56 \alpha_1^3 \alpha_6 - 
       288 \alpha_0^2 \alpha_2 \alpha_6 - 
       228 \alpha_0 \alpha_1 \alpha_2 \alpha_6 - 
       126 \alpha_1^2 \alpha_2 \alpha_6 - 
       84 \alpha_0 \alpha_2^2 \alpha_6 - 
       90 \alpha_1 \alpha_2^2 \alpha_6 - 
       20 \alpha_2^3 \alpha_6 - 
       192 \alpha_0^2 \alpha_3 \alpha_6 - 
       152 \alpha_0 \alpha_1 \alpha_3 \alpha_6 - 
       84 \alpha_1^2 \alpha_3 \alpha_6 - 
       112 \alpha_0 \alpha_2 \alpha_3 \alpha_6 - 
       120 \alpha_1 \alpha_2 \alpha_3 \alpha_6 - 
       40 \alpha_2^2 \alpha_3 \alpha_6 - 
       36 \alpha_0 \alpha_3^2 \alpha_6 - 
       36 \alpha_1 \alpha_3^2 \alpha_6 - 
       24 \alpha_2 \alpha_3^2 \alpha_6 - 4 \alpha_3^3 \alpha_6 - 
       96 \alpha_0^2 \alpha_4 \alpha_6 - 
       76 \alpha_0 \alpha_1 \alpha_4 \alpha_6 - 
       42 \alpha_1^2 \alpha_4 \alpha_6 - 
       56 \alpha_0 \alpha_2 \alpha_4 \alpha_6 - 
       60 \alpha_1 \alpha_2 \alpha_4 \alpha_6 - 
       20 \alpha_2^2 \alpha_4 \alpha_6 - 
       36 \alpha_0 \alpha_3 \alpha_4 \alpha_6 - 
       36 \alpha_1 \alpha_3 \alpha_4 \alpha_6 - 
       24 \alpha_2 \alpha_3 \alpha_4 \alpha_6 - 
       6 \alpha_3^2 \alpha_4 \alpha_6 - 
       8 \alpha_0 \alpha_4^2 \alpha_6 - 
       6 \alpha_1 \alpha_4^2 \alpha_6 - 
       4 \alpha_2 \alpha_4^2 \alpha_6 - 
       2 \alpha_3 \alpha_4^2 \alpha_6 - 
       197 \alpha_0^2 \alpha_6^2 - 
       200 \alpha_0 \alpha_1 \alpha_6^2 - 
       78 \alpha_1^2 \alpha_6^2 - 
       150 \alpha_0 \alpha_2 \alpha_6^2 - 
       117 \alpha_1 \alpha_2 \alpha_6^2 - 
       42 \alpha_2^2 \alpha_6^2 - 
       100 \alpha_0 \alpha_3 \alpha_6^2 - 
       78 \alpha_1 \alpha_3 \alpha_6^2 - 
       56 \alpha_2 \alpha_3 \alpha_6^2 - 
       17 \alpha_3^2 \alpha_6^2 - 
       50 \alpha_0 \alpha_4 \alpha_6^2 - 
       39 \alpha_1 \alpha_4 \alpha_6^2 - 
       28 \alpha_2 \alpha_4 \alpha_6^2 - 
       17 \alpha_3 \alpha_4 \alpha_6^2 - 
       3 \alpha_4^2 \alpha_6^2 - 72 \alpha_0 \alpha_6^3 - 
       48 \alpha_1 \alpha_6^3 - 36 \alpha_2 \alpha_6^3 - 
       24 \alpha_3 \alpha_6^3 - 12 \alpha_4 \alpha_6^3 - 
       11 \alpha_6^4 - 120 \alpha_0^3 \alpha_7 - 
       192 \alpha_0^2 \alpha_1 \alpha_7 - 
       76 \alpha_0 \alpha_1^2 \alpha_7 - 
       28 \alpha_1^3 \alpha_7 - 
       144 \alpha_0^2 \alpha_2 \alpha_7 - 
       114 \alpha_0 \alpha_1 \alpha_2 \alpha_7 - 
       63 \alpha_1^2 \alpha_2 \alpha_7 - 
       42 \alpha_0 \alpha_2^2 \alpha_7 - 
       45 \alpha_1 \alpha_2^2 \alpha_7 - 
       10 \alpha_2^3 \alpha_7 - 
       96 \alpha_0^2 \alpha_3 \alpha_7 - 
       76 \alpha_0 \alpha_1 \alpha_3 \alpha_7 - 
       42 \alpha_1^2 \alpha_3 \alpha_7 - 
       56 \alpha_0 \alpha_2 \alpha_3 \alpha_7 - 
       60 \alpha_1 \alpha_2 \alpha_3 \alpha_7 - 
       20 \alpha_2^2 \alpha_3 \alpha_7 - 
       18 \alpha_0 \alpha_3^2 \alpha_7 - 
       18 \alpha_1 \alpha_3^2 \alpha_7 - 
       12 \alpha_2 \alpha_3^2 \alpha_7 - 2 \alpha_3^3 \alpha_7 - 
       48 \alpha_0^2 \alpha_4 \alpha_7 - 
       38 \alpha_0 \alpha_1 \alpha_4 \alpha_7 - 
       21 \alpha_1^2 \alpha_4 \alpha_7 - 
       28 \alpha_0 \alpha_2 \alpha_4 \alpha_7 - 
       30 \alpha_1 \alpha_2 \alpha_4 \alpha_7 - 
       10 \alpha_2^2 \alpha_4 \alpha_7 - 
       18 \alpha_0 \alpha_3 \alpha_4 \alpha_7 - 
       18 \alpha_1 \alpha_3 \alpha_4 \alpha_7 - 
       12 \alpha_2 \alpha_3 \alpha_4 \alpha_7 - 
       3 \alpha_3^2 \alpha_4 \alpha_7 - 
       4 \alpha_0 \alpha_4^2 \alpha_7 - 
       3 \alpha_1 \alpha_4^2 \alpha_7 - 
       2 \alpha_2 \alpha_4^2 \alpha_7 - \alpha_3 \alpha_4^2 \
\alpha_7 - 197 \alpha_0^2 \alpha_6 \alpha_7 - 
       200 \alpha_0 \alpha_1 \alpha_6 \alpha_7 - 
       78 \alpha_1^2 \alpha_6 \alpha_7 - 
       150 \alpha_0 \alpha_2 \alpha_6 \alpha_7 - 
       117 \alpha_1 \alpha_2 \alpha_6 \alpha_7 - 
       42 \alpha_2^2 \alpha_6 \alpha_7 - 
       100 \alpha_0 \alpha_3 \alpha_6 \alpha_7 - 
       78 \alpha_1 \alpha_3 \alpha_6 \alpha_7 - 
       56 \alpha_2 \alpha_3 \alpha_6 \alpha_7 - 
       17 \alpha_3^2 \alpha_6 \alpha_7 - 
       50 \alpha_0 \alpha_4 \alpha_6 \alpha_7 - 
       39 \alpha_1 \alpha_4 \alpha_6 \alpha_7 - 
       28 \alpha_2 \alpha_4 \alpha_6 \alpha_7 - 
       17 \alpha_3 \alpha_4 \alpha_6 \alpha_7 - 
       3 \alpha_4^2 \alpha_6 \alpha_7 - 
       108 \alpha_0 \alpha_6^2 \alpha_7 - 
       72 \alpha_1 \alpha_6^2 \alpha_7 - 
       54 \alpha_2 \alpha_6^2 \alpha_7 - 
       36 \alpha_3 \alpha_6^2 \alpha_7 - 
       18 \alpha_4 \alpha_6^2 \alpha_7 - 
       22 \alpha_6^3 \alpha_7 - 47 \alpha_0^2 \alpha_7^2 - 
       56 \alpha_0 \alpha_1 \alpha_7^2 - 
       22 \alpha_1^2 \alpha_7^2 - 
       42 \alpha_0 \alpha_2 \alpha_7^2 - 
       33 \alpha_1 \alpha_2 \alpha_7^2 - 
       12 \alpha_2^2 \alpha_7^2 - 
       28 \alpha_0 \alpha_3 \alpha_7^2 - 
       22 \alpha_1 \alpha_3 \alpha_7^2 - 
       16 \alpha_2 \alpha_3 \alpha_7^2 - 
       5 \alpha_3^2 \alpha_7^2 - 
       14 \alpha_0 \alpha_4 \alpha_7^2 - 
       11 \alpha_1 \alpha_4 \alpha_7^2 - 
       8 \alpha_2 \alpha_4 \alpha_7^2 - 
       5 \alpha_3 \alpha_4 \alpha_7^2 - \alpha_4^2 \alpha_7^2 - 
       56 \alpha_0 \alpha_6 \alpha_7^2 - 
       40 \alpha_1 \alpha_6 \alpha_7^2 - 
       30 \alpha_2 \alpha_6 \alpha_7^2 - 
       20 \alpha_3 \alpha_6 \alpha_7^2 - 
       10 \alpha_4 \alpha_6 \alpha_7^2 - 
       18 \alpha_6^2 \alpha_7^2 - 10 \alpha_0 \alpha_7^3 - 
       8 \alpha_1 \alpha_7^3 - 6 \alpha_2 \alpha_7^3 - 
       4 \alpha_3 \alpha_7^3 - 2 \alpha_4 \alpha_7^3 - 
       7 \alpha_6 \alpha_7^3 - \alpha_7^4 - 
       210 \alpha_0^3 \alpha_8 - 
       432 \alpha_0^2 \alpha_1 \alpha_8 - 
       120 \alpha_0 \alpha_1^2 \alpha_8 - 
       12 \alpha_1^3 \alpha_8 - 
       324 \alpha_0^2 \alpha_2 \alpha_8 - 
       180 \alpha_0 \alpha_1 \alpha_2 \alpha_8 - 
       27 \alpha_1^2 \alpha_2 \alpha_8 - 
       72 \alpha_0 \alpha_2^2 \alpha_8 - 
       18 \alpha_1 \alpha_2^2 \alpha_8 - 3 \alpha_2^3 \alpha_8 - 
       216 \alpha_0^2 \alpha_3 \alpha_8 - 
       120 \alpha_0 \alpha_1 \alpha_3 \alpha_8 - 
       18 \alpha_1^2 \alpha_3 \alpha_8 - 
       96 \alpha_0 \alpha_2 \alpha_3 \alpha_8 - 
       24 \alpha_1 \alpha_2 \alpha_3 \alpha_8 - 
       6 \alpha_2^2 \alpha_3 \alpha_8 - 
       36 \alpha_0 \alpha_3^2 \alpha_8 - 
       6 \alpha_1 \alpha_3^2 \alpha_8 - 
       3 \alpha_2 \alpha_3^2 \alpha_8 - 
       108 \alpha_0^2 \alpha_4 \alpha_8 - 
       60 \alpha_0 \alpha_1 \alpha_4 \alpha_8 - 
       9 \alpha_1^2 \alpha_4 \alpha_8 - 
       48 \alpha_0 \alpha_2 \alpha_4 \alpha_8 - 
       12 \alpha_1 \alpha_2 \alpha_4 \alpha_8 - 
       3 \alpha_2^2 \alpha_4 \alpha_8 - 
       36 \alpha_0 \alpha_3 \alpha_4 \alpha_8 - 
       6 \alpha_1 \alpha_3 \alpha_4 \alpha_8 - 
       3 \alpha_2 \alpha_3 \alpha_4 \alpha_8 - 
       12 \alpha_0 \alpha_4^2 \alpha_8 - 
       360 \alpha_0^2 \alpha_6 \alpha_8 - 
       384 \alpha_0 \alpha_1 \alpha_6 \alpha_8 - 
       76 \alpha_1^2 \alpha_6 \alpha_8 - 
       288 \alpha_0 \alpha_2 \alpha_6 \alpha_8 - 
       114 \alpha_1 \alpha_2 \alpha_6 \alpha_8 - 
       42 \alpha_2^2 \alpha_6 \alpha_8 - 
       192 \alpha_0 \alpha_3 \alpha_6 \alpha_8 - 
       76 \alpha_1 \alpha_3 \alpha_6 \alpha_8 - 
       56 \alpha_2 \alpha_3 \alpha_6 \alpha_8 - 
       18 \alpha_3^2 \alpha_6 \alpha_8 - 
       96 \alpha_0 \alpha_4 \alpha_6 \alpha_8 - 
       38 \alpha_1 \alpha_4 \alpha_6 \alpha_8 - 
       28 \alpha_2 \alpha_4 \alpha_6 \alpha_8 - 
       18 \alpha_3 \alpha_4 \alpha_6 \alpha_8 - 
       4 \alpha_4^2 \alpha_6 \alpha_8 - 
       197 \alpha_0 \alpha_6^2 \alpha_8 - 
       100 \alpha_1 \alpha_6^2 \alpha_8 - 
       75 \alpha_2 \alpha_6^2 \alpha_8 - 
       50 \alpha_3 \alpha_6^2 \alpha_8 - 
       25 \alpha_4 \alpha_6^2 \alpha_8 - 
       36 \alpha_6^3 \alpha_8 - 
       180 \alpha_0^2 \alpha_7 \alpha_8 - 
       192 \alpha_0 \alpha_1 \alpha_7 \alpha_8 - 
       38 \alpha_1^2 \alpha_7 \alpha_8 - 
       144 \alpha_0 \alpha_2 \alpha_7 \alpha_8 - 
       57 \alpha_1 \alpha_2 \alpha_7 \alpha_8 - 
       21 \alpha_2^2 \alpha_7 \alpha_8 - 
       96 \alpha_0 \alpha_3 \alpha_7 \alpha_8 - 
       38 \alpha_1 \alpha_3 \alpha_7 \alpha_8 - 
       28 \alpha_2 \alpha_3 \alpha_7 \alpha_8 - 
       9 \alpha_3^2 \alpha_7 \alpha_8 - 
       48 \alpha_0 \alpha_4 \alpha_7 \alpha_8 - 
       19 \alpha_1 \alpha_4 \alpha_7 \alpha_8 - 
       14 \alpha_2 \alpha_4 \alpha_7 \alpha_8 - 
       9 \alpha_3 \alpha_4 \alpha_7 \alpha_8 - 
       2 \alpha_4^2 \alpha_7 \alpha_8 - 
       197 \alpha_0 \alpha_6 \alpha_7 \alpha_8 - 
       100 \alpha_1 \alpha_6 \alpha_7 \alpha_8 - 
       75 \alpha_2 \alpha_6 \alpha_7 \alpha_8 - 
       50 \alpha_3 \alpha_6 \alpha_7 \alpha_8 - 
       25 \alpha_4 \alpha_6 \alpha_7 \alpha_8 - 
       54 \alpha_6^2 \alpha_7 \alpha_8 - 
       47 \alpha_0 \alpha_7^2 \alpha_8 - 
       28 \alpha_1 \alpha_7^2 \alpha_8 - 
       21 \alpha_2 \alpha_7^2 \alpha_8 - 
       14 \alpha_3 \alpha_7^2 \alpha_8 - 
       7 \alpha_4 \alpha_7^2 \alpha_8 - 
       28 \alpha_6 \alpha_7^2 \alpha_8 - 5 \alpha_7^3 \alpha_8 - 
       156 \alpha_0^2 \alpha_8^2 - 
       192 \alpha_0 \alpha_1 \alpha_8^2 - 
       20 \alpha_1^2 \alpha_8^2 - 
       144 \alpha_0 \alpha_2 \alpha_8^2 - 
       30 \alpha_1 \alpha_2 \alpha_8^2 - 
       12 \alpha_2^2 \alpha_8^2 - 
       96 \alpha_0 \alpha_3 \alpha_8^2 - 
       20 \alpha_1 \alpha_3 \alpha_8^2 - 
       16 \alpha_2 \alpha_3 \alpha_8^2 - 
       6 \alpha_3^2 \alpha_8^2 - 
       48 \alpha_0 \alpha_4 \alpha_8^2 - 
       10 \alpha_1 \alpha_4 \alpha_8^2 - 
       8 \alpha_2 \alpha_4 \alpha_8^2 - 
       6 \alpha_3 \alpha_4 \alpha_8^2 - 
       2 \alpha_4^2 \alpha_8^2 - 
       172 \alpha_0 \alpha_6 \alpha_8^2 - 
       80 \alpha_1 \alpha_6 \alpha_8^2 - 
       60 \alpha_2 \alpha_6 \alpha_8^2 - 
       40 \alpha_3 \alpha_6 \alpha_8^2 - 
       20 \alpha_4 \alpha_6 \alpha_8^2 - 
       45 \alpha_6^2 \alpha_8^2 - 
       86 \alpha_0 \alpha_7 \alpha_8^2 - 
       40 \alpha_1 \alpha_7 \alpha_8^2 - 
       30 \alpha_2 \alpha_7 \alpha_8^2 - 
       20 \alpha_3 \alpha_7 \alpha_8^2 - 
       10 \alpha_4 \alpha_7 \alpha_8^2 - 
       45 \alpha_6 \alpha_7 \alpha_8^2 - 
       11 \alpha_7^2 \alpha_8^2 - 51 \alpha_0 \alpha_8^3 - 
       24 \alpha_1 \alpha_8^3 - 18 \alpha_2 \alpha_8^3 - 
       12 \alpha_3 \alpha_8^3 - 6 \alpha_4 \alpha_8^3 - 
       26 \alpha_6 \alpha_8^3 - 13 \alpha_7 \alpha_8^3 - 
       6 \alpha_8^4) + 
    q^3 (30 \alpha_0^4 + 144 \alpha_0^3 \alpha_1 + 
       60 \alpha_0^2 \alpha_1^2 + 108 \alpha_0^3 \alpha_2 + 
       90 \alpha_0^2 \alpha_1 \alpha_2 + 
       36 \alpha_0^2 \alpha_2^2 + 72 \alpha_0^3 \alpha_3 + 
       60 \alpha_0^2 \alpha_1 \alpha_3 + 
       48 \alpha_0^2 \alpha_2 \alpha_3 + 
       18 \alpha_0^2 \alpha_3^2 + 36 \alpha_0^3 \alpha_4 + 
       30 \alpha_0^2 \alpha_1 \alpha_4 + 
       24 \alpha_0^2 \alpha_2 \alpha_4 + 
       18 \alpha_0^2 \alpha_3 \alpha_4 + 
       6 \alpha_0^2 \alpha_4^2 + 80 \alpha_0^3 \alpha_6 + 
       96 \alpha_0^2 \alpha_1 \alpha_6 + 
       72 \alpha_0^2 \alpha_2 \alpha_6 + 
       48 \alpha_0^2 \alpha_3 \alpha_6 + 
       24 \alpha_0^2 \alpha_4 \alpha_6 + 
       36 \alpha_0^2 \alpha_6^2 + 40 \alpha_0^3 \alpha_7 + 
       48 \alpha_0^2 \alpha_1 \alpha_7 + 
       36 \alpha_0^2 \alpha_2 \alpha_7 + 
       24 \alpha_0^2 \alpha_3 \alpha_7 + 
       12 \alpha_0^2 \alpha_4 \alpha_7 + 
       36 \alpha_0^2 \alpha_6 \alpha_7 + 
       6 \alpha_0^2 \alpha_7^2 + 60 \alpha_0^3 \alpha_8 + 
       216 \alpha_0^2 \alpha_1 \alpha_8 + 
       60 \alpha_0 \alpha_1^2 \alpha_8 + 
       162 \alpha_0^2 \alpha_2 \alpha_8 + 
       90 \alpha_0 \alpha_1 \alpha_2 \alpha_8 + 
       36 \alpha_0 \alpha_2^2 \alpha_8 + 
       108 \alpha_0^2 \alpha_3 \alpha_8 + 
       60 \alpha_0 \alpha_1 \alpha_3 \alpha_8 + 
       48 \alpha_0 \alpha_2 \alpha_3 \alpha_8 + 
       18 \alpha_0 \alpha_3^2 \alpha_8 + 
       54 \alpha_0^2 \alpha_4 \alpha_8 + 
       30 \alpha_0 \alpha_1 \alpha_4 \alpha_8 + 
       24 \alpha_0 \alpha_2 \alpha_4 \alpha_8 + 
       18 \alpha_0 \alpha_3 \alpha_4 \alpha_8 + 
       6 \alpha_0 \alpha_4^2 \alpha_8 + 
       120 \alpha_0^2 \alpha_6 \alpha_8 + 
       96 \alpha_0 \alpha_1 \alpha_6 \alpha_8 + 
       72 \alpha_0 \alpha_2 \alpha_6 \alpha_8 + 
       48 \alpha_0 \alpha_3 \alpha_6 \alpha_8 + 
       24 \alpha_0 \alpha_4 \alpha_6 \alpha_8 + 
       36 \alpha_0 \alpha_6^2 \alpha_8 + 
       60 \alpha_0^2 \alpha_7 \alpha_8 + 
       48 \alpha_0 \alpha_1 \alpha_7 \alpha_8 + 
       36 \alpha_0 \alpha_2 \alpha_7 \alpha_8 + 
       24 \alpha_0 \alpha_3 \alpha_7 \alpha_8 + 
       12 \alpha_0 \alpha_4 \alpha_7 \alpha_8 + 
       36 \alpha_0 \alpha_6 \alpha_7 \alpha_8 + 
       6 \alpha_0 \alpha_7^2 \alpha_8 + 
       51 \alpha_0^2 \alpha_8^2 + 
       96 \alpha_0 \alpha_1 \alpha_8^2 + 
       10 \alpha_1^2 \alpha_8^2 + 
       72 \alpha_0 \alpha_2 \alpha_8^2 + 
       15 \alpha_1 \alpha_2 \alpha_8^2 + 
       6 \alpha_2^2 \alpha_8^2 + 
       48 \alpha_0 \alpha_3 \alpha_8^2 + 
       10 \alpha_1 \alpha_3 \alpha_8^2 + 
       8 \alpha_2 \alpha_3 \alpha_8^2 + 
       3 \alpha_3^2 \alpha_8^2 + 
       24 \alpha_0 \alpha_4 \alpha_8^2 + 
       5 \alpha_1 \alpha_4 \alpha_8^2 + 
       4 \alpha_2 \alpha_4 \alpha_8^2 + 
       3 \alpha_3 \alpha_4 \alpha_8^2 + \alpha_4^2 \alpha_8^2 + 
       56 \alpha_0 \alpha_6 \alpha_8^2 + 
       16 \alpha_1 \alpha_6 \alpha_8^2 + 
       12 \alpha_2 \alpha_6 \alpha_8^2 + 
       8 \alpha_3 \alpha_6 \alpha_8^2 + 
       4 \alpha_4 \alpha_6 \alpha_8^2 + 
       6 \alpha_6^2 \alpha_8^2 + 
       28 \alpha_0 \alpha_7 \alpha_8^2 + 
       8 \alpha_1 \alpha_7 \alpha_8^2 + 
       6 \alpha_2 \alpha_7 \alpha_8^2 + 
       4 \alpha_3 \alpha_7 \alpha_8^2 + 
       2 \alpha_4 \alpha_7 \alpha_8^2 + 
       6 \alpha_6 \alpha_7 \alpha_8^2 + \alpha_7^2 \alpha_8^2 + 
       21 \alpha_0 \alpha_8^3 + 12 \alpha_1 \alpha_8^3 + 
       9 \alpha_2 \alpha_8^3 + 6 \alpha_3 \alpha_8^3 + 
       3 \alpha_4 \alpha_8^3 + 8 \alpha_6 \alpha_8^3 + 
       4 \alpha_7 \alpha_8^3 + 3 \alpha_8^4) + 
    q (60 \alpha_0^4 + 144 \alpha_0^3 \alpha_1 + 
       60 \alpha_0^2 \alpha_1^2 + 48 \alpha_0 \alpha_1^3 + 
       30 \alpha_1^4 + 108 \alpha_0^3 \alpha_2 + 
       90 \alpha_0^2 \alpha_1 \alpha_2 + 
       108 \alpha_0 \alpha_1^2 \alpha_2 + 
       90 \alpha_1^3 \alpha_2 + 36 \alpha_0^2 \alpha_2^2 + 
       72 \alpha_0 \alpha_1 \alpha_2^2 + 
       96 \alpha_1^2 \alpha_2^2 + 12 \alpha_0 \alpha_2^3 + 
       42 \alpha_1 \alpha_2^3 + 6 \alpha_2^4 + 
       72 \alpha_0^3 \alpha_3 + 
       60 \alpha_0^2 \alpha_1 \alpha_3 + 
       72 \alpha_0 \alpha_1^2 \alpha_3 + 
       60 \alpha_1^3 \alpha_3 + 
       48 \alpha_0^2 \alpha_2 \alpha_3 + 
       96 \alpha_0 \alpha_1 \alpha_2 \alpha_3 + 
       128 \alpha_1^2 \alpha_2 \alpha_3 + 
       24 \alpha_0 \alpha_2^2 \alpha_3 + 
       84 \alpha_1 \alpha_2^2 \alpha_3 + 
       16 \alpha_2^3 \alpha_3 + 18 \alpha_0^2 \alpha_3^2 + 
       24 \alpha_0 \alpha_1 \alpha_3^2 + 
       38 \alpha_1^2 \alpha_3^2 + 
       12 \alpha_0 \alpha_2 \alpha_3^2 + 
       50 \alpha_1 \alpha_2 \alpha_3^2 + 
       14 \alpha_2^2 \alpha_3^2 + 8 \alpha_1 \alpha_3^3 + 
       4 \alpha_2 \alpha_3^3 + 36 \alpha_0^3 \alpha_4 + 
       30 \alpha_0^2 \alpha_1 \alpha_4 + 
       36 \alpha_0 \alpha_1^2 \alpha_4 + 
       30 \alpha_1^3 \alpha_4 + 
       24 \alpha_0^2 \alpha_2 \alpha_4 + 
       48 \alpha_0 \alpha_1 \alpha_2 \alpha_4 + 
       64 \alpha_1^2 \alpha_2 \alpha_4 + 
       12 \alpha_0 \alpha_2^2 \alpha_4 + 
       42 \alpha_1 \alpha_2^2 \alpha_4 + 8 \alpha_2^3 \alpha_4 + 
       18 \alpha_0^2 \alpha_3 \alpha_4 + 
       24 \alpha_0 \alpha_1 \alpha_3 \alpha_4 + 
       38 \alpha_1^2 \alpha_3 \alpha_4 + 
       12 \alpha_0 \alpha_2 \alpha_3 \alpha_4 + 
       50 \alpha_1 \alpha_2 \alpha_3 \alpha_4 + 
       14 \alpha_2^2 \alpha_3 \alpha_4 + 
       12 \alpha_1 \alpha_3^2 \alpha_4 + 
       6 \alpha_2 \alpha_3^2 \alpha_4 + 
       6 \alpha_0^2 \alpha_4^2 + 6 \alpha_1^2 \alpha_4^2 + 
       8 \alpha_1 \alpha_2 \alpha_4^2 + 
       2 \alpha_2^2 \alpha_4^2 + 
       4 \alpha_1 \alpha_3 \alpha_4^2 + 
       2 \alpha_2 \alpha_3 \alpha_4^2 + 
       160 \alpha_0^3 \alpha_6 + 
       288 \alpha_0^2 \alpha_1 \alpha_6 + 
       152 \alpha_0 \alpha_1^2 \alpha_6 + 
       72 \alpha_1^3 \alpha_6 + 
       216 \alpha_0^2 \alpha_2 \alpha_6 + 
       228 \alpha_0 \alpha_1 \alpha_2 \alpha_6 + 
       162 \alpha_1^2 \alpha_2 \alpha_6 + 
       84 \alpha_0 \alpha_2^2 \alpha_6 + 
       114 \alpha_1 \alpha_2^2 \alpha_6 + 
       24 \alpha_2^3 \alpha_6 + 
       144 \alpha_0^2 \alpha_3 \alpha_6 + 
       152 \alpha_0 \alpha_1 \alpha_3 \alpha_6 + 
       108 \alpha_1^2 \alpha_3 \alpha_6 + 
       112 \alpha_0 \alpha_2 \alpha_3 \alpha_6 + 
       152 \alpha_1 \alpha_2 \alpha_3 \alpha_6 + 
       48 \alpha_2^2 \alpha_3 \alpha_6 + 
       36 \alpha_0 \alpha_3^2 \alpha_6 + 
       44 \alpha_1 \alpha_3^2 \alpha_6 + 
       28 \alpha_2 \alpha_3^2 \alpha_6 + 4 \alpha_3^3 \alpha_6 + 
       72 \alpha_0^2 \alpha_4 \alpha_6 + 
       76 \alpha_0 \alpha_1 \alpha_4 \alpha_6 + 
       54 \alpha_1^2 \alpha_4 \alpha_6 + 
       56 \alpha_0 \alpha_2 \alpha_4 \alpha_6 + 
       76 \alpha_1 \alpha_2 \alpha_4 \alpha_6 + 
       24 \alpha_2^2 \alpha_4 \alpha_6 + 
       36 \alpha_0 \alpha_3 \alpha_4 \alpha_6 + 
       44 \alpha_1 \alpha_3 \alpha_4 \alpha_6 + 
       28 \alpha_2 \alpha_3 \alpha_4 \alpha_6 + 
       6 \alpha_3^2 \alpha_4 \alpha_6 + 
       8 \alpha_0 \alpha_4^2 \alpha_6 + 
       6 \alpha_1 \alpha_4^2 \alpha_6 + 
       4 \alpha_2 \alpha_4^2 \alpha_6 + 
       2 \alpha_3 \alpha_4^2 \alpha_6 + 
       161 \alpha_0^2 \alpha_6^2 + 
       200 \alpha_0 \alpha_1 \alpha_6^2 + 
       78 \alpha_1^2 \alpha_6^2 + 
       150 \alpha_0 \alpha_2 \alpha_6^2 + 
       117 \alpha_1 \alpha_2 \alpha_6^2 + 
       42 \alpha_2^2 \alpha_6^2 + 
       100 \alpha_0 \alpha_3 \alpha_6^2 + 
       78 \alpha_1 \alpha_3 \alpha_6^2 + 
       56 \alpha_2 \alpha_3 \alpha_6^2 + 
       17 \alpha_3^2 \alpha_6^2 + 
       50 \alpha_0 \alpha_4 \alpha_6^2 + 
       39 \alpha_1 \alpha_4 \alpha_6^2 + 
       28 \alpha_2 \alpha_4 \alpha_6^2 + 
       17 \alpha_3 \alpha_4 \alpha_6^2 + 
       3 \alpha_4^2 \alpha_6^2 + 72 \alpha_0 \alpha_6^3 + 
       48 \alpha_1 \alpha_6^3 + 36 \alpha_2 \alpha_6^3 + 
       24 \alpha_3 \alpha_6^3 + 12 \alpha_4 \alpha_6^3 + 
       12 \alpha_6^4 + 80 \alpha_0^3 \alpha_7 + 
       144 \alpha_0^2 \alpha_1 \alpha_7 + 
       76 \alpha_0 \alpha_1^2 \alpha_7 + 
       36 \alpha_1^3 \alpha_7 + 
       108 \alpha_0^2 \alpha_2 \alpha_7 + 
       114 \alpha_0 \alpha_1 \alpha_2 \alpha_7 + 
       81 \alpha_1^2 \alpha_2 \alpha_7 + 
       42 \alpha_0 \alpha_2^2 \alpha_7 + 
       57 \alpha_1 \alpha_2^2 \alpha_7 + 
       12 \alpha_2^3 \alpha_7 + 
       72 \alpha_0^2 \alpha_3 \alpha_7 + 
       76 \alpha_0 \alpha_1 \alpha_3 \alpha_7 + 
       54 \alpha_1^2 \alpha_3 \alpha_7 + 
       56 \alpha_0 \alpha_2 \alpha_3 \alpha_7 + 
       76 \alpha_1 \alpha_2 \alpha_3 \alpha_7 + 
       24 \alpha_2^2 \alpha_3 \alpha_7 + 
       18 \alpha_0 \alpha_3^2 \alpha_7 + 
       22 \alpha_1 \alpha_3^2 \alpha_7 + 
       14 \alpha_2 \alpha_3^2 \alpha_7 + 2 \alpha_3^3 \alpha_7 + 
       36 \alpha_0^2 \alpha_4 \alpha_7 + 
       38 \alpha_0 \alpha_1 \alpha_4 \alpha_7 + 
       27 \alpha_1^2 \alpha_4 \alpha_7 + 
       28 \alpha_0 \alpha_2 \alpha_4 \alpha_7 + 
       38 \alpha_1 \alpha_2 \alpha_4 \alpha_7 + 
       12 \alpha_2^2 \alpha_4 \alpha_7 + 
       18 \alpha_0 \alpha_3 \alpha_4 \alpha_7 + 
       22 \alpha_1 \alpha_3 \alpha_4 \alpha_7 + 
       14 \alpha_2 \alpha_3 \alpha_4 \alpha_7 + 
       3 \alpha_3^2 \alpha_4 \alpha_7 + 
       4 \alpha_0 \alpha_4^2 \alpha_7 + 
       3 \alpha_1 \alpha_4^2 \alpha_7 + 
       2 \alpha_2 \alpha_4^2 \alpha_7 + \alpha_3 \alpha_4^2 \
\alpha_7 + 161 \alpha_0^2 \alpha_6 \alpha_7 + 
       200 \alpha_0 \alpha_1 \alpha_6 \alpha_7 + 
       78 \alpha_1^2 \alpha_6 \alpha_7 + 
       150 \alpha_0 \alpha_2 \alpha_6 \alpha_7 + 
       117 \alpha_1 \alpha_2 \alpha_6 \alpha_7 + 
       42 \alpha_2^2 \alpha_6 \alpha_7 + 
       100 \alpha_0 \alpha_3 \alpha_6 \alpha_7 + 
       78 \alpha_1 \alpha_3 \alpha_6 \alpha_7 + 
       56 \alpha_2 \alpha_3 \alpha_6 \alpha_7 + 
       17 \alpha_3^2 \alpha_6 \alpha_7 + 
       50 \alpha_0 \alpha_4 \alpha_6 \alpha_7 + 
       39 \alpha_1 \alpha_4 \alpha_6 \alpha_7 + 
       28 \alpha_2 \alpha_4 \alpha_6 \alpha_7 + 
       17 \alpha_3 \alpha_4 \alpha_6 \alpha_7 + 
       3 \alpha_4^2 \alpha_6 \alpha_7 + 
       108 \alpha_0 \alpha_6^2 \alpha_7 + 
       72 \alpha_1 \alpha_6^2 \alpha_7 + 
       54 \alpha_2 \alpha_6^2 \alpha_7 + 
       36 \alpha_3 \alpha_6^2 \alpha_7 + 
       18 \alpha_4 \alpha_6^2 \alpha_7 + 
       24 \alpha_6^3 \alpha_7 + 41 \alpha_0^2 \alpha_7^2 + 
       56 \alpha_0 \alpha_1 \alpha_7^2 + 
       22 \alpha_1^2 \alpha_7^2 + 
       42 \alpha_0 \alpha_2 \alpha_7^2 + 
       33 \alpha_1 \alpha_2 \alpha_7^2 + 
       12 \alpha_2^2 \alpha_7^2 + 
       28 \alpha_0 \alpha_3 \alpha_7^2 + 
       22 \alpha_1 \alpha_3 \alpha_7^2 + 
       16 \alpha_2 \alpha_3 \alpha_7^2 + 
       5 \alpha_3^2 \alpha_7^2 + 
       14 \alpha_0 \alpha_4 \alpha_7^2 + 
       11 \alpha_1 \alpha_4 \alpha_7^2 + 
       8 \alpha_2 \alpha_4 \alpha_7^2 + 
       5 \alpha_3 \alpha_4 \alpha_7^2 + \alpha_4^2 \alpha_7^2 + 
       56 \alpha_0 \alpha_6 \alpha_7^2 + 
       40 \alpha_1 \alpha_6 \alpha_7^2 + 
       30 \alpha_2 \alpha_6 \alpha_7^2 + 
       20 \alpha_3 \alpha_6 \alpha_7^2 + 
       10 \alpha_4 \alpha_6 \alpha_7^2 + 
       19 \alpha_6^2 \alpha_7^2 + 10 \alpha_0 \alpha_7^3 + 
       8 \alpha_1 \alpha_7^3 + 6 \alpha_2 \alpha_7^3 + 
       4 \alpha_3 \alpha_7^3 + 2 \alpha_4 \alpha_7^3 + 
       7 \alpha_6 \alpha_7^3 + \alpha_7^4 + 
       120 \alpha_0^3 \alpha_8 + 
       216 \alpha_0^2 \alpha_1 \alpha_8 + 
       60 \alpha_0 \alpha_1^2 \alpha_8 + 
       24 \alpha_1^3 \alpha_8 + 
       162 \alpha_0^2 \alpha_2 \alpha_8 + 
       90 \alpha_0 \alpha_1 \alpha_2 \alpha_8 + 
       54 \alpha_1^2 \alpha_2 \alpha_8 + 
       36 \alpha_0 \alpha_2^2 \alpha_8 + 
       36 \alpha_1 \alpha_2^2 \alpha_8 + 6 \alpha_2^3 \alpha_8 + 
       108 \alpha_0^2 \alpha_3 \alpha_8 + 
       60 \alpha_0 \alpha_1 \alpha_3 \alpha_8 + 
       36 \alpha_1^2 \alpha_3 \alpha_8 + 
       48 \alpha_0 \alpha_2 \alpha_3 \alpha_8 + 
       48 \alpha_1 \alpha_2 \alpha_3 \alpha_8 + 
       12 \alpha_2^2 \alpha_3 \alpha_8 + 
       18 \alpha_0 \alpha_3^2 \alpha_8 + 
       12 \alpha_1 \alpha_3^2 \alpha_8 + 
       6 \alpha_2 \alpha_3^2 \alpha_8 + 
       54 \alpha_0^2 \alpha_4 \alpha_8 + 
       30 \alpha_0 \alpha_1 \alpha_4 \alpha_8 + 
       18 \alpha_1^2 \alpha_4 \alpha_8 + 
       24 \alpha_0 \alpha_2 \alpha_4 \alpha_8 + 
       24 \alpha_1 \alpha_2 \alpha_4 \alpha_8 + 
       6 \alpha_2^2 \alpha_4 \alpha_8 + 
       18 \alpha_0 \alpha_3 \alpha_4 \alpha_8 + 
       12 \alpha_1 \alpha_3 \alpha_4 \alpha_8 + 
       6 \alpha_2 \alpha_3 \alpha_4 \alpha_8 + 
       6 \alpha_0 \alpha_4^2 \alpha_8 + 
       240 \alpha_0^2 \alpha_6 \alpha_8 + 
       288 \alpha_0 \alpha_1 \alpha_6 \alpha_8 + 
       76 \alpha_1^2 \alpha_6 \alpha_8 + 
       216 \alpha_0 \alpha_2 \alpha_6 \alpha_8 + 
       114 \alpha_1 \alpha_2 \alpha_6 \alpha_8 + 
       42 \alpha_2^2 \alpha_6 \alpha_8 + 
       144 \alpha_0 \alpha_3 \alpha_6 \alpha_8 + 
       76 \alpha_1 \alpha_3 \alpha_6 \alpha_8 + 
       56 \alpha_2 \alpha_3 \alpha_6 \alpha_8 + 
       18 \alpha_3^2 \alpha_6 \alpha_8 + 
       72 \alpha_0 \alpha_4 \alpha_6 \alpha_8 + 
       38 \alpha_1 \alpha_4 \alpha_6 \alpha_8 + 
       28 \alpha_2 \alpha_4 \alpha_6 \alpha_8 + 
       18 \alpha_3 \alpha_4 \alpha_6 \alpha_8 + 
       4 \alpha_4^2 \alpha_6 \alpha_8 + 
       161 \alpha_0 \alpha_6^2 \alpha_8 + 
       100 \alpha_1 \alpha_6^2 \alpha_8 + 
       75 \alpha_2 \alpha_6^2 \alpha_8 + 
       50 \alpha_3 \alpha_6^2 \alpha_8 + 
       25 \alpha_4 \alpha_6^2 \alpha_8 + 
       36 \alpha_6^3 \alpha_8 + 
       120 \alpha_0^2 \alpha_7 \alpha_8 + 
       144 \alpha_0 \alpha_1 \alpha_7 \alpha_8 + 
       38 \alpha_1^2 \alpha_7 \alpha_8 + 
       108 \alpha_0 \alpha_2 \alpha_7 \alpha_8 + 
       57 \alpha_1 \alpha_2 \alpha_7 \alpha_8 + 
       21 \alpha_2^2 \alpha_7 \alpha_8 + 
       72 \alpha_0 \alpha_3 \alpha_7 \alpha_8 + 
       38 \alpha_1 \alpha_3 \alpha_7 \alpha_8 + 
       28 \alpha_2 \alpha_3 \alpha_7 \alpha_8 + 
       9 \alpha_3^2 \alpha_7 \alpha_8 + 
       36 \alpha_0 \alpha_4 \alpha_7 \alpha_8 + 
       19 \alpha_1 \alpha_4 \alpha_7 \alpha_8 + 
       14 \alpha_2 \alpha_4 \alpha_7 \alpha_8 + 
       9 \alpha_3 \alpha_4 \alpha_7 \alpha_8 + 
       2 \alpha_4^2 \alpha_7 \alpha_8 + 
       161 \alpha_0 \alpha_6 \alpha_7 \alpha_8 + 
       100 \alpha_1 \alpha_6 \alpha_7 \alpha_8 + 
       75 \alpha_2 \alpha_6 \alpha_7 \alpha_8 + 
       50 \alpha_3 \alpha_6 \alpha_7 \alpha_8 + 
       25 \alpha_4 \alpha_6 \alpha_7 \alpha_8 + 
       54 \alpha_6^2 \alpha_7 \alpha_8 + 
       41 \alpha_0 \alpha_7^2 \alpha_8 + 
       28 \alpha_1 \alpha_7^2 \alpha_8 + 
       21 \alpha_2 \alpha_7^2 \alpha_8 + 
       14 \alpha_3 \alpha_7^2 \alpha_8 + 
       7 \alpha_4 \alpha_7^2 \alpha_8 + 
       28 \alpha_6 \alpha_7^2 \alpha_8 + 5 \alpha_7^3 \alpha_8 + 
       87 \alpha_0^2 \alpha_8^2 + 
       96 \alpha_0 \alpha_1 \alpha_8^2 + 
       10 \alpha_1^2 \alpha_8^2 + 
       72 \alpha_0 \alpha_2 \alpha_8^2 + 
       15 \alpha_1 \alpha_2 \alpha_8^2 + 
       6 \alpha_2^2 \alpha_8^2 + 
       48 \alpha_0 \alpha_3 \alpha_8^2 + 
       10 \alpha_1 \alpha_3 \alpha_8^2 + 
       8 \alpha_2 \alpha_3 \alpha_8^2 + 
       3 \alpha_3^2 \alpha_8^2 + 
       24 \alpha_0 \alpha_4 \alpha_8^2 + 
       5 \alpha_1 \alpha_4 \alpha_8^2 + 
       4 \alpha_2 \alpha_4 \alpha_8^2 + 
       3 \alpha_3 \alpha_4 \alpha_8^2 + \alpha_4^2 \alpha_8^2 + 
       116 \alpha_0 \alpha_6 \alpha_8^2 + 
       64 \alpha_1 \alpha_6 \alpha_8^2 + 
       48 \alpha_2 \alpha_6 \alpha_8^2 + 
       32 \alpha_3 \alpha_6 \alpha_8^2 + 
       16 \alpha_4 \alpha_6 \alpha_8^2 + 
       39 \alpha_6^2 \alpha_8^2 + 
       58 \alpha_0 \alpha_7 \alpha_8^2 + 
       32 \alpha_1 \alpha_7 \alpha_8^2 + 
       24 \alpha_2 \alpha_7 \alpha_8^2 + 
       16 \alpha_3 \alpha_7 \alpha_8^2 + 
       8 \alpha_4 \alpha_7 \alpha_8^2 + 
       39 \alpha_6 \alpha_7 \alpha_8^2 + 
       10 \alpha_7^2 \alpha_8^2 + 27 \alpha_0 \alpha_8^3 + 
       12 \alpha_1 \alpha_8^3 + 9 \alpha_2 \alpha_8^3 + 
       6 \alpha_3 \alpha_8^3 + 3 \alpha_4 \alpha_8^3 + 
       18 \alpha_6 \alpha_8^3 + 9 \alpha_7 \alpha_8^3 + 
       3 \alpha_8^4)) + 
 p (-\alpha_1 (\alpha_1 + \alpha_2) (\alpha_1 + \alpha_2 + \
\alpha_3) (\alpha_1 + \alpha_2 + \alpha_3 + \alpha_4) (\alpha_1 \
+ \alpha_2 + \alpha_3 + \alpha_4 + \alpha_5) + 
    3 q^4 \alpha_0^2 (\alpha_0 + \alpha_8)^2 (2 \alpha_0 + \
\alpha_8) + 
    2 q^3 \alpha_0 (\alpha_0 + \alpha_8) (18 \alpha_0^3 + 
       48 \alpha_0^2 \alpha_1 + 20 \alpha_0 \alpha_1^2 + 
       36 \alpha_0^2 \alpha_2 + 30 \alpha_0 \alpha_1 \alpha_2 + 
       12 \alpha_0 \alpha_2^2 + 24 \alpha_0^2 \alpha_3 + 
       20 \alpha_0 \alpha_1 \alpha_3 + 
       16 \alpha_0 \alpha_2 \alpha_3 + 6 \alpha_0 \alpha_3^2 + 
       12 \alpha_0^2 \alpha_4 + 10 \alpha_0 \alpha_1 \alpha_4 + 
       8 \alpha_0 \alpha_2 \alpha_4 + 
       6 \alpha_0 \alpha_3 \alpha_4 + 2 \alpha_0 \alpha_4^2 + 
       30 \alpha_0^2 \alpha_6 + 32 \alpha_0 \alpha_1 \alpha_6 + 
       24 \alpha_0 \alpha_2 \alpha_6 + 
       16 \alpha_0 \alpha_3 \alpha_6 + 
       8 \alpha_0 \alpha_4 \alpha_6 + 12 \alpha_0 \alpha_6^2 + 
       15 \alpha_0^2 \alpha_7 + 16 \alpha_0 \alpha_1 \alpha_7 + 
       12 \alpha_0 \alpha_2 \alpha_7 + 
       8 \alpha_0 \alpha_3 \alpha_7 + 
       4 \alpha_0 \alpha_4 \alpha_7 + 
       12 \alpha_0 \alpha_6 \alpha_7 + 2 \alpha_0 \alpha_7^2 + 
       27 \alpha_0^2 \alpha_8 + 48 \alpha_0 \alpha_1 \alpha_8 + 
       10 \alpha_1^2 \alpha_8 + 36 \alpha_0 \alpha_2 \alpha_8 + 
       15 \alpha_1 \alpha_2 \alpha_8 + 6 \alpha_2^2 \alpha_8 + 
       24 \alpha_0 \alpha_3 \alpha_8 + 
       10 \alpha_1 \alpha_3 \alpha_8 + 
       8 \alpha_2 \alpha_3 \alpha_8 + 3 \alpha_3^2 \alpha_8 + 
       12 \alpha_0 \alpha_4 \alpha_8 + 
       5 \alpha_1 \alpha_4 \alpha_8 + 
       4 \alpha_2 \alpha_4 \alpha_8 + 
       3 \alpha_3 \alpha_4 \alpha_8 + \alpha_4^2 \alpha_8 + 
       30 \alpha_0 \alpha_6 \alpha_8 + 
       16 \alpha_1 \alpha_6 \alpha_8 + 
       12 \alpha_2 \alpha_6 \alpha_8 + 
       8 \alpha_3 \alpha_6 \alpha_8 + 
       4 \alpha_4 \alpha_6 \alpha_8 + 6 \alpha_6^2 \alpha_8 + 
       15 \alpha_0 \alpha_7 \alpha_8 + 
       8 \alpha_1 \alpha_7 \alpha_8 + 
       6 \alpha_2 \alpha_7 \alpha_8 + 
       4 \alpha_3 \alpha_7 \alpha_8 + 
       2 \alpha_4 \alpha_7 \alpha_8 + 
       6 \alpha_6 \alpha_7 \alpha_8 + \alpha_7^2 \alpha_8 + 
       15 \alpha_0 \alpha_8^2 + 12 \alpha_1 \alpha_8^2 + 
       9 \alpha_2 \alpha_8^2 + 6 \alpha_3 \alpha_8^2 + 
       3 \alpha_4 \alpha_8^2 + 8 \alpha_6 \alpha_8^2 + 
       4 \alpha_7 \alpha_8^2 + 3 \alpha_8^3) + 
    q^2 (-114 \alpha_0^5 - 288 \alpha_0^4 \alpha_1 - 
       228 \alpha_0^3 \alpha_1^2 - 108 \alpha_0^2 \alpha_1^3 - 
       30 \alpha_0 \alpha_1^4 - 216 \alpha_0^4 \alpha_2 - 
       342 \alpha_0^3 \alpha_1 \alpha_2 - 
       243 \alpha_0^2 \alpha_1^2 \alpha_2 - 
       90 \alpha_0 \alpha_1^3 \alpha_2 - 
       126 \alpha_0^3 \alpha_2^2 - 
       171 \alpha_0^2 \alpha_1 \alpha_2^2 - 
       96 \alpha_0 \alpha_1^2 \alpha_2^2 - 
       36 \alpha_0^2 \alpha_2^3 - 
       42 \alpha_0 \alpha_1 \alpha_2^3 - 6 \alpha_0 \alpha_2^4 - 
       144 \alpha_0^4 \alpha_3 - 
       228 \alpha_0^3 \alpha_1 \alpha_3 - 
       162 \alpha_0^2 \alpha_1^2 \alpha_3 - 
       60 \alpha_0 \alpha_1^3 \alpha_3 - 
       168 \alpha_0^3 \alpha_2 \alpha_3 - 
       228 \alpha_0^2 \alpha_1 \alpha_2 \alpha_3 - 
       128 \alpha_0 \alpha_1^2 \alpha_2 \alpha_3 - 
       72 \alpha_0^2 \alpha_2^2 \alpha_3 - 
       84 \alpha_0 \alpha_1 \alpha_2^2 \alpha_3 - 
       16 \alpha_0 \alpha_2^3 \alpha_3 - 
       54 \alpha_0^3 \alpha_3^2 - 
       66 \alpha_0^2 \alpha_1 \alpha_3^2 - 
       38 \alpha_0 \alpha_1^2 \alpha_3^2 - 
       42 \alpha_0^2 \alpha_2 \alpha_3^2 - 
       50 \alpha_0 \alpha_1 \alpha_2 \alpha_3^2 - 
       14 \alpha_0 \alpha_2^2 \alpha_3^2 - 
       6 \alpha_0^2 \alpha_3^3 - 
       8 \alpha_0 \alpha_1 \alpha_3^3 - 
       4 \alpha_0 \alpha_2 \alpha_3^3 - 72 \alpha_0^4 \alpha_4 - 
       114 \alpha_0^3 \alpha_1 \alpha_4 - 
       81 \alpha_0^2 \alpha_1^2 \alpha_4 - 
       30 \alpha_0 \alpha_1^3 \alpha_4 - 
       84 \alpha_0^3 \alpha_2 \alpha_4 - 
       114 \alpha_0^2 \alpha_1 \alpha_2 \alpha_4 - 
       64 \alpha_0 \alpha_1^2 \alpha_2 \alpha_4 - 
       36 \alpha_0^2 \alpha_2^2 \alpha_4 - 
       42 \alpha_0 \alpha_1 \alpha_2^2 \alpha_4 - 
       8 \alpha_0 \alpha_2^3 \alpha_4 - 
       54 \alpha_0^3 \alpha_3 \alpha_4 - 
       66 \alpha_0^2 \alpha_1 \alpha_3 \alpha_4 - 
       38 \alpha_0 \alpha_1^2 \alpha_3 \alpha_4 - 
       42 \alpha_0^2 \alpha_2 \alpha_3 \alpha_4 - 
       50 \alpha_0 \alpha_1 \alpha_2 \alpha_3 \alpha_4 - 
       14 \alpha_0 \alpha_2^2 \alpha_3 \alpha_4 - 
       9 \alpha_0^2 \alpha_3^2 \alpha_4 - 
       12 \alpha_0 \alpha_1 \alpha_3^2 \alpha_4 - 
       6 \alpha_0 \alpha_2 \alpha_3^2 \alpha_4 - 
       12 \alpha_0^3 \alpha_4^2 - 
       9 \alpha_0^2 \alpha_1 \alpha_4^2 - 
       6 \alpha_0 \alpha_1^2 \alpha_4^2 - 
       6 \alpha_0^2 \alpha_2 \alpha_4^2 - 
       8 \alpha_0 \alpha_1 \alpha_2 \alpha_4^2 - 
       2 \alpha_0 \alpha_2^2 \alpha_4^2 - 
       3 \alpha_0^2 \alpha_3 \alpha_4^2 - 
       4 \alpha_0 \alpha_1 \alpha_3 \alpha_4^2 - 
       2 \alpha_0 \alpha_2 \alpha_3 \alpha_4^2 - 
       300 \alpha_0^4 \alpha_6 - 
       576 \alpha_0^3 \alpha_1 \alpha_6 - 
       376 \alpha_0^2 \alpha_1^2 \alpha_6 - 
       112 \alpha_0 \alpha_1^3 \alpha_6 - 
       432 \alpha_0^3 \alpha_2 \alpha_6 - 
       564 \alpha_0^2 \alpha_1 \alpha_2 \alpha_6 - 
       252 \alpha_0 \alpha_1^2 \alpha_2 \alpha_6 - 
       204 \alpha_0^2 \alpha_2^2 \alpha_6 - 
       180 \alpha_0 \alpha_1 \alpha_2^2 \alpha_6 - 
       40 \alpha_0 \alpha_2^3 \alpha_6 - 
       288 \alpha_0^3 \alpha_3 \alpha_6 - 
       376 \alpha_0^2 \alpha_1 \alpha_3 \alpha_6 - 
       168 \alpha_0 \alpha_1^2 \alpha_3 \alpha_6 - 
       272 \alpha_0^2 \alpha_2 \alpha_3 \alpha_6 - 
       240 \alpha_0 \alpha_1 \alpha_2 \alpha_3 \alpha_6 - 
       80 \alpha_0 \alpha_2^2 \alpha_3 \alpha_6 - 
       84 \alpha_0^2 \alpha_3^2 \alpha_6 - 
       72 \alpha_0 \alpha_1 \alpha_3^2 \alpha_6 - 
       48 \alpha_0 \alpha_2 \alpha_3^2 \alpha_6 - 
       8 \alpha_0 \alpha_3^3 \alpha_6 - 
       144 \alpha_0^3 \alpha_4 \alpha_6 - 
       188 \alpha_0^2 \alpha_1 \alpha_4 \alpha_6 - 
       84 \alpha_0 \alpha_1^2 \alpha_4 \alpha_6 - 
       136 \alpha_0^2 \alpha_2 \alpha_4 \alpha_6 - 
       120 \alpha_0 \alpha_1 \alpha_2 \alpha_4 \alpha_6 - 
       40 \alpha_0 \alpha_2^2 \alpha_4 \alpha_6 - 
       84 \alpha_0^2 \alpha_3 \alpha_4 \alpha_6 - 
       72 \alpha_0 \alpha_1 \alpha_3 \alpha_4 \alpha_6 - 
       48 \alpha_0 \alpha_2 \alpha_3 \alpha_4 \alpha_6 - 
       12 \alpha_0 \alpha_3^2 \alpha_4 \alpha_6 - 
       16 \alpha_0^2 \alpha_4^2 \alpha_6 - 
       12 \alpha_0 \alpha_1 \alpha_4^2 \alpha_6 - 
       8 \alpha_0 \alpha_2 \alpha_4^2 \alpha_6 - 
       4 \alpha_0 \alpha_3 \alpha_4^2 \alpha_6 - 
       298 \alpha_0^3 \alpha_6^2 - 
       400 \alpha_0^2 \alpha_1 \alpha_6^2 - 
       156 \alpha_0 \alpha_1^2 \alpha_6^2 - 
       300 \alpha_0^2 \alpha_2 \alpha_6^2 - 
       234 \alpha_0 \alpha_1 \alpha_2 \alpha_6^2 - 
       84 \alpha_0 \alpha_2^2 \alpha_6^2 - 
       200 \alpha_0^2 \alpha_3 \alpha_6^2 - 
       156 \alpha_0 \alpha_1 \alpha_3 \alpha_6^2 - 
       112 \alpha_0 \alpha_2 \alpha_3 \alpha_6^2 - 
       34 \alpha_0 \alpha_3^2 \alpha_6^2 - 
       100 \alpha_0^2 \alpha_4 \alpha_6^2 - 
       78 \alpha_0 \alpha_1 \alpha_4 \alpha_6^2 - 
       56 \alpha_0 \alpha_2 \alpha_4 \alpha_6^2 - 
       34 \alpha_0 \alpha_3 \alpha_4 \alpha_6^2 - 
       6 \alpha_0 \alpha_4^2 \alpha_6^2 - 
       132 \alpha_0^2 \alpha_6^3 - 
       96 \alpha_0 \alpha_1 \alpha_6^3 - 
       72 \alpha_0 \alpha_2 \alpha_6^3 - 
       48 \alpha_0 \alpha_3 \alpha_6^3 - 
       24 \alpha_0 \alpha_4 \alpha_6^3 - 
       22 \alpha_0 \alpha_6^4 - 150 \alpha_0^4 \alpha_7 - 
       288 \alpha_0^3 \alpha_1 \alpha_7 - 
       188 \alpha_0^2 \alpha_1^2 \alpha_7 - 
       56 \alpha_0 \alpha_1^3 \alpha_7 - 
       216 \alpha_0^3 \alpha_2 \alpha_7 - 
       282 \alpha_0^2 \alpha_1 \alpha_2 \alpha_7 - 
       126 \alpha_0 \alpha_1^2 \alpha_2 \alpha_7 - 
       102 \alpha_0^2 \alpha_2^2 \alpha_7 - 
       90 \alpha_0 \alpha_1 \alpha_2^2 \alpha_7 - 
       20 \alpha_0 \alpha_2^3 \alpha_7 - 
       144 \alpha_0^3 \alpha_3 \alpha_7 - 
       188 \alpha_0^2 \alpha_1 \alpha_3 \alpha_7 - 
       84 \alpha_0 \alpha_1^2 \alpha_3 \alpha_7 - 
       136 \alpha_0^2 \alpha_2 \alpha_3 \alpha_7 - 
       120 \alpha_0 \alpha_1 \alpha_2 \alpha_3 \alpha_7 - 
       40 \alpha_0 \alpha_2^2 \alpha_3 \alpha_7 - 
       42 \alpha_0^2 \alpha_3^2 \alpha_7 - 
       36 \alpha_0 \alpha_1 \alpha_3^2 \alpha_7 - 
       24 \alpha_0 \alpha_2 \alpha_3^2 \alpha_7 - 
       4 \alpha_0 \alpha_3^3 \alpha_7 - 
       72 \alpha_0^3 \alpha_4 \alpha_7 - 
       94 \alpha_0^2 \alpha_1 \alpha_4 \alpha_7 - 
       42 \alpha_0 \alpha_1^2 \alpha_4 \alpha_7 - 
       68 \alpha_0^2 \alpha_2 \alpha_4 \alpha_7 - 
       60 \alpha_0 \alpha_1 \alpha_2 \alpha_4 \alpha_7 - 
       20 \alpha_0 \alpha_2^2 \alpha_4 \alpha_7 - 
       42 \alpha_0^2 \alpha_3 \alpha_4 \alpha_7 - 
       36 \alpha_0 \alpha_1 \alpha_3 \alpha_4 \alpha_7 - 
       24 \alpha_0 \alpha_2 \alpha_3 \alpha_4 \alpha_7 - 
       6 \alpha_0 \alpha_3^2 \alpha_4 \alpha_7 - 
       8 \alpha_0^2 \alpha_4^2 \alpha_7 - 
       6 \alpha_0 \alpha_1 \alpha_4^2 \alpha_7 - 
       4 \alpha_0 \alpha_2 \alpha_4^2 \alpha_7 - 
       2 \alpha_0 \alpha_3 \alpha_4^2 \alpha_7 - 
       298 \alpha_0^3 \alpha_6 \alpha_7 - 
       400 \alpha_0^2 \alpha_1 \alpha_6 \alpha_7 - 
       156 \alpha_0 \alpha_1^2 \alpha_6 \alpha_7 - 
       300 \alpha_0^2 \alpha_2 \alpha_6 \alpha_7 - 
       234 \alpha_0 \alpha_1 \alpha_2 \alpha_6 \alpha_7 - 
       84 \alpha_0 \alpha_2^2 \alpha_6 \alpha_7 - 
       200 \alpha_0^2 \alpha_3 \alpha_6 \alpha_7 - 
       156 \alpha_0 \alpha_1 \alpha_3 \alpha_6 \alpha_7 - 
       112 \alpha_0 \alpha_2 \alpha_3 \alpha_6 \alpha_7 - 
       34 \alpha_0 \alpha_3^2 \alpha_6 \alpha_7 - 
       100 \alpha_0^2 \alpha_4 \alpha_6 \alpha_7 - 
       78 \alpha_0 \alpha_1 \alpha_4 \alpha_6 \alpha_7 - 
       56 \alpha_0 \alpha_2 \alpha_4 \alpha_6 \alpha_7 - 
       34 \alpha_0 \alpha_3 \alpha_4 \alpha_6 \alpha_7 - 
       6 \alpha_0 \alpha_4^2 \alpha_6 \alpha_7 - 
       198 \alpha_0^2 \alpha_6^2 \alpha_7 - 
       144 \alpha_0 \alpha_1 \alpha_6^2 \alpha_7 - 
       108 \alpha_0 \alpha_2 \alpha_6^2 \alpha_7 - 
       72 \alpha_0 \alpha_3 \alpha_6^2 \alpha_7 - 
       36 \alpha_0 \alpha_4 \alpha_6^2 \alpha_7 - 
       44 \alpha_0 \alpha_6^3 \alpha_7 - 
       78 \alpha_0^3 \alpha_7^2 - 
       112 \alpha_0^2 \alpha_1 \alpha_7^2 - 
       44 \alpha_0 \alpha_1^2 \alpha_7^2 - 
       84 \alpha_0^2 \alpha_2 \alpha_7^2 - 
       66 \alpha_0 \alpha_1 \alpha_2 \alpha_7^2 - 
       24 \alpha_0 \alpha_2^2 \alpha_7^2 - 
       56 \alpha_0^2 \alpha_3 \alpha_7^2 - 
       44 \alpha_0 \alpha_1 \alpha_3 \alpha_7^2 - 
       32 \alpha_0 \alpha_2 \alpha_3 \alpha_7^2 - 
       10 \alpha_0 \alpha_3^2 \alpha_7^2 - 
       28 \alpha_0^2 \alpha_4 \alpha_7^2 - 
       22 \alpha_0 \alpha_1 \alpha_4 \alpha_7^2 - 
       16 \alpha_0 \alpha_2 \alpha_4 \alpha_7^2 - 
       10 \alpha_0 \alpha_3 \alpha_4 \alpha_7^2 - 
       2 \alpha_0 \alpha_4^2 \alpha_7^2 - 
       106 \alpha_0^2 \alpha_6 \alpha_7^2 - 
       80 \alpha_0 \alpha_1 \alpha_6 \alpha_7^2 - 
       60 \alpha_0 \alpha_2 \alpha_6 \alpha_7^2 - 
       40 \alpha_0 \alpha_3 \alpha_6 \alpha_7^2 - 
       20 \alpha_0 \alpha_4 \alpha_6 \alpha_7^2 - 
       36 \alpha_0 \alpha_6^2 \alpha_7^2 - 
       20 \alpha_0^2 \alpha_7^3 - 
       16 \alpha_0 \alpha_1 \alpha_7^3 - 
       12 \alpha_0 \alpha_2 \alpha_7^3 - 
       8 \alpha_0 \alpha_3 \alpha_7^3 - 
       4 \alpha_0 \alpha_4 \alpha_7^3 - 
       14 \alpha_0 \alpha_6 \alpha_7^3 - 2 \alpha_0 \alpha_7^4 - 
       285 \alpha_0^4 \alpha_8 - 
       576 \alpha_0^3 \alpha_1 \alpha_8 - 
       342 \alpha_0^2 \alpha_1^2 \alpha_8 - 
       108 \alpha_0 \alpha_1^3 \alpha_8 - 
       15 \alpha_1^4 \alpha_8 - 
       432 \alpha_0^3 \alpha_2 \alpha_8 - 
       513 \alpha_0^2 \alpha_1 \alpha_2 \alpha_8 - 
       243 \alpha_0 \alpha_1^2 \alpha_2 \alpha_8 - 
       45 \alpha_1^3 \alpha_2 \alpha_8 - 
       189 \alpha_0^2 \alpha_2^2 \alpha_8 - 
       171 \alpha_0 \alpha_1 \alpha_2^2 \alpha_8 - 
       48 \alpha_1^2 \alpha_2^2 \alpha_8 - 
       36 \alpha_0 \alpha_2^3 \alpha_8 - 
       21 \alpha_1 \alpha_2^3 \alpha_8 - 3 \alpha_2^4 \alpha_8 - 
       288 \alpha_0^3 \alpha_3 \alpha_8 - 
       342 \alpha_0^2 \alpha_1 \alpha_3 \alpha_8 - 
       162 \alpha_0 \alpha_1^2 \alpha_3 \alpha_8 - 
       30 \alpha_1^3 \alpha_3 \alpha_8 - 
       252 \alpha_0^2 \alpha_2 \alpha_3 \alpha_8 - 
       228 \alpha_0 \alpha_1 \alpha_2 \alpha_3 \alpha_8 - 
       64 \alpha_1^2 \alpha_2 \alpha_3 \alpha_8 - 
       72 \alpha_0 \alpha_2^2 \alpha_3 \alpha_8 - 
       42 \alpha_1 \alpha_2^2 \alpha_3 \alpha_8 - 
       8 \alpha_2^3 \alpha_3 \alpha_8 - 
       81 \alpha_0^2 \alpha_3^2 \alpha_8 - 
       66 \alpha_0 \alpha_1 \alpha_3^2 \alpha_8 - 
       19 \alpha_1^2 \alpha_3^2 \alpha_8 - 
       42 \alpha_0 \alpha_2 \alpha_3^2 \alpha_8 - 
       25 \alpha_1 \alpha_2 \alpha_3^2 \alpha_8 - 
       7 \alpha_2^2 \alpha_3^2 \alpha_8 - 
       6 \alpha_0 \alpha_3^3 \alpha_8 - 
       4 \alpha_1 \alpha_3^3 \alpha_8 - 
       2 \alpha_2 \alpha_3^3 \alpha_8 - 
       144 \alpha_0^3 \alpha_4 \alpha_8 - 
       171 \alpha_0^2 \alpha_1 \alpha_4 \alpha_8 - 
       81 \alpha_0 \alpha_1^2 \alpha_4 \alpha_8 - 
       15 \alpha_1^3 \alpha_4 \alpha_8 - 
       126 \alpha_0^2 \alpha_2 \alpha_4 \alpha_8 - 
       114 \alpha_0 \alpha_1 \alpha_2 \alpha_4 \alpha_8 - 
       32 \alpha_1^2 \alpha_2 \alpha_4 \alpha_8 - 
       36 \alpha_0 \alpha_2^2 \alpha_4 \alpha_8 - 
       21 \alpha_1 \alpha_2^2 \alpha_4 \alpha_8 - 
       4 \alpha_2^3 \alpha_4 \alpha_8 - 
       81 \alpha_0^2 \alpha_3 \alpha_4 \alpha_8 - 
       66 \alpha_0 \alpha_1 \alpha_3 \alpha_4 \alpha_8 - 
       19 \alpha_1^2 \alpha_3 \alpha_4 \alpha_8 - 
       42 \alpha_0 \alpha_2 \alpha_3 \alpha_4 \alpha_8 - 
       25 \alpha_1 \alpha_2 \alpha_3 \alpha_4 \alpha_8 - 
       7 \alpha_2^2 \alpha_3 \alpha_4 \alpha_8 - 
       9 \alpha_0 \alpha_3^2 \alpha_4 \alpha_8 - 
       6 \alpha_1 \alpha_3^2 \alpha_4 \alpha_8 - 
       3 \alpha_2 \alpha_3^2 \alpha_4 \alpha_8 - 
       18 \alpha_0^2 \alpha_4^2 \alpha_8 - 
       9 \alpha_0 \alpha_1 \alpha_4^2 \alpha_8 - 
       3 \alpha_1^2 \alpha_4^2 \alpha_8 - 
       6 \alpha_0 \alpha_2 \alpha_4^2 \alpha_8 - 
       4 \alpha_1 \alpha_2 \alpha_4^2 \alpha_8 - \alpha_2^2 \
\alpha_4^2 \alpha_8 - 3 \alpha_0 \alpha_3 \alpha_4^2 \alpha_8 - 
       2 \alpha_1 \alpha_3 \alpha_4^2 \alpha_8 - \alpha_2 \
\alpha_3 \alpha_4^2 \alpha_8 - 
       600 \alpha_0^3 \alpha_6 \alpha_8 - 
       864 \alpha_0^2 \alpha_1 \alpha_6 \alpha_8 - 
       376 \alpha_0 \alpha_1^2 \alpha_6 \alpha_8 - 
       56 \alpha_1^3 \alpha_6 \alpha_8 - 
       648 \alpha_0^2 \alpha_2 \alpha_6 \alpha_8 - 
       564 \alpha_0 \alpha_1 \alpha_2 \alpha_6 \alpha_8 - 
       126 \alpha_1^2 \alpha_2 \alpha_6 \alpha_8 - 
       204 \alpha_0 \alpha_2^2 \alpha_6 \alpha_8 - 
       90 \alpha_1 \alpha_2^2 \alpha_6 \alpha_8 - 
       20 \alpha_2^3 \alpha_6 \alpha_8 - 
       432 \alpha_0^2 \alpha_3 \alpha_6 \alpha_8 - 
       376 \alpha_0 \alpha_1 \alpha_3 \alpha_6 \alpha_8 - 
       84 \alpha_1^2 \alpha_3 \alpha_6 \alpha_8 - 
       272 \alpha_0 \alpha_2 \alpha_3 \alpha_6 \alpha_8 - 
       120 \alpha_1 \alpha_2 \alpha_3 \alpha_6 \alpha_8 - 
       40 \alpha_2^2 \alpha_3 \alpha_6 \alpha_8 - 
       84 \alpha_0 \alpha_3^2 \alpha_6 \alpha_8 - 
       36 \alpha_1 \alpha_3^2 \alpha_6 \alpha_8 - 
       24 \alpha_2 \alpha_3^2 \alpha_6 \alpha_8 - 
       4 \alpha_3^3 \alpha_6 \alpha_8 - 
       216 \alpha_0^2 \alpha_4 \alpha_6 \alpha_8 - 
       188 \alpha_0 \alpha_1 \alpha_4 \alpha_6 \alpha_8 - 
       42 \alpha_1^2 \alpha_4 \alpha_6 \alpha_8 - 
       136 \alpha_0 \alpha_2 \alpha_4 \alpha_6 \alpha_8 - 
       60 \alpha_1 \alpha_2 \alpha_4 \alpha_6 \alpha_8 - 
       20 \alpha_2^2 \alpha_4 \alpha_6 \alpha_8 - 
       84 \alpha_0 \alpha_3 \alpha_4 \alpha_6 \alpha_8 - 
       36 \alpha_1 \alpha_3 \alpha_4 \alpha_6 \alpha_8 - 
       24 \alpha_2 \alpha_3 \alpha_4 \alpha_6 \alpha_8 - 
       6 \alpha_3^2 \alpha_4 \alpha_6 \alpha_8 - 
       16 \alpha_0 \alpha_4^2 \alpha_6 \alpha_8 - 
       6 \alpha_1 \alpha_4^2 \alpha_6 \alpha_8 - 
       4 \alpha_2 \alpha_4^2 \alpha_6 \alpha_8 - 
       2 \alpha_3 \alpha_4^2 \alpha_6 \alpha_8 - 
       447 \alpha_0^2 \alpha_6^2 \alpha_8 - 
       400 \alpha_0 \alpha_1 \alpha_6^2 \alpha_8 - 
       78 \alpha_1^2 \alpha_6^2 \alpha_8 - 
       300 \alpha_0 \alpha_2 \alpha_6^2 \alpha_8 - 
       117 \alpha_1 \alpha_2 \alpha_6^2 \alpha_8 - 
       42 \alpha_2^2 \alpha_6^2 \alpha_8 - 
       200 \alpha_0 \alpha_3 \alpha_6^2 \alpha_8 - 
       78 \alpha_1 \alpha_3 \alpha_6^2 \alpha_8 - 
       56 \alpha_2 \alpha_3 \alpha_6^2 \alpha_8 - 
       17 \alpha_3^2 \alpha_6^2 \alpha_8 - 
       100 \alpha_0 \alpha_4 \alpha_6^2 \alpha_8 - 
       39 \alpha_1 \alpha_4 \alpha_6^2 \alpha_8 - 
       28 \alpha_2 \alpha_4 \alpha_6^2 \alpha_8 - 
       17 \alpha_3 \alpha_4 \alpha_6^2 \alpha_8 - 
       3 \alpha_4^2 \alpha_6^2 \alpha_8 - 
       132 \alpha_0 \alpha_6^3 \alpha_8 - 
       48 \alpha_1 \alpha_6^3 \alpha_8 - 
       36 \alpha_2 \alpha_6^3 \alpha_8 - 
       24 \alpha_3 \alpha_6^3 \alpha_8 - 
       12 \alpha_4 \alpha_6^3 \alpha_8 - 
       11 \alpha_6^4 \alpha_8 - 
       300 \alpha_0^3 \alpha_7 \alpha_8 - 
       432 \alpha_0^2 \alpha_1 \alpha_7 \alpha_8 - 
       188 \alpha_0 \alpha_1^2 \alpha_7 \alpha_8 - 
       28 \alpha_1^3 \alpha_7 \alpha_8 - 
       324 \alpha_0^2 \alpha_2 \alpha_7 \alpha_8 - 
       282 \alpha_0 \alpha_1 \alpha_2 \alpha_7 \alpha_8 - 
       63 \alpha_1^2 \alpha_2 \alpha_7 \alpha_8 - 
       102 \alpha_0 \alpha_2^2 \alpha_7 \alpha_8 - 
       45 \alpha_1 \alpha_2^2 \alpha_7 \alpha_8 - 
       10 \alpha_2^3 \alpha_7 \alpha_8 - 
       216 \alpha_0^2 \alpha_3 \alpha_7 \alpha_8 - 
       188 \alpha_0 \alpha_1 \alpha_3 \alpha_7 \alpha_8 - 
       42 \alpha_1^2 \alpha_3 \alpha_7 \alpha_8 - 
       136 \alpha_0 \alpha_2 \alpha_3 \alpha_7 \alpha_8 - 
       60 \alpha_1 \alpha_2 \alpha_3 \alpha_7 \alpha_8 - 
       20 \alpha_2^2 \alpha_3 \alpha_7 \alpha_8 - 
       42 \alpha_0 \alpha_3^2 \alpha_7 \alpha_8 - 
       18 \alpha_1 \alpha_3^2 \alpha_7 \alpha_8 - 
       12 \alpha_2 \alpha_3^2 \alpha_7 \alpha_8 - 
       2 \alpha_3^3 \alpha_7 \alpha_8 - 
       108 \alpha_0^2 \alpha_4 \alpha_7 \alpha_8 - 
       94 \alpha_0 \alpha_1 \alpha_4 \alpha_7 \alpha_8 - 
       21 \alpha_1^2 \alpha_4 \alpha_7 \alpha_8 - 
       68 \alpha_0 \alpha_2 \alpha_4 \alpha_7 \alpha_8 - 
       30 \alpha_1 \alpha_2 \alpha_4 \alpha_7 \alpha_8 - 
       10 \alpha_2^2 \alpha_4 \alpha_7 \alpha_8 - 
       42 \alpha_0 \alpha_3 \alpha_4 \alpha_7 \alpha_8 - 
       18 \alpha_1 \alpha_3 \alpha_4 \alpha_7 \alpha_8 - 
       12 \alpha_2 \alpha_3 \alpha_4 \alpha_7 \alpha_8 - 
       3 \alpha_3^2 \alpha_4 \alpha_7 \alpha_8 - 
       8 \alpha_0 \alpha_4^2 \alpha_7 \alpha_8 - 
       3 \alpha_1 \alpha_4^2 \alpha_7 \alpha_8 - 
       2 \alpha_2 \alpha_4^2 \alpha_7 \alpha_8 - \alpha_3 \
\alpha_4^2 \alpha_7 \alpha_8 - 
       447 \alpha_0^2 \alpha_6 \alpha_7 \alpha_8 - 
       400 \alpha_0 \alpha_1 \alpha_6 \alpha_7 \alpha_8 - 
       78 \alpha_1^2 \alpha_6 \alpha_7 \alpha_8 - 
       300 \alpha_0 \alpha_2 \alpha_6 \alpha_7 \alpha_8 - 
       117 \alpha_1 \alpha_2 \alpha_6 \alpha_7 \alpha_8 - 
       42 \alpha_2^2 \alpha_6 \alpha_7 \alpha_8 - 
       200 \alpha_0 \alpha_3 \alpha_6 \alpha_7 \alpha_8 - 
       78 \alpha_1 \alpha_3 \alpha_6 \alpha_7 \alpha_8 - 
       56 \alpha_2 \alpha_3 \alpha_6 \alpha_7 \alpha_8 - 
       17 \alpha_3^2 \alpha_6 \alpha_7 \alpha_8 - 
       100 \alpha_0 \alpha_4 \alpha_6 \alpha_7 \alpha_8 - 
       39 \alpha_1 \alpha_4 \alpha_6 \alpha_7 \alpha_8 - 
       28 \alpha_2 \alpha_4 \alpha_6 \alpha_7 \alpha_8 - 
       17 \alpha_3 \alpha_4 \alpha_6 \alpha_7 \alpha_8 - 
       3 \alpha_4^2 \alpha_6 \alpha_7 \alpha_8 - 
       198 \alpha_0 \alpha_6^2 \alpha_7 \alpha_8 - 
       72 \alpha_1 \alpha_6^2 \alpha_7 \alpha_8 - 
       54 \alpha_2 \alpha_6^2 \alpha_7 \alpha_8 - 
       36 \alpha_3 \alpha_6^2 \alpha_7 \alpha_8 - 
       18 \alpha_4 \alpha_6^2 \alpha_7 \alpha_8 - 
       22 \alpha_6^3 \alpha_7 \alpha_8 - 
       117 \alpha_0^2 \alpha_7^2 \alpha_8 - 
       112 \alpha_0 \alpha_1 \alpha_7^2 \alpha_8 - 
       22 \alpha_1^2 \alpha_7^2 \alpha_8 - 
       84 \alpha_0 \alpha_2 \alpha_7^2 \alpha_8 - 
       33 \alpha_1 \alpha_2 \alpha_7^2 \alpha_8 - 
       12 \alpha_2^2 \alpha_7^2 \alpha_8 - 
       56 \alpha_0 \alpha_3 \alpha_7^2 \alpha_8 - 
       22 \alpha_1 \alpha_3 \alpha_7^2 \alpha_8 - 
       16 \alpha_2 \alpha_3 \alpha_7^2 \alpha_8 - 
       5 \alpha_3^2 \alpha_7^2 \alpha_8 - 
       28 \alpha_0 \alpha_4 \alpha_7^2 \alpha_8 - 
       11 \alpha_1 \alpha_4 \alpha_7^2 \alpha_8 - 
       8 \alpha_2 \alpha_4 \alpha_7^2 \alpha_8 - 
       5 \alpha_3 \alpha_4 \alpha_7^2 \alpha_8 - \alpha_4^2 \
\alpha_7^2 \alpha_8 - 
       106 \alpha_0 \alpha_6 \alpha_7^2 \alpha_8 - 
       40 \alpha_1 \alpha_6 \alpha_7^2 \alpha_8 - 
       30 \alpha_2 \alpha_6 \alpha_7^2 \alpha_8 - 
       20 \alpha_3 \alpha_6 \alpha_7^2 \alpha_8 - 
       10 \alpha_4 \alpha_6 \alpha_7^2 \alpha_8 - 
       18 \alpha_6^2 \alpha_7^2 \alpha_8 - 
       20 \alpha_0 \alpha_7^3 \alpha_8 - 
       8 \alpha_1 \alpha_7^3 \alpha_8 - 
       6 \alpha_2 \alpha_7^3 \alpha_8 - 
       4 \alpha_3 \alpha_7^3 \alpha_8 - 
       2 \alpha_4 \alpha_7^3 \alpha_8 - 
       7 \alpha_6 \alpha_7^3 \alpha_8 - \alpha_7^4 \alpha_8 - 
       274 \alpha_0^3 \alpha_8^2 - 
       408 \alpha_0^2 \alpha_1 \alpha_8^2 - 
       170 \alpha_0 \alpha_1^2 \alpha_8^2 - 
       32 \alpha_1^3 \alpha_8^2 - 
       306 \alpha_0^2 \alpha_2 \alpha_8^2 - 
       255 \alpha_0 \alpha_1 \alpha_2 \alpha_8^2 - 
       72 \alpha_1^2 \alpha_2 \alpha_8^2 - 
       93 \alpha_0 \alpha_2^2 \alpha_8^2 - 
       51 \alpha_1 \alpha_2^2 \alpha_8^2 - 
       11 \alpha_2^3 \alpha_8^2 - 
       204 \alpha_0^2 \alpha_3 \alpha_8^2 - 
       170 \alpha_0 \alpha_1 \alpha_3 \alpha_8^2 - 
       48 \alpha_1^2 \alpha_3 \alpha_8^2 - 
       124 \alpha_0 \alpha_2 \alpha_3 \alpha_8^2 - 
       68 \alpha_1 \alpha_2 \alpha_3 \alpha_8^2 - 
       22 \alpha_2^2 \alpha_3 \alpha_8^2 - 
       39 \alpha_0 \alpha_3^2 \alpha_8^2 - 
       20 \alpha_1 \alpha_3^2 \alpha_8^2 - 
       13 \alpha_2 \alpha_3^2 \alpha_8^2 - 
       2 \alpha_3^3 \alpha_8^2 - 
       102 \alpha_0^2 \alpha_4 \alpha_8^2 - 
       85 \alpha_0 \alpha_1 \alpha_4 \alpha_8^2 - 
       24 \alpha_1^2 \alpha_4 \alpha_8^2 - 
       62 \alpha_0 \alpha_2 \alpha_4 \alpha_8^2 - 
       34 \alpha_1 \alpha_2 \alpha_4 \alpha_8^2 - 
       11 \alpha_2^2 \alpha_4 \alpha_8^2 - 
       39 \alpha_0 \alpha_3 \alpha_4 \alpha_8^2 - 
       20 \alpha_1 \alpha_3 \alpha_4 \alpha_8^2 - 
       13 \alpha_2 \alpha_3 \alpha_4 \alpha_8^2 - 
       3 \alpha_3^2 \alpha_4 \alpha_8^2 - 
       8 \alpha_0 \alpha_4^2 \alpha_8^2 - 
       3 \alpha_1 \alpha_4^2 \alpha_8^2 - 
       2 \alpha_2 \alpha_4^2 \alpha_8^2 - \alpha_3 \alpha_4^2 \
\alpha_8^2 - 432 \alpha_0^2 \alpha_6 \alpha_8^2 - 
       416 \alpha_0 \alpha_1 \alpha_6 \alpha_8^2 - 
       100 \alpha_1^2 \alpha_6 \alpha_8^2 - 
       312 \alpha_0 \alpha_2 \alpha_6 \alpha_8^2 - 
       150 \alpha_1 \alpha_2 \alpha_6 \alpha_8^2 - 
       54 \alpha_2^2 \alpha_6 \alpha_8^2 - 
       208 \alpha_0 \alpha_3 \alpha_6 \alpha_8^2 - 
       100 \alpha_1 \alpha_3 \alpha_6 \alpha_8^2 - 
       72 \alpha_2 \alpha_3 \alpha_6 \alpha_8^2 - 
       22 \alpha_3^2 \alpha_6 \alpha_8^2 - 
       104 \alpha_0 \alpha_4 \alpha_6 \alpha_8^2 - 
       50 \alpha_1 \alpha_4 \alpha_6 \alpha_8^2 - 
       36 \alpha_2 \alpha_4 \alpha_6 \alpha_8^2 - 
       22 \alpha_3 \alpha_4 \alpha_6 \alpha_8^2 - 
       4 \alpha_4^2 \alpha_6 \alpha_8^2 - 
       215 \alpha_0 \alpha_6^2 \alpha_8^2 - 
       100 \alpha_1 \alpha_6^2 \alpha_8^2 - 
       75 \alpha_2 \alpha_6^2 \alpha_8^2 - 
       50 \alpha_3 \alpha_6^2 \alpha_8^2 - 
       25 \alpha_4 \alpha_6^2 \alpha_8^2 - 
       32 \alpha_6^3 \alpha_8^2 - 
       216 \alpha_0^2 \alpha_7 \alpha_8^2 - 
       208 \alpha_0 \alpha_1 \alpha_7 \alpha_8^2 - 
       50 \alpha_1^2 \alpha_7 \alpha_8^2 - 
       156 \alpha_0 \alpha_2 \alpha_7 \alpha_8^2 - 
       75 \alpha_1 \alpha_2 \alpha_7 \alpha_8^2 - 
       27 \alpha_2^2 \alpha_7 \alpha_8^2 - 
       104 \alpha_0 \alpha_3 \alpha_7 \alpha_8^2 - 
       50 \alpha_1 \alpha_3 \alpha_7 \alpha_8^2 - 
       36 \alpha_2 \alpha_3 \alpha_7 \alpha_8^2 - 
       11 \alpha_3^2 \alpha_7 \alpha_8^2 - 
       52 \alpha_0 \alpha_4 \alpha_7 \alpha_8^2 - 
       25 \alpha_1 \alpha_4 \alpha_7 \alpha_8^2 - 
       18 \alpha_2 \alpha_4 \alpha_7 \alpha_8^2 - 
       11 \alpha_3 \alpha_4 \alpha_7 \alpha_8^2 - 
       2 \alpha_4^2 \alpha_7 \alpha_8^2 - 
       215 \alpha_0 \alpha_6 \alpha_7 \alpha_8^2 - 
       100 \alpha_1 \alpha_6 \alpha_7 \alpha_8^2 - 
       75 \alpha_2 \alpha_6 \alpha_7 \alpha_8^2 - 
       50 \alpha_3 \alpha_6 \alpha_7 \alpha_8^2 - 
       25 \alpha_4 \alpha_6 \alpha_7 \alpha_8^2 - 
       48 \alpha_6^2 \alpha_7 \alpha_8^2 - 
       57 \alpha_0 \alpha_7^2 \alpha_8^2 - 
       28 \alpha_1 \alpha_7^2 \alpha_8^2 - 
       21 \alpha_2 \alpha_7^2 \alpha_8^2 - 
       14 \alpha_3 \alpha_7^2 \alpha_8^2 - 
       7 \alpha_4 \alpha_7^2 \alpha_8^2 - 
       26 \alpha_6 \alpha_7^2 \alpha_8^2 - 
       5 \alpha_7^3 \alpha_8^2 - 126 \alpha_0^2 \alpha_8^3 - 
       120 \alpha_0 \alpha_1 \alpha_8^3 - 
       28 \alpha_1^2 \alpha_8^3 - 
       90 \alpha_0 \alpha_2 \alpha_8^3 - 
       42 \alpha_1 \alpha_2 \alpha_8^3 - 
       15 \alpha_2^2 \alpha_8^3 - 
       60 \alpha_0 \alpha_3 \alpha_8^3 - 
       28 \alpha_1 \alpha_3 \alpha_8^3 - 
       20 \alpha_2 \alpha_3 \alpha_8^3 - 
       6 \alpha_3^2 \alpha_8^3 - 
       30 \alpha_0 \alpha_4 \alpha_8^3 - 
       14 \alpha_1 \alpha_4 \alpha_8^3 - 
       10 \alpha_2 \alpha_4 \alpha_8^3 - 
       6 \alpha_3 \alpha_4 \alpha_8^3 - \alpha_4^2 \alpha_8^3 - 
       132 \alpha_0 \alpha_6 \alpha_8^3 - 
       64 \alpha_1 \alpha_6 \alpha_8^3 - 
       48 \alpha_2 \alpha_6 \alpha_8^3 - 
       32 \alpha_3 \alpha_6 \alpha_8^3 - 
       16 \alpha_4 \alpha_6 \alpha_8^3 - 
       33 \alpha_6^2 \alpha_8^3 - 
       66 \alpha_0 \alpha_7 \alpha_8^3 - 
       32 \alpha_1 \alpha_7 \alpha_8^3 - 
       24 \alpha_2 \alpha_7 \alpha_8^3 - 
       16 \alpha_3 \alpha_7 \alpha_8^3 - 
       8 \alpha_4 \alpha_7 \alpha_8^3 - 
       33 \alpha_6 \alpha_7 \alpha_8^3 - 
       9 \alpha_7^2 \alpha_8^3 - 27 \alpha_0 \alpha_8^4 - 
       12 \alpha_1 \alpha_8^4 - 9 \alpha_2 \alpha_8^4 - 
       6 \alpha_3 \alpha_8^4 - 3 \alpha_4 \alpha_8^4 - 
       14 \alpha_6 \alpha_8^4 - 7 \alpha_7 \alpha_8^4 - 
       2 \alpha_8^5) + 
    q (72 \alpha_0^5 + 192 \alpha_0^4 \alpha_1 + 
       188 \alpha_0^3 \alpha_1^2 + 108 \alpha_0^2 \alpha_1^3 + 
       24 \alpha_0 \alpha_1^4 - 4 \alpha_1^5 + 
       144 \alpha_0^4 \alpha_2 + 
       282 \alpha_0^3 \alpha_1 \alpha_2 + 
       243 \alpha_0^2 \alpha_1^2 \alpha_2 + 
       72 \alpha_0 \alpha_1^3 \alpha_2 - 
       15 \alpha_1^4 \alpha_2 + 102 \alpha_0^3 \alpha_2^2 + 
       171 \alpha_0^2 \alpha_1 \alpha_2^2 + 
       78 \alpha_0 \alpha_1^2 \alpha_2^2 - 
       21 \alpha_1^3 \alpha_2^2 + 36 \alpha_0^2 \alpha_2^3 + 
       36 \alpha_0 \alpha_1 \alpha_2^3 - 
       13 \alpha_1^2 \alpha_2^3 + 6 \alpha_0 \alpha_2^4 - 
       3 \alpha_1 \alpha_2^4 + 96 \alpha_0^4 \alpha_3 + 
       188 \alpha_0^3 \alpha_1 \alpha_3 + 
       162 \alpha_0^2 \alpha_1^2 \alpha_3 + 
       48 \alpha_0 \alpha_1^3 \alpha_3 - 
       10 \alpha_1^4 \alpha_3 + 
       136 \alpha_0^3 \alpha_2 \alpha_3 + 
       228 \alpha_0^2 \alpha_1 \alpha_2 \alpha_3 + 
       104 \alpha_0 \alpha_1^2 \alpha_2 \alpha_3 - 
       28 \alpha_1^3 \alpha_2 \alpha_3 + 
       72 \alpha_0^2 \alpha_2^2 \alpha_3 + 
       72 \alpha_0 \alpha_1 \alpha_2^2 \alpha_3 - 
       26 \alpha_1^2 \alpha_2^2 \alpha_3 + 
       16 \alpha_0 \alpha_2^3 \alpha_3 - 
       8 \alpha_1 \alpha_2^3 \alpha_3 + 
       42 \alpha_0^3 \alpha_3^2 + 
       66 \alpha_0^2 \alpha_1 \alpha_3^2 + 
       32 \alpha_0 \alpha_1^2 \alpha_3^2 - 
       8 \alpha_1^3 \alpha_3^2 + 
       42 \alpha_0^2 \alpha_2 \alpha_3^2 + 
       44 \alpha_0 \alpha_1 \alpha_2 \alpha_3^2 - 
       15 \alpha_1^2 \alpha_2 \alpha_3^2 + 
       14 \alpha_0 \alpha_2^2 \alpha_3^2 - 
       7 \alpha_1 \alpha_2^2 \alpha_3^2 + 
       6 \alpha_0^2 \alpha_3^3 + 
       8 \alpha_0 \alpha_1 \alpha_3^3 - 
       2 \alpha_1^2 \alpha_3^3 + 
       4 \alpha_0 \alpha_2 \alpha_3^3 - 
       2 \alpha_1 \alpha_2 \alpha_3^3 + 48 \alpha_0^4 \alpha_4 + 
       94 \alpha_0^3 \alpha_1 \alpha_4 + 
       81 \alpha_0^2 \alpha_1^2 \alpha_4 + 
       24 \alpha_0 \alpha_1^3 \alpha_4 - 5 \alpha_1^4 \alpha_4 + 
       68 \alpha_0^3 \alpha_2 \alpha_4 + 
       114 \alpha_0^2 \alpha_1 \alpha_2 \alpha_4 + 
       52 \alpha_0 \alpha_1^2 \alpha_2 \alpha_4 - 
       14 \alpha_1^3 \alpha_2 \alpha_4 + 
       36 \alpha_0^2 \alpha_2^2 \alpha_4 + 
       36 \alpha_0 \alpha_1 \alpha_2^2 \alpha_4 - 
       13 \alpha_1^2 \alpha_2^2 \alpha_4 + 
       8 \alpha_0 \alpha_2^3 \alpha_4 - 
       4 \alpha_1 \alpha_2^3 \alpha_4 + 
       42 \alpha_0^3 \alpha_3 \alpha_4 + 
       66 \alpha_0^2 \alpha_1 \alpha_3 \alpha_4 + 
       32 \alpha_0 \alpha_1^2 \alpha_3 \alpha_4 - 
       8 \alpha_1^3 \alpha_3 \alpha_4 + 
       42 \alpha_0^2 \alpha_2 \alpha_3 \alpha_4 + 
       44 \alpha_0 \alpha_1 \alpha_2 \alpha_3 \alpha_4 - 
       15 \alpha_1^2 \alpha_2 \alpha_3 \alpha_4 + 
       14 \alpha_0 \alpha_2^2 \alpha_3 \alpha_4 - 
       7 \alpha_1 \alpha_2^2 \alpha_3 \alpha_4 + 
       9 \alpha_0^2 \alpha_3^2 \alpha_4 + 
       12 \alpha_0 \alpha_1 \alpha_3^2 \alpha_4 - 
       3 \alpha_1^2 \alpha_3^2 \alpha_4 + 
       6 \alpha_0 \alpha_2 \alpha_3^2 \alpha_4 - 
       3 \alpha_1 \alpha_2 \alpha_3^2 \alpha_4 + 
       8 \alpha_0^3 \alpha_4^2 + 
       9 \alpha_0^2 \alpha_1 \alpha_4^2 + 
       6 \alpha_0 \alpha_1^2 \alpha_4^2 - \alpha_1^3 \alpha_4^2 \
+ 6 \alpha_0^2 \alpha_2 \alpha_4^2 + 
       8 \alpha_0 \alpha_1 \alpha_2 \alpha_4^2 - 
       2 \alpha_1^2 \alpha_2 \alpha_4^2 + 
       2 \alpha_0 \alpha_2^2 \alpha_4^2 - \alpha_1 \alpha_2^2 \
\alpha_4^2 + 3 \alpha_0^2 \alpha_3 \alpha_4^2 + 
       4 \alpha_0 \alpha_1 \alpha_3 \alpha_4^2 - \alpha_1^2 \
\alpha_3 \alpha_4^2 + 
       2 \alpha_0 \alpha_2 \alpha_3 \alpha_4^2 - \alpha_1 \
\alpha_2 \alpha_3 \alpha_4^2 + 240 \alpha_0^4 \alpha_6 + 
       512 \alpha_0^3 \alpha_1 \alpha_6 + 
       376 \alpha_0^2 \alpha_1^2 \alpha_6 + 
       112 \alpha_0 \alpha_1^3 \alpha_6 - 
       4 \alpha_1^4 \alpha_6 + 
       384 \alpha_0^3 \alpha_2 \alpha_6 + 
       564 \alpha_0^2 \alpha_1 \alpha_2 \alpha_6 + 
       252 \alpha_0 \alpha_1^2 \alpha_2 \alpha_6 - 
       12 \alpha_1^3 \alpha_2 \alpha_6 + 
       204 \alpha_0^2 \alpha_2^2 \alpha_6 + 
       180 \alpha_0 \alpha_1 \alpha_2^2 \alpha_6 - 
       12 \alpha_1^2 \alpha_2^2 \alpha_6 + 
       40 \alpha_0 \alpha_2^3 \alpha_6 - 
       4 \alpha_1 \alpha_2^3 \alpha_6 + 
       256 \alpha_0^3 \alpha_3 \alpha_6 + 
       376 \alpha_0^2 \alpha_1 \alpha_3 \alpha_6 + 
       168 \alpha_0 \alpha_1^2 \alpha_3 \alpha_6 - 
       8 \alpha_1^3 \alpha_3 \alpha_6 + 
       272 \alpha_0^2 \alpha_2 \alpha_3 \alpha_6 + 
       240 \alpha_0 \alpha_1 \alpha_2 \alpha_3 \alpha_6 - 
       16 \alpha_1^2 \alpha_2 \alpha_3 \alpha_6 + 
       80 \alpha_0 \alpha_2^2 \alpha_3 \alpha_6 - 
       8 \alpha_1 \alpha_2^2 \alpha_3 \alpha_6 + 
       84 \alpha_0^2 \alpha_3^2 \alpha_6 + 
       72 \alpha_0 \alpha_1 \alpha_3^2 \alpha_6 - 
       4 \alpha_1^2 \alpha_3^2 \alpha_6 + 
       48 \alpha_0 \alpha_2 \alpha_3^2 \alpha_6 - 
       4 \alpha_1 \alpha_2 \alpha_3^2 \alpha_6 + 
       8 \alpha_0 \alpha_3^3 \alpha_6 + 
       128 \alpha_0^3 \alpha_4 \alpha_6 + 
       188 \alpha_0^2 \alpha_1 \alpha_4 \alpha_6 + 
       84 \alpha_0 \alpha_1^2 \alpha_4 \alpha_6 - 
       4 \alpha_1^3 \alpha_4 \alpha_6 + 
       136 \alpha_0^2 \alpha_2 \alpha_4 \alpha_6 + 
       120 \alpha_0 \alpha_1 \alpha_2 \alpha_4 \alpha_6 - 
       8 \alpha_1^2 \alpha_2 \alpha_4 \alpha_6 + 
       40 \alpha_0 \alpha_2^2 \alpha_4 \alpha_6 - 
       4 \alpha_1 \alpha_2^2 \alpha_4 \alpha_6 + 
       84 \alpha_0^2 \alpha_3 \alpha_4 \alpha_6 + 
       72 \alpha_0 \alpha_1 \alpha_3 \alpha_4 \alpha_6 - 
       4 \alpha_1^2 \alpha_3 \alpha_4 \alpha_6 + 
       48 \alpha_0 \alpha_2 \alpha_3 \alpha_4 \alpha_6 - 
       4 \alpha_1 \alpha_2 \alpha_3 \alpha_4 \alpha_6 + 
       12 \alpha_0 \alpha_3^2 \alpha_4 \alpha_6 + 
       16 \alpha_0^2 \alpha_4^2 \alpha_6 + 
       12 \alpha_0 \alpha_1 \alpha_4^2 \alpha_6 + 
       8 \alpha_0 \alpha_2 \alpha_4^2 \alpha_6 + 
       4 \alpha_0 \alpha_3 \alpha_4^2 \alpha_6 + 
       314 \alpha_0^3 \alpha_6^2 + 
       496 \alpha_0^2 \alpha_1 \alpha_6^2 + 
       232 \alpha_0 \alpha_1^2 \alpha_6^2 + 
       20 \alpha_1^3 \alpha_6^2 + 
       372 \alpha_0^2 \alpha_2 \alpha_6^2 + 
       348 \alpha_0 \alpha_1 \alpha_2 \alpha_6^2 + 
       45 \alpha_1^2 \alpha_2 \alpha_6^2 + 
       126 \alpha_0 \alpha_2^2 \alpha_6^2 + 
       33 \alpha_1 \alpha_2^2 \alpha_6^2 + 
       8 \alpha_2^3 \alpha_6^2 + 
       248 \alpha_0^2 \alpha_3 \alpha_6^2 + 
       232 \alpha_0 \alpha_1 \alpha_3 \alpha_6^2 + 
       30 \alpha_1^2 \alpha_3 \alpha_6^2 + 
       168 \alpha_0 \alpha_2 \alpha_3 \alpha_6^2 + 
       44 \alpha_1 \alpha_2 \alpha_3 \alpha_6^2 + 
       16 \alpha_2^2 \alpha_3 \alpha_6^2 + 
       52 \alpha_0 \alpha_3^2 \alpha_6^2 + 
       14 \alpha_1 \alpha_3^2 \alpha_6^2 + 
       10 \alpha_2 \alpha_3^2 \alpha_6^2 + 
       2 \alpha_3^3 \alpha_6^2 + 
       124 \alpha_0^2 \alpha_4 \alpha_6^2 + 
       116 \alpha_0 \alpha_1 \alpha_4 \alpha_6^2 + 
       15 \alpha_1^2 \alpha_4 \alpha_6^2 + 
       84 \alpha_0 \alpha_2 \alpha_4 \alpha_6^2 + 
       22 \alpha_1 \alpha_2 \alpha_4 \alpha_6^2 + 
       8 \alpha_2^2 \alpha_4 \alpha_6^2 + 
       52 \alpha_0 \alpha_3 \alpha_4 \alpha_6^2 + 
       14 \alpha_1 \alpha_3 \alpha_4 \alpha_6^2 + 
       10 \alpha_2 \alpha_3 \alpha_4 \alpha_6^2 + 
       3 \alpha_3^2 \alpha_4 \alpha_6^2 + 
       10 \alpha_0 \alpha_4^2 \alpha_6^2 + 
       3 \alpha_1 \alpha_4^2 \alpha_6^2 + 
       2 \alpha_2 \alpha_4^2 \alpha_6^2 + \alpha_3 \alpha_4^2 \
\alpha_6^2 + 202 \alpha_0^2 \alpha_6^3 + 
       208 \alpha_0 \alpha_1 \alpha_6^3 + 
       44 \alpha_1^2 \alpha_6^3 + 
       156 \alpha_0 \alpha_2 \alpha_6^3 + 
       66 \alpha_1 \alpha_2 \alpha_6^3 + 
       24 \alpha_2^2 \alpha_6^3 + 
       104 \alpha_0 \alpha_3 \alpha_6^3 + 
       44 \alpha_1 \alpha_3 \alpha_6^3 + 
       32 \alpha_2 \alpha_3 \alpha_6^3 + 
       10 \alpha_3^2 \alpha_6^3 + 
       52 \alpha_0 \alpha_4 \alpha_6^3 + 
       22 \alpha_1 \alpha_4 \alpha_6^3 + 
       16 \alpha_2 \alpha_4 \alpha_6^3 + 
       10 \alpha_3 \alpha_4 \alpha_6^3 + 
       2 \alpha_4^2 \alpha_6^3 + 64 \alpha_0 \alpha_6^4 + 
       32 \alpha_1 \alpha_6^4 + 24 \alpha_2 \alpha_6^4 + 
       16 \alpha_3 \alpha_6^4 + 8 \alpha_4 \alpha_6^4 + 
       8 \alpha_6^5 + 120 \alpha_0^4 \alpha_7 + 
       256 \alpha_0^3 \alpha_1 \alpha_7 + 
       188 \alpha_0^2 \alpha_1^2 \alpha_7 + 
       56 \alpha_0 \alpha_1^3 \alpha_7 - 2 \alpha_1^4 \alpha_7 + 
       192 \alpha_0^3 \alpha_2 \alpha_7 + 
       282 \alpha_0^2 \alpha_1 \alpha_2 \alpha_7 + 
       126 \alpha_0 \alpha_1^2 \alpha_2 \alpha_7 - 
       6 \alpha_1^3 \alpha_2 \alpha_7 + 
       102 \alpha_0^2 \alpha_2^2 \alpha_7 + 
       90 \alpha_0 \alpha_1 \alpha_2^2 \alpha_7 - 
       6 \alpha_1^2 \alpha_2^2 \alpha_7 + 
       20 \alpha_0 \alpha_2^3 \alpha_7 - 
       2 \alpha_1 \alpha_2^3 \alpha_7 + 
       128 \alpha_0^3 \alpha_3 \alpha_7 + 
       188 \alpha_0^2 \alpha_1 \alpha_3 \alpha_7 + 
       84 \alpha_0 \alpha_1^2 \alpha_3 \alpha_7 - 
       4 \alpha_1^3 \alpha_3 \alpha_7 + 
       136 \alpha_0^2 \alpha_2 \alpha_3 \alpha_7 + 
       120 \alpha_0 \alpha_1 \alpha_2 \alpha_3 \alpha_7 - 
       8 \alpha_1^2 \alpha_2 \alpha_3 \alpha_7 + 
       40 \alpha_0 \alpha_2^2 \alpha_3 \alpha_7 - 
       4 \alpha_1 \alpha_2^2 \alpha_3 \alpha_7 + 
       42 \alpha_0^2 \alpha_3^2 \alpha_7 + 
       36 \alpha_0 \alpha_1 \alpha_3^2 \alpha_7 - 
       2 \alpha_1^2 \alpha_3^2 \alpha_7 + 
       24 \alpha_0 \alpha_2 \alpha_3^2 \alpha_7 - 
       2 \alpha_1 \alpha_2 \alpha_3^2 \alpha_7 + 
       4 \alpha_0 \alpha_3^3 \alpha_7 + 
       64 \alpha_0^3 \alpha_4 \alpha_7 + 
       94 \alpha_0^2 \alpha_1 \alpha_4 \alpha_7 + 
       42 \alpha_0 \alpha_1^2 \alpha_4 \alpha_7 - 
       2 \alpha_1^3 \alpha_4 \alpha_7 + 
       68 \alpha_0^2 \alpha_2 \alpha_4 \alpha_7 + 
       60 \alpha_0 \alpha_1 \alpha_2 \alpha_4 \alpha_7 - 
       4 \alpha_1^2 \alpha_2 \alpha_4 \alpha_7 + 
       20 \alpha_0 \alpha_2^2 \alpha_4 \alpha_7 - 
       2 \alpha_1 \alpha_2^2 \alpha_4 \alpha_7 + 
       42 \alpha_0^2 \alpha_3 \alpha_4 \alpha_7 + 
       36 \alpha_0 \alpha_1 \alpha_3 \alpha_4 \alpha_7 - 
       2 \alpha_1^2 \alpha_3 \alpha_4 \alpha_7 + 
       24 \alpha_0 \alpha_2 \alpha_3 \alpha_4 \alpha_7 - 
       2 \alpha_1 \alpha_2 \alpha_3 \alpha_4 \alpha_7 + 
       6 \alpha_0 \alpha_3^2 \alpha_4 \alpha_7 + 
       8 \alpha_0^2 \alpha_4^2 \alpha_7 + 
       6 \alpha_0 \alpha_1 \alpha_4^2 \alpha_7 + 
       4 \alpha_0 \alpha_2 \alpha_4^2 \alpha_7 + 
       2 \alpha_0 \alpha_3 \alpha_4^2 \alpha_7 + 
       314 \alpha_0^3 \alpha_6 \alpha_7 + 
       496 \alpha_0^2 \alpha_1 \alpha_6 \alpha_7 + 
       232 \alpha_0 \alpha_1^2 \alpha_6 \alpha_7 + 
       20 \alpha_1^3 \alpha_6 \alpha_7 + 
       372 \alpha_0^2 \alpha_2 \alpha_6 \alpha_7 + 
       348 \alpha_0 \alpha_1 \alpha_2 \alpha_6 \alpha_7 + 
       45 \alpha_1^2 \alpha_2 \alpha_6 \alpha_7 + 
       126 \alpha_0 \alpha_2^2 \alpha_6 \alpha_7 + 
       33 \alpha_1 \alpha_2^2 \alpha_6 \alpha_7 + 
       8 \alpha_2^3 \alpha_6 \alpha_7 + 
       248 \alpha_0^2 \alpha_3 \alpha_6 \alpha_7 + 
       232 \alpha_0 \alpha_1 \alpha_3 \alpha_6 \alpha_7 + 
       30 \alpha_1^2 \alpha_3 \alpha_6 \alpha_7 + 
       168 \alpha_0 \alpha_2 \alpha_3 \alpha_6 \alpha_7 + 
       44 \alpha_1 \alpha_2 \alpha_3 \alpha_6 \alpha_7 + 
       16 \alpha_2^2 \alpha_3 \alpha_6 \alpha_7 + 
       52 \alpha_0 \alpha_3^2 \alpha_6 \alpha_7 + 
       14 \alpha_1 \alpha_3^2 \alpha_6 \alpha_7 + 
       10 \alpha_2 \alpha_3^2 \alpha_6 \alpha_7 + 
       2 \alpha_3^3 \alpha_6 \alpha_7 + 
       124 \alpha_0^2 \alpha_4 \alpha_6 \alpha_7 + 
       116 \alpha_0 \alpha_1 \alpha_4 \alpha_6 \alpha_7 + 
       15 \alpha_1^2 \alpha_4 \alpha_6 \alpha_7 + 
       84 \alpha_0 \alpha_2 \alpha_4 \alpha_6 \alpha_7 + 
       22 \alpha_1 \alpha_2 \alpha_4 \alpha_6 \alpha_7 + 
       8 \alpha_2^2 \alpha_4 \alpha_6 \alpha_7 + 
       52 \alpha_0 \alpha_3 \alpha_4 \alpha_6 \alpha_7 + 
       14 \alpha_1 \alpha_3 \alpha_4 \alpha_6 \alpha_7 + 
       10 \alpha_2 \alpha_3 \alpha_4 \alpha_6 \alpha_7 + 
       3 \alpha_3^2 \alpha_4 \alpha_6 \alpha_7 + 
       10 \alpha_0 \alpha_4^2 \alpha_6 \alpha_7 + 
       3 \alpha_1 \alpha_4^2 \alpha_6 \alpha_7 + 
       2 \alpha_2 \alpha_4^2 \alpha_6 \alpha_7 + \alpha_3 \
\alpha_4^2 \alpha_6 \alpha_7 + 
       303 \alpha_0^2 \alpha_6^2 \alpha_7 + 
       312 \alpha_0 \alpha_1 \alpha_6^2 \alpha_7 + 
       66 \alpha_1^2 \alpha_6^2 \alpha_7 + 
       234 \alpha_0 \alpha_2 \alpha_6^2 \alpha_7 + 
       99 \alpha_1 \alpha_2 \alpha_6^2 \alpha_7 + 
       36 \alpha_2^2 \alpha_6^2 \alpha_7 + 
       156 \alpha_0 \alpha_3 \alpha_6^2 \alpha_7 + 
       66 \alpha_1 \alpha_3 \alpha_6^2 \alpha_7 + 
       48 \alpha_2 \alpha_3 \alpha_6^2 \alpha_7 + 
       15 \alpha_3^2 \alpha_6^2 \alpha_7 + 
       78 \alpha_0 \alpha_4 \alpha_6^2 \alpha_7 + 
       33 \alpha_1 \alpha_4 \alpha_6^2 \alpha_7 + 
       24 \alpha_2 \alpha_4 \alpha_6^2 \alpha_7 + 
       15 \alpha_3 \alpha_4 \alpha_6^2 \alpha_7 + 
       3 \alpha_4^2 \alpha_6^2 \alpha_7 + 
       128 \alpha_0 \alpha_6^3 \alpha_7 + 
       64 \alpha_1 \alpha_6^3 \alpha_7 + 
       48 \alpha_2 \alpha_6^3 \alpha_7 + 
       32 \alpha_3 \alpha_6^3 \alpha_7 + 
       16 \alpha_4 \alpha_6^3 \alpha_7 + 
       20 \alpha_6^4 \alpha_7 + 74 \alpha_0^3 \alpha_7^2 + 
       112 \alpha_0^2 \alpha_1 \alpha_7^2 + 
       44 \alpha_0 \alpha_1^2 \alpha_7^2 + 
       84 \alpha_0^2 \alpha_2 \alpha_7^2 + 
       66 \alpha_0 \alpha_1 \alpha_2 \alpha_7^2 + 
       24 \alpha_0 \alpha_2^2 \alpha_7^2 + 
       56 \alpha_0^2 \alpha_3 \alpha_7^2 + 
       44 \alpha_0 \alpha_1 \alpha_3 \alpha_7^2 + 
       32 \alpha_0 \alpha_2 \alpha_3 \alpha_7^2 + 
       10 \alpha_0 \alpha_3^2 \alpha_7^2 + 
       28 \alpha_0^2 \alpha_4 \alpha_7^2 + 
       22 \alpha_0 \alpha_1 \alpha_4 \alpha_7^2 + 
       16 \alpha_0 \alpha_2 \alpha_4 \alpha_7^2 + 
       10 \alpha_0 \alpha_3 \alpha_4 \alpha_7^2 + 
       2 \alpha_0 \alpha_4^2 \alpha_7^2 + 
       141 \alpha_0^2 \alpha_6 \alpha_7^2 + 
       136 \alpha_0 \alpha_1 \alpha_6 \alpha_7^2 + 
       22 \alpha_1^2 \alpha_6 \alpha_7^2 + 
       102 \alpha_0 \alpha_2 \alpha_6 \alpha_7^2 + 
       33 \alpha_1 \alpha_2 \alpha_6 \alpha_7^2 + 
       12 \alpha_2^2 \alpha_6 \alpha_7^2 + 
       68 \alpha_0 \alpha_3 \alpha_6 \alpha_7^2 + 
       22 \alpha_1 \alpha_3 \alpha_6 \alpha_7^2 + 
       16 \alpha_2 \alpha_3 \alpha_6 \alpha_7^2 + 
       5 \alpha_3^2 \alpha_6 \alpha_7^2 + 
       34 \alpha_0 \alpha_4 \alpha_6 \alpha_7^2 + 
       11 \alpha_1 \alpha_4 \alpha_6 \alpha_7^2 + 
       8 \alpha_2 \alpha_4 \alpha_6 \alpha_7^2 + 
       5 \alpha_3 \alpha_4 \alpha_6 \alpha_7^2 + \alpha_4^2 \
\alpha_6 \alpha_7^2 + 88 \alpha_0 \alpha_6^2 \alpha_7^2 + 
       40 \alpha_1 \alpha_6^2 \alpha_7^2 + 
       30 \alpha_2 \alpha_6^2 \alpha_7^2 + 
       20 \alpha_3 \alpha_6^2 \alpha_7^2 + 
       10 \alpha_4 \alpha_6^2 \alpha_7^2 + 
       18 \alpha_6^3 \alpha_7^2 + 20 \alpha_0^2 \alpha_7^3 + 
       16 \alpha_0 \alpha_1 \alpha_7^3 + 
       12 \alpha_0 \alpha_2 \alpha_7^3 + 
       8 \alpha_0 \alpha_3 \alpha_7^3 + 
       4 \alpha_0 \alpha_4 \alpha_7^3 + 
       24 \alpha_0 \alpha_6 \alpha_7^3 + 
       8 \alpha_1 \alpha_6 \alpha_7^3 + 
       6 \alpha_2 \alpha_6 \alpha_7^3 + 
       4 \alpha_3 \alpha_6 \alpha_7^3 + 
       2 \alpha_4 \alpha_6 \alpha_7^3 + 
       7 \alpha_6^2 \alpha_7^3 + 
       2 \alpha_0 \alpha_7^4 + \alpha_6 \alpha_7^4 + 
       180 \alpha_0^4 \alpha_8 + 
       384 \alpha_0^3 \alpha_1 \alpha_8 + 
       282 \alpha_0^2 \alpha_1^2 \alpha_8 + 
       108 \alpha_0 \alpha_1^3 \alpha_8 + 
       12 \alpha_1^4 \alpha_8 + 
       288 \alpha_0^3 \alpha_2 \alpha_8 + 
       423 \alpha_0^2 \alpha_1 \alpha_2 \alpha_8 + 
       243 \alpha_0 \alpha_1^2 \alpha_2 \alpha_8 + 
       36 \alpha_1^3 \alpha_2 \alpha_8 + 
       153 \alpha_0^2 \alpha_2^2 \alpha_8 + 
       171 \alpha_0 \alpha_1 \alpha_2^2 \alpha_8 + 
       39 \alpha_1^2 \alpha_2^2 \alpha_8 + 
       36 \alpha_0 \alpha_2^3 \alpha_8 + 
       18 \alpha_1 \alpha_2^3 \alpha_8 + 3 \alpha_2^4 \alpha_8 + 
       192 \alpha_0^3 \alpha_3 \alpha_8 + 
       282 \alpha_0^2 \alpha_1 \alpha_3 \alpha_8 + 
       162 \alpha_0 \alpha_1^2 \alpha_3 \alpha_8 + 
       24 \alpha_1^3 \alpha_3 \alpha_8 + 
       204 \alpha_0^2 \alpha_2 \alpha_3 \alpha_8 + 
       228 \alpha_0 \alpha_1 \alpha_2 \alpha_3 \alpha_8 + 
       52 \alpha_1^2 \alpha_2 \alpha_3 \alpha_8 + 
       72 \alpha_0 \alpha_2^2 \alpha_3 \alpha_8 + 
       36 \alpha_1 \alpha_2^2 \alpha_3 \alpha_8 + 
       8 \alpha_2^3 \alpha_3 \alpha_8 + 
       63 \alpha_0^2 \alpha_3^2 \alpha_8 + 
       66 \alpha_0 \alpha_1 \alpha_3^2 \alpha_8 + 
       16 \alpha_1^2 \alpha_3^2 \alpha_8 + 
       42 \alpha_0 \alpha_2 \alpha_3^2 \alpha_8 + 
       22 \alpha_1 \alpha_2 \alpha_3^2 \alpha_8 + 
       7 \alpha_2^2 \alpha_3^2 \alpha_8 + 
       6 \alpha_0 \alpha_3^3 \alpha_8 + 
       4 \alpha_1 \alpha_3^3 \alpha_8 + 
       2 \alpha_2 \alpha_3^3 \alpha_8 + 
       96 \alpha_0^3 \alpha_4 \alpha_8 + 
       141 \alpha_0^2 \alpha_1 \alpha_4 \alpha_8 + 
       81 \alpha_0 \alpha_1^2 \alpha_4 \alpha_8 + 
       12 \alpha_1^3 \alpha_4 \alpha_8 + 
       102 \alpha_0^2 \alpha_2 \alpha_4 \alpha_8 + 
       114 \alpha_0 \alpha_1 \alpha_2 \alpha_4 \alpha_8 + 
       26 \alpha_1^2 \alpha_2 \alpha_4 \alpha_8 + 
       36 \alpha_0 \alpha_2^2 \alpha_4 \alpha_8 + 
       18 \alpha_1 \alpha_2^2 \alpha_4 \alpha_8 + 
       4 \alpha_2^3 \alpha_4 \alpha_8 + 
       63 \alpha_0^2 \alpha_3 \alpha_4 \alpha_8 + 
       66 \alpha_0 \alpha_1 \alpha_3 \alpha_4 \alpha_8 + 
       16 \alpha_1^2 \alpha_3 \alpha_4 \alpha_8 + 
       42 \alpha_0 \alpha_2 \alpha_3 \alpha_4 \alpha_8 + 
       22 \alpha_1 \alpha_2 \alpha_3 \alpha_4 \alpha_8 + 
       7 \alpha_2^2 \alpha_3 \alpha_4 \alpha_8 + 
       9 \alpha_0 \alpha_3^2 \alpha_4 \alpha_8 + 
       6 \alpha_1 \alpha_3^2 \alpha_4 \alpha_8 + 
       3 \alpha_2 \alpha_3^2 \alpha_4 \alpha_8 + 
       12 \alpha_0^2 \alpha_4^2 \alpha_8 + 
       9 \alpha_0 \alpha_1 \alpha_4^2 \alpha_8 + 
       3 \alpha_1^2 \alpha_4^2 \alpha_8 + 
       6 \alpha_0 \alpha_2 \alpha_4^2 \alpha_8 + 
       4 \alpha_1 \alpha_2 \alpha_4^2 \alpha_8 + \alpha_2^2 \
\alpha_4^2 \alpha_8 + 3 \alpha_0 \alpha_3 \alpha_4^2 \alpha_8 + 
       2 \alpha_1 \alpha_3 \alpha_4^2 \alpha_8 + \alpha_2 \
\alpha_3 \alpha_4^2 \alpha_8 + 
       480 \alpha_0^3 \alpha_6 \alpha_8 + 
       768 \alpha_0^2 \alpha_1 \alpha_6 \alpha_8 + 
       376 \alpha_0 \alpha_1^2 \alpha_6 \alpha_8 + 
       56 \alpha_1^3 \alpha_6 \alpha_8 + 
       576 \alpha_0^2 \alpha_2 \alpha_6 \alpha_8 + 
       564 \alpha_0 \alpha_1 \alpha_2 \alpha_6 \alpha_8 + 
       126 \alpha_1^2 \alpha_2 \alpha_6 \alpha_8 + 
       204 \alpha_0 \alpha_2^2 \alpha_6 \alpha_8 + 
       90 \alpha_1 \alpha_2^2 \alpha_6 \alpha_8 + 
       20 \alpha_2^3 \alpha_6 \alpha_8 + 
       384 \alpha_0^2 \alpha_3 \alpha_6 \alpha_8 + 
       376 \alpha_0 \alpha_1 \alpha_3 \alpha_6 \alpha_8 + 
       84 \alpha_1^2 \alpha_3 \alpha_6 \alpha_8 + 
       272 \alpha_0 \alpha_2 \alpha_3 \alpha_6 \alpha_8 + 
       120 \alpha_1 \alpha_2 \alpha_3 \alpha_6 \alpha_8 + 
       40 \alpha_2^2 \alpha_3 \alpha_6 \alpha_8 + 
       84 \alpha_0 \alpha_3^2 \alpha_6 \alpha_8 + 
       36 \alpha_1 \alpha_3^2 \alpha_6 \alpha_8 + 
       24 \alpha_2 \alpha_3^2 \alpha_6 \alpha_8 + 
       4 \alpha_3^3 \alpha_6 \alpha_8 + 
       192 \alpha_0^2 \alpha_4 \alpha_6 \alpha_8 + 
       188 \alpha_0 \alpha_1 \alpha_4 \alpha_6 \alpha_8 + 
       42 \alpha_1^2 \alpha_4 \alpha_6 \alpha_8 + 
       136 \alpha_0 \alpha_2 \alpha_4 \alpha_6 \alpha_8 + 
       60 \alpha_1 \alpha_2 \alpha_4 \alpha_6 \alpha_8 + 
       20 \alpha_2^2 \alpha_4 \alpha_6 \alpha_8 + 
       84 \alpha_0 \alpha_3 \alpha_4 \alpha_6 \alpha_8 + 
       36 \alpha_1 \alpha_3 \alpha_4 \alpha_6 \alpha_8 + 
       24 \alpha_2 \alpha_3 \alpha_4 \alpha_6 \alpha_8 + 
       6 \alpha_3^2 \alpha_4 \alpha_6 \alpha_8 + 
       16 \alpha_0 \alpha_4^2 \alpha_6 \alpha_8 + 
       6 \alpha_1 \alpha_4^2 \alpha_6 \alpha_8 + 
       4 \alpha_2 \alpha_4^2 \alpha_6 \alpha_8 + 
       2 \alpha_3 \alpha_4^2 \alpha_6 \alpha_8 + 
       471 \alpha_0^2 \alpha_6^2 \alpha_8 + 
       496 \alpha_0 \alpha_1 \alpha_6^2 \alpha_8 + 
       116 \alpha_1^2 \alpha_6^2 \alpha_8 + 
       372 \alpha_0 \alpha_2 \alpha_6^2 \alpha_8 + 
       174 \alpha_1 \alpha_2 \alpha_6^2 \alpha_8 + 
       63 \alpha_2^2 \alpha_6^2 \alpha_8 + 
       248 \alpha_0 \alpha_3 \alpha_6^2 \alpha_8 + 
       116 \alpha_1 \alpha_3 \alpha_6^2 \alpha_8 + 
       84 \alpha_2 \alpha_3 \alpha_6^2 \alpha_8 + 
       26 \alpha_3^2 \alpha_6^2 \alpha_8 + 
       124 \alpha_0 \alpha_4 \alpha_6^2 \alpha_8 + 
       58 \alpha_1 \alpha_4 \alpha_6^2 \alpha_8 + 
       42 \alpha_2 \alpha_4 \alpha_6^2 \alpha_8 + 
       26 \alpha_3 \alpha_4 \alpha_6^2 \alpha_8 + 
       5 \alpha_4^2 \alpha_6^2 \alpha_8 + 
       202 \alpha_0 \alpha_6^3 \alpha_8 + 
       104 \alpha_1 \alpha_6^3 \alpha_8 + 
       78 \alpha_2 \alpha_6^3 \alpha_8 + 
       52 \alpha_3 \alpha_6^3 \alpha_8 + 
       26 \alpha_4 \alpha_6^3 \alpha_8 + 
       32 \alpha_6^4 \alpha_8 + 
       240 \alpha_0^3 \alpha_7 \alpha_8 + 
       384 \alpha_0^2 \alpha_1 \alpha_7 \alpha_8 + 
       188 \alpha_0 \alpha_1^2 \alpha_7 \alpha_8 + 
       28 \alpha_1^3 \alpha_7 \alpha_8 + 
       288 \alpha_0^2 \alpha_2 \alpha_7 \alpha_8 + 
       282 \alpha_0 \alpha_1 \alpha_2 \alpha_7 \alpha_8 + 
       63 \alpha_1^2 \alpha_2 \alpha_7 \alpha_8 + 
       102 \alpha_0 \alpha_2^2 \alpha_7 \alpha_8 + 
       45 \alpha_1 \alpha_2^2 \alpha_7 \alpha_8 + 
       10 \alpha_2^3 \alpha_7 \alpha_8 + 
       192 \alpha_0^2 \alpha_3 \alpha_7 \alpha_8 + 
       188 \alpha_0 \alpha_1 \alpha_3 \alpha_7 \alpha_8 + 
       42 \alpha_1^2 \alpha_3 \alpha_7 \alpha_8 + 
       136 \alpha_0 \alpha_2 \alpha_3 \alpha_7 \alpha_8 + 
       60 \alpha_1 \alpha_2 \alpha_3 \alpha_7 \alpha_8 + 
       20 \alpha_2^2 \alpha_3 \alpha_7 \alpha_8 + 
       42 \alpha_0 \alpha_3^2 \alpha_7 \alpha_8 + 
       18 \alpha_1 \alpha_3^2 \alpha_7 \alpha_8 + 
       12 \alpha_2 \alpha_3^2 \alpha_7 \alpha_8 + 
       2 \alpha_3^3 \alpha_7 \alpha_8 + 
       96 \alpha_0^2 \alpha_4 \alpha_7 \alpha_8 + 
       94 \alpha_0 \alpha_1 \alpha_4 \alpha_7 \alpha_8 + 
       21 \alpha_1^2 \alpha_4 \alpha_7 \alpha_8 + 
       68 \alpha_0 \alpha_2 \alpha_4 \alpha_7 \alpha_8 + 
       30 \alpha_1 \alpha_2 \alpha_4 \alpha_7 \alpha_8 + 
       10 \alpha_2^2 \alpha_4 \alpha_7 \alpha_8 + 
       42 \alpha_0 \alpha_3 \alpha_4 \alpha_7 \alpha_8 + 
       18 \alpha_1 \alpha_3 \alpha_4 \alpha_7 \alpha_8 + 
       12 \alpha_2 \alpha_3 \alpha_4 \alpha_7 \alpha_8 + 
       3 \alpha_3^2 \alpha_4 \alpha_7 \alpha_8 + 
       8 \alpha_0 \alpha_4^2 \alpha_7 \alpha_8 + 
       3 \alpha_1 \alpha_4^2 \alpha_7 \alpha_8 + 
       2 \alpha_2 \alpha_4^2 \alpha_7 \alpha_8 + \alpha_3 \
\alpha_4^2 \alpha_7 \alpha_8 + 
       471 \alpha_0^2 \alpha_6 \alpha_7 \alpha_8 + 
       496 \alpha_0 \alpha_1 \alpha_6 \alpha_7 \alpha_8 + 
       116 \alpha_1^2 \alpha_6 \alpha_7 \alpha_8 + 
       372 \alpha_0 \alpha_2 \alpha_6 \alpha_7 \alpha_8 + 
       174 \alpha_1 \alpha_2 \alpha_6 \alpha_7 \alpha_8 + 
       63 \alpha_2^2 \alpha_6 \alpha_7 \alpha_8 + 
       248 \alpha_0 \alpha_3 \alpha_6 \alpha_7 \alpha_8 + 
       116 \alpha_1 \alpha_3 \alpha_6 \alpha_7 \alpha_8 + 
       84 \alpha_2 \alpha_3 \alpha_6 \alpha_7 \alpha_8 + 
       26 \alpha_3^2 \alpha_6 \alpha_7 \alpha_8 + 
       124 \alpha_0 \alpha_4 \alpha_6 \alpha_7 \alpha_8 + 
       58 \alpha_1 \alpha_4 \alpha_6 \alpha_7 \alpha_8 + 
       42 \alpha_2 \alpha_4 \alpha_6 \alpha_7 \alpha_8 + 
       26 \alpha_3 \alpha_4 \alpha_6 \alpha_7 \alpha_8 + 
       5 \alpha_4^2 \alpha_6 \alpha_7 \alpha_8 + 
       303 \alpha_0 \alpha_6^2 \alpha_7 \alpha_8 + 
       156 \alpha_1 \alpha_6^2 \alpha_7 \alpha_8 + 
       117 \alpha_2 \alpha_6^2 \alpha_7 \alpha_8 + 
       78 \alpha_3 \alpha_6^2 \alpha_7 \alpha_8 + 
       39 \alpha_4 \alpha_6^2 \alpha_7 \alpha_8 + 
       64 \alpha_6^3 \alpha_7 \alpha_8 + 
       111 \alpha_0^2 \alpha_7^2 \alpha_8 + 
       112 \alpha_0 \alpha_1 \alpha_7^2 \alpha_8 + 
       22 \alpha_1^2 \alpha_7^2 \alpha_8 + 
       84 \alpha_0 \alpha_2 \alpha_7^2 \alpha_8 + 
       33 \alpha_1 \alpha_2 \alpha_7^2 \alpha_8 + 
       12 \alpha_2^2 \alpha_7^2 \alpha_8 + 
       56 \alpha_0 \alpha_3 \alpha_7^2 \alpha_8 + 
       22 \alpha_1 \alpha_3 \alpha_7^2 \alpha_8 + 
       16 \alpha_2 \alpha_3 \alpha_7^2 \alpha_8 + 
       5 \alpha_3^2 \alpha_7^2 \alpha_8 + 
       28 \alpha_0 \alpha_4 \alpha_7^2 \alpha_8 + 
       11 \alpha_1 \alpha_4 \alpha_7^2 \alpha_8 + 
       8 \alpha_2 \alpha_4 \alpha_7^2 \alpha_8 + 
       5 \alpha_3 \alpha_4 \alpha_7^2 \alpha_8 + \alpha_4^2 \
\alpha_7^2 \alpha_8 + 
       141 \alpha_0 \alpha_6 \alpha_7^2 \alpha_8 + 
       68 \alpha_1 \alpha_6 \alpha_7^2 \alpha_8 + 
       51 \alpha_2 \alpha_6 \alpha_7^2 \alpha_8 + 
       34 \alpha_3 \alpha_6 \alpha_7^2 \alpha_8 + 
       17 \alpha_4 \alpha_6 \alpha_7^2 \alpha_8 + 
       44 \alpha_6^2 \alpha_7^2 \alpha_8 + 
       20 \alpha_0 \alpha_7^3 \alpha_8 + 
       8 \alpha_1 \alpha_7^3 \alpha_8 + 
       6 \alpha_2 \alpha_7^3 \alpha_8 + 
       4 \alpha_3 \alpha_7^3 \alpha_8 + 
       2 \alpha_4 \alpha_7^3 \alpha_8 + 
       12 \alpha_6 \alpha_7^3 \alpha_8 + \alpha_7^4 \alpha_8 + 
       178 \alpha_0^3 \alpha_8^2 + 
       288 \alpha_0^2 \alpha_1 \alpha_8^2 + 
       150 \alpha_0 \alpha_1^2 \alpha_8^2 + 
       32 \alpha_1^3 \alpha_8^2 + 
       216 \alpha_0^2 \alpha_2 \alpha_8^2 + 
       225 \alpha_0 \alpha_1 \alpha_2 \alpha_8^2 + 
       72 \alpha_1^2 \alpha_2 \alpha_8^2 + 
       81 \alpha_0 \alpha_2^2 \alpha_8^2 + 
       51 \alpha_1 \alpha_2^2 \alpha_8^2 + 
       11 \alpha_2^3 \alpha_8^2 + 
       144 \alpha_0^2 \alpha_3 \alpha_8^2 + 
       150 \alpha_0 \alpha_1 \alpha_3 \alpha_8^2 + 
       48 \alpha_1^2 \alpha_3 \alpha_8^2 + 
       108 \alpha_0 \alpha_2 \alpha_3 \alpha_8^2 + 
       68 \alpha_1 \alpha_2 \alpha_3 \alpha_8^2 + 
       22 \alpha_2^2 \alpha_3 \alpha_8^2 + 
       33 \alpha_0 \alpha_3^2 \alpha_8^2 + 
       20 \alpha_1 \alpha_3^2 \alpha_8^2 + 
       13 \alpha_2 \alpha_3^2 \alpha_8^2 + 
       2 \alpha_3^3 \alpha_8^2 + 
       72 \alpha_0^2 \alpha_4 \alpha_8^2 + 
       75 \alpha_0 \alpha_1 \alpha_4 \alpha_8^2 + 
       24 \alpha_1^2 \alpha_4 \alpha_8^2 + 
       54 \alpha_0 \alpha_2 \alpha_4 \alpha_8^2 + 
       34 \alpha_1 \alpha_2 \alpha_4 \alpha_8^2 + 
       11 \alpha_2^2 \alpha_4 \alpha_8^2 + 
       33 \alpha_0 \alpha_3 \alpha_4 \alpha_8^2 + 
       20 \alpha_1 \alpha_3 \alpha_4 \alpha_8^2 + 
       13 \alpha_2 \alpha_3 \alpha_4 \alpha_8^2 + 
       3 \alpha_3^2 \alpha_4 \alpha_8^2 + 
       6 \alpha_0 \alpha_4^2 \alpha_8^2 + 
       3 \alpha_1 \alpha_4^2 \alpha_8^2 + 
       2 \alpha_2 \alpha_4^2 \alpha_8^2 + \alpha_3 \alpha_4^2 \
\alpha_8^2 + 356 \alpha_0^2 \alpha_6 \alpha_8^2 + 
       384 \alpha_0 \alpha_1 \alpha_6 \alpha_8^2 + 
       100 \alpha_1^2 \alpha_6 \alpha_8^2 + 
       288 \alpha_0 \alpha_2 \alpha_6 \alpha_8^2 + 
       150 \alpha_1 \alpha_2 \alpha_6 \alpha_8^2 + 
       54 \alpha_2^2 \alpha_6 \alpha_8^2 + 
       192 \alpha_0 \alpha_3 \alpha_6 \alpha_8^2 + 
       100 \alpha_1 \alpha_3 \alpha_6 \alpha_8^2 + 
       72 \alpha_2 \alpha_3 \alpha_6 \alpha_8^2 + 
       22 \alpha_3^2 \alpha_6 \alpha_8^2 + 
       96 \alpha_0 \alpha_4 \alpha_6 \alpha_8^2 + 
       50 \alpha_1 \alpha_4 \alpha_6 \alpha_8^2 + 
       36 \alpha_2 \alpha_4 \alpha_6 \alpha_8^2 + 
       22 \alpha_3 \alpha_4 \alpha_6 \alpha_8^2 + 
       4 \alpha_4^2 \alpha_6 \alpha_8^2 + 
       233 \alpha_0 \alpha_6^2 \alpha_8^2 + 
       124 \alpha_1 \alpha_6^2 \alpha_8^2 + 
       93 \alpha_2 \alpha_6^2 \alpha_8^2 + 
       62 \alpha_3 \alpha_6^2 \alpha_8^2 + 
       31 \alpha_4 \alpha_6^2 \alpha_8^2 + 
       50 \alpha_6^3 \alpha_8^2 + 
       178 \alpha_0^2 \alpha_7 \alpha_8^2 + 
       192 \alpha_0 \alpha_1 \alpha_7 \alpha_8^2 + 
       50 \alpha_1^2 \alpha_7 \alpha_8^2 + 
       144 \alpha_0 \alpha_2 \alpha_7 \alpha_8^2 + 
       75 \alpha_1 \alpha_2 \alpha_7 \alpha_8^2 + 
       27 \alpha_2^2 \alpha_7 \alpha_8^2 + 
       96 \alpha_0 \alpha_3 \alpha_7 \alpha_8^2 + 
       50 \alpha_1 \alpha_3 \alpha_7 \alpha_8^2 + 
       36 \alpha_2 \alpha_3 \alpha_7 \alpha_8^2 + 
       11 \alpha_3^2 \alpha_7 \alpha_8^2 + 
       48 \alpha_0 \alpha_4 \alpha_7 \alpha_8^2 + 
       25 \alpha_1 \alpha_4 \alpha_7 \alpha_8^2 + 
       18 \alpha_2 \alpha_4 \alpha_7 \alpha_8^2 + 
       11 \alpha_3 \alpha_4 \alpha_7 \alpha_8^2 + 
       2 \alpha_4^2 \alpha_7 \alpha_8^2 + 
       233 \alpha_0 \alpha_6 \alpha_7 \alpha_8^2 + 
       124 \alpha_1 \alpha_6 \alpha_7 \alpha_8^2 + 
       93 \alpha_2 \alpha_6 \alpha_7 \alpha_8^2 + 
       62 \alpha_3 \alpha_6 \alpha_7 \alpha_8^2 + 
       31 \alpha_4 \alpha_6 \alpha_7 \alpha_8^2 + 
       75 \alpha_6^2 \alpha_7 \alpha_8^2 + 
       55 \alpha_0 \alpha_7^2 \alpha_8^2 + 
       28 \alpha_1 \alpha_7^2 \alpha_8^2 + 
       21 \alpha_2 \alpha_7^2 \alpha_8^2 + 
       14 \alpha_3 \alpha_7^2 \alpha_8^2 + 
       7 \alpha_4 \alpha_7^2 \alpha_8^2 + 
       35 \alpha_6 \alpha_7^2 \alpha_8^2 + 
       5 \alpha_7^3 \alpha_8^2 + 87 \alpha_0^2 \alpha_8^3 + 
       96 \alpha_0 \alpha_1 \alpha_8^3 + 
       28 \alpha_1^2 \alpha_8^3 + 
       72 \alpha_0 \alpha_2 \alpha_8^3 + 
       42 \alpha_1 \alpha_2 \alpha_8^3 + 
       15 \alpha_2^2 \alpha_8^3 + 
       48 \alpha_0 \alpha_3 \alpha_8^3 + 
       28 \alpha_1 \alpha_3 \alpha_8^3 + 
       20 \alpha_2 \alpha_3 \alpha_8^3 + 
       6 \alpha_3^2 \alpha_8^3 + 
       24 \alpha_0 \alpha_4 \alpha_8^3 + 
       14 \alpha_1 \alpha_4 \alpha_8^3 + 
       10 \alpha_2 \alpha_4 \alpha_8^3 + 
       6 \alpha_3 \alpha_4 \alpha_8^3 + \alpha_4^2 \alpha_8^3 + 
       116 \alpha_0 \alpha_6 \alpha_8^3 + 
       64 \alpha_1 \alpha_6 \alpha_8^3 + 
       48 \alpha_2 \alpha_6 \alpha_8^3 + 
       32 \alpha_3 \alpha_6 \alpha_8^3 + 
       16 \alpha_4 \alpha_6 \alpha_8^3 + 
       38 \alpha_6^2 \alpha_8^3 + 
       58 \alpha_0 \alpha_7 \alpha_8^3 + 
       32 \alpha_1 \alpha_7 \alpha_8^3 + 
       24 \alpha_2 \alpha_7 \alpha_8^3 + 
       16 \alpha_3 \alpha_7 \alpha_8^3 + 
       8 \alpha_4 \alpha_7 \alpha_8^3 + 
       38 \alpha_6 \alpha_7 \alpha_8^3 + 
       9 \alpha_7^2 \alpha_8^3 + 21 \alpha_0 \alpha_8^4 + 
       12 \alpha_1 \alpha_8^4 + 9 \alpha_2 \alpha_8^4 + 
       6 \alpha_3 \alpha_8^4 + 3 \alpha_4 \alpha_8^4 + 
       14 \alpha_6 \alpha_8^4 + 7 \alpha_7 \alpha_8^4 + 
       2 \alpha_8^5)).
$

Here, the constant parameters $\alpha_i$ satisfy the relation:
\begin{equation}
6\alpha_0+5\alpha_1+4\alpha_2+3\alpha_3+2\alpha_4+\alpha_5+4\alpha_6+2\alpha_7+3\alpha_8=0.
\end{equation}
\end{Theorem}

The holomorphy conditions $(C2),(C3)$ are new. Theorem 11.1 can be checked by a direct calculation.

\begin{Proposition}
The Hamiltonian $I$ is its first integral.
\end{Proposition}

\begin{Remark}
For the Hamiltonian system in each coordinate system $(x_i,y_i) \ (i=0,1,\ldots,8)$ given by $(C2)$ and $(C3)$ in Theorem 11.1, by eliminating $x_i$ or $y_i$, we obtain the second-order ordinary differential equation. However, its form is not normal \rm{(cf. \cite{Cosgrove1,Cosgrove2})}.
\end{Remark}

\section{Symmetry}

\begin{Theorem}
The system \eqref{E811} admits the affine Weyl group symmetry of type $E_8^{(1)}$ as the group of its B{\"a}cklund transformations whose generators $s_i, \ i=0,1,\ldots,8$ are explicitly given as follows{\rm : \rm}with the notation $(*):=(q,p,t;\alpha_0,\alpha_1,\ldots,\alpha_8)$,
\begin{align*}
        s_{0}: (*) &\rightarrow \left(q+\frac{\alpha_0}{p},p,t;-\alpha_0,\alpha_1+\alpha_0,\alpha_2,\alpha_3,\alpha_4,\alpha_5,\alpha_6+\alpha_0,\alpha_7,\alpha_8+\alpha_0 \right),\\
        s_{1}: (*) &\rightarrow \left(q,p-\frac{\alpha_1}{q},t;\alpha_0+\alpha_1,-\alpha_1,\alpha_2+\alpha_1,\alpha_3,\alpha_4,\alpha_5,\alpha_6,\alpha_7,\alpha_8 \right), \\
        s_{2}: (*) &\rightarrow (q,p,t;\alpha_0,\alpha_1+\alpha_2,-\alpha_2,\alpha_3+\alpha_2,\alpha_4,\alpha_5,\alpha_6,\alpha_7,\alpha_8), \\
        s_{3}: (*) &\rightarrow (q,p,t;\alpha_0,\alpha_1,\alpha_2+\alpha_3,-\alpha_3,\alpha_4+\alpha_3,\alpha_5,\alpha_6,\alpha_7,\alpha_8), \\
        s_{4}: (*) &\rightarrow \left(q,p,t;\alpha_0,\alpha_1,\alpha_2,\alpha_3+\alpha_4,-\alpha_4,\alpha_5+\alpha_4,\alpha_6,\alpha_7,\alpha_8 \right),\\
        s_{5}: (*) &\rightarrow (q,p,t;\alpha_0,\alpha_1,\alpha_2,\alpha_3,\alpha_4+\alpha_5,-\alpha_5,\alpha_6,\alpha_7,\alpha_8),\\
        s_{6}: (*) &\rightarrow \left(q,p-\frac{\alpha_6}{q-1},t;\alpha_0+\alpha_6,\alpha_1,\alpha_2,\alpha_3,\alpha_4,\alpha_5,-\alpha_6,\alpha_7+\alpha_6,\alpha_8 \right),\\
        s_{7}: (*) &\rightarrow (q,p,t;\alpha_0,\alpha_1,\alpha_2,\alpha_3,\alpha_4,\alpha_5,\alpha_6+\alpha_7,-\alpha_7,\alpha_8),\\
        s_{8}: (*) &\rightarrow (q,p,t;\alpha_0+\alpha_8,\alpha_1,\alpha_2,\alpha_3,\alpha_4,\alpha_5,\alpha_6,\alpha_7,-\alpha_8).
\end{align*}
\end{Theorem}
Theorem 12.1 can be checked by a direct calculation.

\section{Space of initial conditions}

\begin{Theorem}\label{3.1}
After a series of explicit blowing-ups at eleven points including the infinitely near points of ${\Sigma_2}$ and successive blowing-down along the $(-1)$-curves ${D^{(0)}}' \cong {\Bbb P}^1$, $D_{1}^{(1)} \cong {\Bbb P}^1$ and $D_{\infty}^{(1)} \cong {\Bbb P}^1$, we obtain the rational surface $\tilde{S}$ of the system \eqref{E811} and a birational morphism $\varphi:\tilde{S} \cdots \rightarrow {\Sigma_2}$. Its canonical divisor $K_{\tilde{S}}$ of $\tilde{S}$ is given by
\begin{align}
\begin{split}
K_{\tilde{S}}&=-D_{0}^{(1)}, \quad (D_{0}^{(1)})^2=-3, \ D_{0}^{(1)} \cong {\Bbb P}^1,
\end{split}
\end{align}
where the symbol ${D^{(0)}}'$ denotes the strict transform of $D^{(0)}$, $D_{\nu}^{(1)}$ denote the exceptional divisors and $-K_{\Sigma_2}=2D^{(0)}, \ D^{(0)} \cong {\Bbb P}^1, \ (D^{(0)})^2=2$.
\end{Theorem}

\begin{Theorem}
The space of initial conditions $S$ of the system \eqref{E811} is obtained by gluing ten copies of ${\Bbb C}^2$:
\begin{align}
\begin{split}
S&={\tilde{S}}-(-{K_{\tilde{S}}})_{red}\\
&={\Bbb C}^2 \cup \bigcup_{i=0}^{8} U_j,\\
&{\Bbb C}^2 \ni (q,p), \quad U_j \cong {\Bbb C}^2 \ni (x_j,y_j) \ (j=0,1,\ldots,8)
\end{split}
\end{align}
via the birational and symplectic transformations $r_j$ \rm{(see Theorem 11.1)}.
\end{Theorem}

{\bf Proof of Theorems 13.1 and 13.2.}

By a direct calculation, we see that the system \eqref{E811} has three accessible singular points $a_{\nu}^{(0)} \in D^{(0)} \quad (\nu=0,1,\infty)$:
\begin{align}
\begin{split}
&a_{\nu}^{(0)}=\{(z_2,w_2)=(\nu,0)\} \in U_2 \cap D^{(0)} \ (\nu=0,1),\\
&a_{\infty}^{(0)}=\{(z_3,w_3)=(0,0)\} \in U_3 \cap D^{(0)}.
\end{split}
\end{align}
We perform blowing-ups in ${\Sigma_2}$ at $a_{\nu}^{(0)}$, and let $D_{\nu}^{(1)}$ be the exceptional curves of the blowing-ups at $a_{\nu}^{(0)}$ for $\nu=0,1,\infty$. We can take three coordinate systems $(u_{\nu},v_{\nu})$ around the points at infinity of the exceptional curves $D_{\nu}^{(1)} \quad (\nu=0,1,\infty)$, where
\begin{align}
\begin{split}
&(u_{\nu},v_{\nu})=\left(\frac{z_2-\nu}{w_2},w_2 \right) \ (\nu=0,1),\\
&(u_{\infty},v_{\infty})=\left(\frac{z_3}{w_3},w_3 \right).
\end{split}
\end{align}

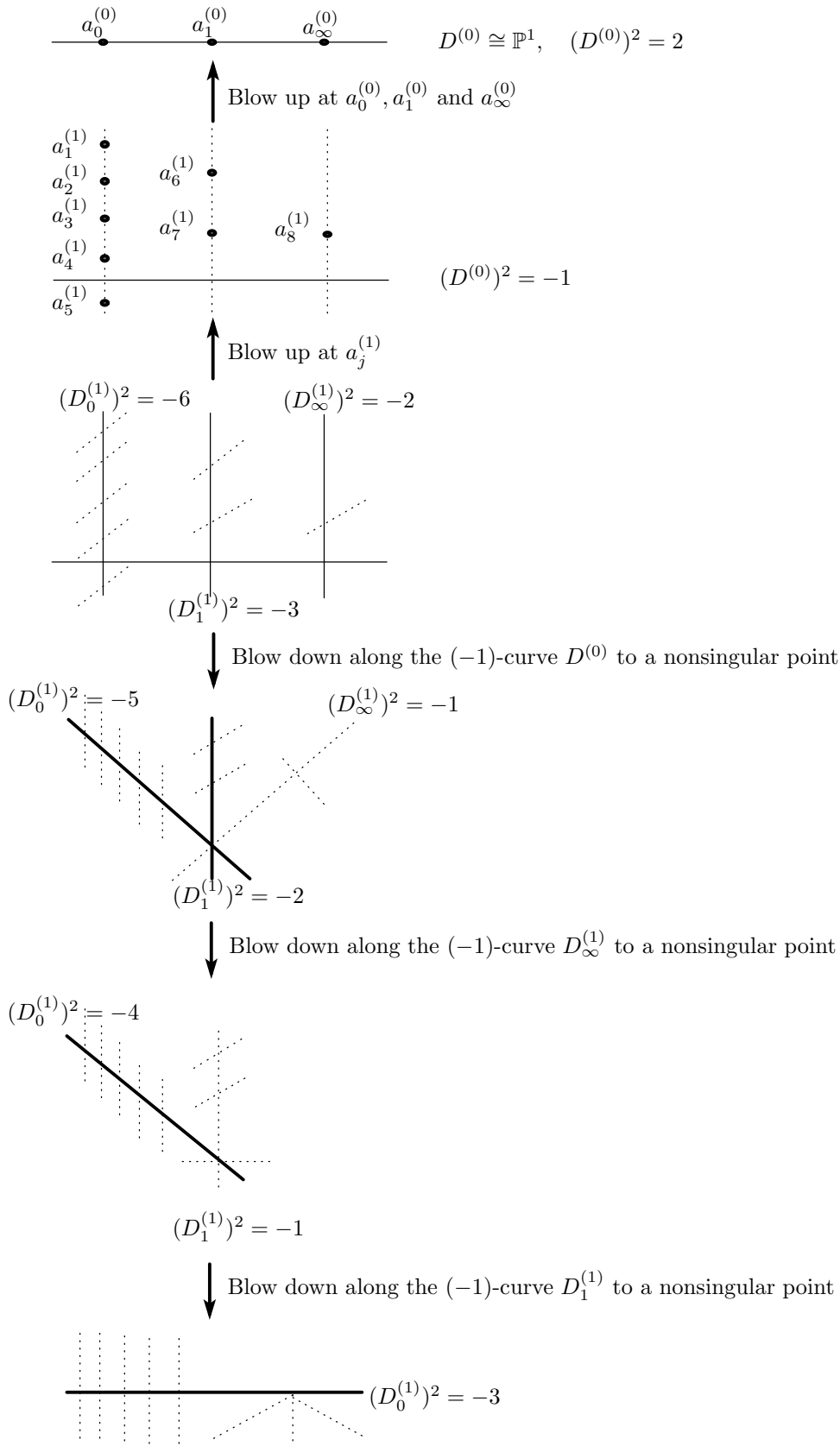
\begin{figure}
\unitlength 0.1in
\begin{picture}( 26.2000, 87.0700)( 18.4000,-87.1100)
\put(44.5000,-2.7800){\makebox(0,0)[lb]{$D^{(0)} \cong {\Bbb P}^1, \quad (D^{(0)})^2=2$}}%
%
\special{pn 8}%
\special{pa 2110 208}%
\special{pa 4140 208}%
\special{fp}%
%
\special{pn 20}%
\special{sh 0.600}%
\special{ar 2420 208 20 12  0.0000000 6.2831853}%
%
\special{pn 20}%
\special{sh 0.600}%
\special{ar 3080 208 20 12  0.0000000 6.2831853}%
%
\special{pn 20}%
\special{sh 0.600}%
\special{ar 3760 208 20 12  0.0000000 6.2831853}%
\put(22.9000,-1.7400){\makebox(0,0)[lb]{$a_0^{(0)}$}}%
\put(29.6000,-1.7400){\makebox(0,0)[lb]{$a_1^{(0)}$}}%
\put(36.3000,-1.7400){\makebox(0,0)[lb]{$a_{\infty}^{(0)}$}}%
%
\special{pn 20}%
\special{pa 3086 688}%
\special{pa 3086 350}%
\special{fp}%
\special{sh 1}%
\special{pa 3086 350}%
\special{pa 3066 416}%
\special{pa 3086 402}%
\special{pa 3106 416}%
\special{pa 3086 350}%
\special{fp}%
\put(31.7500,-6.2400){\makebox(0,0)[lb]{Blow up at $a_0^{(0)},a_1^{(0)}$ and $a_{\infty}^{(0)}$}}%
%
\special{pn 8}%
\special{pa 2120 1652}%
\special{pa 4150 1652}%
\special{fp}%
%
\special{pn 8}%
\special{pa 2430 1852}%
\special{pa 2430 738}%
\special{dt 0.045}%
%
\special{pn 20}%
\special{sh 0.600}%
\special{ar 2430 828 20 16  0.0000000 6.2831853}%
\put(21.1000,-9.1300){\makebox(0,0)[lb]{$a_1^{(1)}$}}%
%
\special{pn 8}%
\special{pa 3080 1844}%
\special{pa 3080 732}%
\special{dt 0.045}%
%
\special{pn 20}%
\special{sh 0.600}%
\special{ar 3080 1366 20 14  0.0000000 6.2831853}%
%
\special{pn 20}%
\special{sh 0.600}%
\special{ar 3080 1000 20 14  0.0000000 6.2831853}%
\put(27.6000,-10.7700){\makebox(0,0)[lb]{$a_6^{(1)}$}}%
\put(27.6000,-14.0800){\makebox(0,0)[lb]{$a_7^{(1)}$}}%
%
\special{pn 8}%
\special{pa 3780 1852}%
\special{pa 3780 738}%
\special{dt 0.045}%
%
\special{pn 20}%
\special{sh 0.600}%
\special{ar 3780 1374 20 12  0.0000000 6.2831853}%
\put(34.6000,-14.1700){\makebox(0,0)[lb]{$a_8^{(1)}$}}%
%
\special{pn 20}%
\special{pa 3086 2250}%
\special{pa 3086 1912}%
\special{fp}%
\special{sh 1}%
\special{pa 3086 1912}%
\special{pa 3066 1980}%
\special{pa 3086 1966}%
\special{pa 3106 1980}%
\special{pa 3086 1912}%
\special{fp}%
\put(31.7500,-21.8900){\makebox(0,0)[lb]{Blow up at $a_j^{(1)}$}}%
%
\special{pn 8}%
\special{pa 2110 3362}%
\special{pa 4140 3362}%
\special{fp}%
%
\special{pn 8}%
\special{pa 2420 3564}%
\special{pa 2420 2450}%
\special{fp}%
%
\special{pn 8}%
\special{pa 3070 3572}%
\special{pa 3070 2458}%
\special{fp}%
%
\special{pn 8}%
\special{pa 2970 3182}%
\special{pa 3320 2982}%
\special{dt 0.045}%
%
\special{pn 8}%
\special{pa 2970 2860}%
\special{pa 3270 2626}%
\special{dt 0.045}%
%
\special{pn 8}%
\special{pa 3760 3580}%
\special{pa 3760 2468}%
\special{fp}%
%
\special{pn 8}%
\special{pa 3660 3190}%
\special{pa 4010 2988}%
\special{dt 0.045}%
\put(44.6000,-17.1200){\makebox(0,0)[lb]{$(D^{(0)})^2=-1$}}%
\put(21.5000,-24.5000){\makebox(0,0)[lb]{$(D_{0}^{(1)})^2=-6$}}%
\put(28.1000,-37.2700){\makebox(0,0)[lb]{$(D_{1}^{(1)})^2=-3$}}%
\put(35.1000,-24.5800){\makebox(0,0)[lb]{$(D_{\infty}^{(1)})^2=-2$}}%
%
\special{pn 20}%
\special{pa 3090 3798}%
\special{pa 3090 4102}%
\special{fp}%
\special{sh 1}%
\special{pa 3090 4102}%
\special{pa 3110 4036}%
\special{pa 3090 4050}%
\special{pa 3070 4036}%
\special{pa 3090 4102}%
\special{fp}%
\put(32.0000,-40.1500){\makebox(0,0)[lb]{Blow down along the $(-1)$-curve $D^{(0)}$ to a nonsingular point}}%
%
\special{pn 20}%
\special{pa 2210 4318}%
\special{pa 3310 5284}%
\special{fp}%
\special{pa 3080 4310}%
\special{pa 3080 5284}%
\special{fp}%
%
\special{pn 8}%
\special{pa 2840 5294}%
\special{pa 3940 4338}%
\special{dt 0.045}%
%
\special{pn 8}%
\special{pa 2970 4528}%
\special{pa 3260 4354}%
\special{dt 0.045}%
%
\special{pn 8}%
\special{pa 2980 4764}%
\special{pa 3280 4590}%
\special{dt 0.045}%
%
\special{pn 8}%
\special{pa 3510 4556}%
\special{pa 3760 4832}%
\special{dt 0.045}%
%
\special{pn 20}%
\special{sh 0.600}%
\special{ar 2430 1052 20 16  0.0000000 6.2831853}%
%
\special{pn 20}%
\special{sh 0.600}%
\special{ar 2430 1278 20 14  0.0000000 6.2831853}%
%
\special{pn 20}%
\special{sh 0.600}%
\special{ar 2430 1520 20 16  0.0000000 6.2831853}%
%
\special{pn 20}%
\special{sh 0.600}%
\special{ar 2430 1790 20 14  0.0000000 6.2831853}%
\put(21.1000,-11.2700){\makebox(0,0)[lb]{$a_2^{(1)}$}}%
\put(21.1000,-13.3200){\makebox(0,0)[lb]{$a_3^{(1)}$}}%
\put(21.1000,-15.8500){\makebox(0,0)[lb]{$a_4^{(1)}$}}%
\put(21.1000,-18.6500){\makebox(0,0)[lb]{$a_5^{(1)}$}}%
%
\special{pn 8}%
\special{pa 2260 2696}%
\special{pa 2560 2458}%
\special{dt 0.045}%
%
\special{pn 8}%
\special{pa 2260 2872}%
\special{pa 2560 2636}%
\special{dt 0.045}%
%
\special{pn 8}%
\special{pa 2260 3124}%
\special{pa 2560 2888}%
\special{dt 0.045}%
%
\special{pn 8}%
\special{pa 2270 3338}%
\special{pa 2570 3104}%
\special{dt 0.045}%
%
\special{pn 8}%
\special{pa 2270 3628}%
\special{pa 2570 3392}%
\special{dt 0.045}%
%
\special{pn 8}%
\special{pa 2310 4170}%
\special{pa 2310 4608}%
\special{dt 0.045}%
%
\special{pn 8}%
\special{pa 2410 4272}%
\special{pa 2410 4712}%
\special{dt 0.045}%
%
\special{pn 8}%
\special{pa 2520 4374}%
\special{pa 2520 4814}%
\special{dt 0.045}%
%
\special{pn 8}%
\special{pa 2640 4516}%
\special{pa 2640 4954}%
\special{dt 0.045}%
%
\special{pn 8}%
\special{pa 2780 4598}%
\special{pa 2780 5038}%
\special{dt 0.045}%
\put(18.4000,-42.7200){\makebox(0,0)[lb]{$(D_{0}^{(1)})^2=-5$}}%
\put(28.4000,-54.6700){\makebox(0,0)[lb]{$(D_{1}^{(1)})^2=-2$}}%
\put(37.8000,-43.0100){\makebox(0,0)[lb]{$(D_{\infty}^{(1)})^2=-1$}}%
%
\special{pn 20}%
\special{pa 3076 5552}%
\special{pa 3076 5858}%
\special{fp}%
\special{sh 1}%
\special{pa 3076 5858}%
\special{pa 3096 5792}%
\special{pa 3076 5806}%
\special{pa 3056 5792}%
\special{pa 3076 5858}%
\special{fp}%
\put(31.8500,-57.7100){\makebox(0,0)[lb]{Blow down along the $(-1)$-curve $D_{\infty}^{(1)}$ to a nonsingular point}}%
%
\special{pn 8}%
\special{pa 2970 6434}%
\special{pa 3260 6256}%
\special{dt 0.045}%
%
\special{pn 8}%
\special{pa 2980 6668}%
\special{pa 3280 6494}%
\special{dt 0.045}%
%
\special{pn 8}%
\special{pa 2310 6074}%
\special{pa 2310 6512}%
\special{dt 0.045}%
%
\special{pn 8}%
\special{pa 2410 6176}%
\special{pa 2410 6614}%
\special{dt 0.045}%
%
\special{pn 8}%
\special{pa 2520 6278}%
\special{pa 2520 6718}%
\special{dt 0.045}%
%
\special{pn 8}%
\special{pa 2640 6418}%
\special{pa 2640 6858}%
\special{dt 0.045}%
%
\special{pn 8}%
\special{pa 2780 6502}%
\special{pa 2780 6940}%
\special{dt 0.045}%
\put(18.4000,-61.7500){\makebox(0,0)[lb]{$(D_{0}^{(1)})^2=-4$}}%
\put(28.4000,-74.8200){\makebox(0,0)[lb]{$(D_{1}^{(1)})^2=-1$}}%
%
\special{pn 8}%
\special{pa 2900 7000}%
\special{pa 3430 7000}%
\special{dt 0.045}%
%
\special{pn 20}%
\special{pa 3066 7622}%
\special{pa 3066 7926}%
\special{fp}%
\special{sh 1}%
\special{pa 3066 7926}%
\special{pa 3086 7860}%
\special{pa 3066 7874}%
\special{pa 3046 7860}%
\special{pa 3066 7926}%
\special{fp}%
\put(31.7500,-78.3900){\makebox(0,0)[lb]{Blow down along the $(-1)$-curve $D_{1}^{(1)}$ to a nonsingular point}}%
%
\special{pn 20}%
\special{pa 2200 8400}%
\special{pa 3990 8400}%
\special{fp}%
%
\special{pn 8}%
\special{pa 2280 8042}%
\special{pa 2280 8680}%
\special{dt 0.045}%
%
\special{pn 8}%
\special{pa 2400 8042}%
\special{pa 2400 8680}%
\special{dt 0.045}%
%
\special{pn 8}%
\special{pa 2550 8058}%
\special{pa 2550 8696}%
\special{dt 0.045}%
%
\special{pn 8}%
\special{pa 2700 8074}%
\special{pa 2700 8712}%
\special{dt 0.045}%
%
\special{pn 8}%
\special{pa 2880 8074}%
\special{pa 2880 8712}%
\special{dt 0.045}%
\put(40.3000,-85.0800){\makebox(0,0)[lb]{$(D_{0}^{(1)})^2=-3$}}%
%
\special{pn 20}%
\special{pa 2200 6238}%
\special{pa 3270 7110}%
\special{fp}%
%
\special{pn 8}%
\special{pa 3120 6208}%
\special{pa 3120 7156}%
\special{dt 0.045}%
%
\special{pn 8}%
\special{pa 3570 8416}%
\special{pa 3090 8680}%
\special{dt 0.045}%
\special{pa 3570 8416}%
\special{pa 3570 8696}%
\special{dt 0.045}%
%
\special{pn 8}%
\special{pa 3570 8432}%
\special{pa 3990 8664}%
\special{dt 0.045}%
\end{picture}%
\label{fig:E8figure2}
\caption{Resolution of accessible singular points}
\end{figure}

Note that $\{(u_{\nu},v_{\nu})|v_{\nu}=0\} \subset D_{\nu}^{(1)}$ for $\nu=0,1,\infty$. By a direct calculation, we see that the system \eqref{E811} has eight accessible singular points $a_{\nu}^{(1)}$ for $\nu=1,2,3,4,5,6,7,8$ in $D_{\nu}^{(1)} \cong {\Bbb P}^1 \ (\nu=0,1,\infty)$.
\begin{align}
\begin{split}
&a_{1}^{(1)}=\{(u_{0},v_{0})=(\alpha_1,0)\} \in D_{0}^{(1)}, \quad a_{2}^{(1)}=\{(u_{0},v_{0})=(\alpha_1+\alpha_2,0)\} \in D_{0}^{(1)},\\
&a_{3}^{(1)}=\{(u_{0},v_{0})=(\alpha_1+\alpha_2+\alpha_3,0)\} \in D_{0}^{(1)}, \quad a_{4}^{(1)}=\{(u_{0},v_{0})=(\alpha_1+\alpha_2+\alpha_3+\alpha_4,0)\} \in D_{0}^{(1)},\\
&a_{5}^{(1)}=\{(u_{0},v_{0})=(\alpha_1+\alpha_2+\alpha_3+\alpha_4+\alpha_5,0)\} \in D_{0}^{(1)}, \quad a_{6}^{(1)}=\{(u_{1},v_{1})=(\alpha_6,0)\} \in D_{1}^{(1)},\\
&a_{7}^{(1)}=\{(u_{1},v_{1})=(\alpha_6+\alpha_7,0)\} \in D_{1}^{(1)},\\
&a_{8}^{(1)}=\{(u_{\infty},v_{\infty})=(\alpha_8,0)\} \in D_{\infty}^{(1)}
\end{split}
\end{align}
Let us perform blowing-ups at $a_{j}^{(1)}$, and denote $D_{j}^{(2)}$ for the exceptional curves, respectively. We take seven coordinate systems $(W_j,V_j)$ around the points at infinity of $D_{j}^{(2)}$ for $j=1,2,3,4,5,6,7,8$, where
\begin{align}
\begin{split}
&(W_{1},V_{1})=\left(\frac{u_0-\alpha_1}{v_0},v_0 \right),\\
&(W_{2},V_{2})=\left(\frac{u_0-(\alpha_1+\alpha_2)}{v_0},v_0 \right),\\
&(W_{3},V_{3})=\left(\frac{u_0-(\alpha_1+\alpha_2+\alpha_3)}{v_0},v_0 \right),\\
&(W_{4},V_{4})=\left(\frac{u_0-(\alpha_1+\alpha_2+\alpha_3+\alpha_4)}{v_0},v_1 \right),\\
&(W_{5},V_{5})=\left(\frac{u_0-(\alpha_1+\alpha_2+\alpha_3+\alpha_4+\alpha_5)}{v_0},v_0 \right),\\
&(W_{6},V_{6})=\left(\frac{u_1-\alpha_6}{v_1},v_1 \right),\\
&(W_{7},V_{7})=\left(\frac{u_{1}-(\alpha_6+\alpha_7)}{v_{1}},v_{1} \right),\\
&(W_{8},V_{8})=\left(\frac{u_{\infty}-\alpha_8}{v_{\infty}},v_{\infty} \right).
\end{split}
\end{align}
For the strict transform of $D^{(0)}$, $D_{\nu}^{(1)}$ and $D_{j}^{(2)}$ by the blowing-ups, we also denote by same symbol, respectively.  Here, the self-intersection number of $D^{(0)}, D_{\nu}^{(1)}$ is given by
\begin{equation}
(D^{(0)})^2=-1. \quad (D_{0}^{(1)})^2=-6, \quad (D_{1}^{(1)})^2=-3, \quad (D_{\infty}^{(1)})^2=-2.
\end{equation}
In order to obtain a minimal compactification of the space of initial conditions, we must blow down along the $(-1)$-curves $D^{(0)} \cong {\Bbb P}^1$ to a nonsingular point. For the strict transform of $D_{\nu}^{(1)}$ and $D_{j}^{(2)}$ by the blowing-down, we also denote by same symbol, respectively. Here, the self-intersection number of $D_{\nu}^{(1)}$ is given by
\begin{equation}
(D_{0}^{(1)})^2=-5, \quad (D_{1}^{(1)})^2=-2, \quad (D_{\infty}^{(1)})^2=-1.
\end{equation}
We must blow down again along the $(-1)$-curve $D_{\infty}^{(1)} \cong {\Bbb P}^1$ to a nonsingular point. For the strict transform of $D_{\nu}^{(1)}$ and $D_{j}^{(2)}$ by the blowing-down, we also denote by same symbol, respectively. Here, the self-intersection number of $D_{\nu}^{(1)}$ is given by
\begin{equation}
(D_{0}^{(1)})^2=-4, \quad (D_{1}^{(1)})^2=-1.
\end{equation}
We must blow down again along the $(-1)$-curve $D_{1}^{(1)} \cong {\Bbb P}^1$ to a nonsingular point. For the strict transform of $D_{\nu}^{(1)}$ and $D_{j}^{(2)}$ by the blowing-down, we also denote by same symbol, respectively. Here, the self-intersection number of $D_{0}^{(1)}$ is given by
\begin{equation}
(D_{0}^{(1)})^2=-3.
\end{equation}

Let ${\tilde S} \cdots \rightarrow {\Sigma_2}$ be the composition of above eleven times blowing-ups and three times blowing-downs. Then, we see that the canonical divisor class $K_{{\tilde S}}$ of ${\tilde S}$ is given by
\begin{equation}
K_{{\tilde S}}:=-D_{0}^{(1)},
\end{equation}
where the self-intersection number of $D_{0}^{(1)} \cong {\Bbb P}^1$ is given by
\begin{equation}
(D_{0}^{(1)})^2=-3.
\end{equation}
We see that ${\tilde S}-(-K_{{\tilde S}})_{red}$ is covered by ten Zariski open sets
\begin{align}
\begin{split}
& \rm{Spec} \ {\Bbb C}[W_{j},V_{j}] \quad (j=1,2,3,4,5,6,7,8),\\
& \rm{Spec} \ {\Bbb C}[z_0,w_0],\\
& \rm{Spec} \ {\Bbb C}[z_1,w_1].
\end{split}
\end{align}
The relations between $(W_{j},V_{j})$ and $(x_j,y_j)$ are given by
\begin{equation}
(-W_{j},V_{j})=(x_j,y_j) \quad (j=1,2,3,4,5,6,7,8).
\end{equation}
We see that the pole divisor of the symplectic 2-form $dp \wedge dq$ coincides with $(-K_{{\tilde S}})_{red}$. Thus, we have completed the proof of Theorems 13.1 and 13.2. \qed


\begin{thebibliography}{99}

\bibitem[1]{1} P. Painlev\'e, {\em M\'emoire sur les \'equations diff\'erentielles dont l'int\'egrale g\'en\'erale est uniforme}, Bull. Soci\'et\'e Math\'ematique de France. {\bf 28} (1900),  201--261.

\bibitem[2]{2} P. Painlev\'e, {\em Sur les \'equations diff\'erentielles du second ordre et d'ordre sup\'erieur dont l'int\'egrale est uniforme}, Acta Math. {\bf 25} (1902), 1--85. 

\bibitem[3]{3} B. Gambier, {\em Sur les \'equations diff\'erentielles du second ordre et du premier degr\'e dont l'int\'egrale g\'en\'erale est \`a points critiques fixes}, Acta Math. {\bf 33} (1910), 1--55.

\bibitem[4]{O2} K. Okamoto, {\em Studies on the Painlev\'e equations, I}, Ann. Mat. Pura Appl., {\bf 146} (1987), 337--381; II, Jap. J. Math., {\bf 13} (1987), 47--76; III, Math. Ann., {\bf 275} (1986), 221--256; IV, Funkcial. Ekvac., {\bf 30} (1987), 305--332.


\bibitem[5]{O3} K. Okamoto, {\em Sur les 
feuilletages associ\'es aux \'equations du second ordre \`a points critiques fixes de P. Painlev\'e, Espaces des conditions initiales}, Japan. 
J. Math. {\bf 5} (1979), 1--79.  

\bibitem[6]{O4} K. Okamoto, {\em Polynomial Hamiltonians associated with Painlev\'e equations}, I, II, Proc. Japan Acad. {\bf 56} (1980),  264--268; ibid, 367--371.


\bibitem[7]{Sakai} H. Sakai, {\em Rational surfaces associated with affine root systems and geometry of the Painlev\'e equations}, Commun. Math. Phys. {\bf 220} (2001), 165--229. 

 

\bibitem[8]{Cosgrove1} C. M. Cosgrove and G. Scoufis,
{\em Painlev\'e classification of a class of differential equations of the second order and second degree}, Studies in Applied Mathematics. {\bf 88} (1993), 25-87.

\bibitem[9]{Cosgrove2} C. M. Cosgrove,
{\em All binomial-type Painlev\'e equations of the second order and degree three or higher}, Studies in Applied Mathematics. {\bf 90} (1993), 119-187.

\bibitem[10]{Sasano} Y. Sasano,
{\em Symmetry in the Painlev\'e systems and their extensions to four-dimensional systems}, Funkcialaj Ekvacioj, {\bf 51} (2008), 351-369.

\bibitem[11]{Joshi} Joshi, N, Kitaev, A.V. and Treharne, P.A.,
{\em On the linearization of the Painlev\'e III-VI equations and reductions of the three-wave resonant system}, (2007), arXiv:0706.1750.

\bibitem[12]{Mazzocco} Mazzocco, M.,
{\em Painlev\'e sixth equation as isomonodromic deformations equation of an irregular system},  in The Kowalevski property, CRM Proc. Lecture Notes {\bf 32} (2002), Providence, RI 219-238.



\bibitem[13]{13} F. Bureau, 
{\em Integration of some nonlinear systems of ordinary differential equations}, 
Annali di Matematica. {\bf 94} (1972), 345--359. 

\bibitem[14]{14} J. Chazy, 
{\em Sur les \'equations diff\'erentielles dont l'int\'egrale g\'en\'erale est uniforme et admet des singularit\'es essentielles mobiles}, 
Comptes Rendus de l'Acad\'emie des Sciences, Paris. {\bf 149} (1909), 563--565. 

\bibitem[15]{15} J. Chazy, 
{\em Sur les \'equations diff\'erentielles dont l'int\'egrale g\'en\'erale poss\'ede une coupure essentielle mobile }, 
Comptes Rendus de l'Acad\'emie des Sciences, Paris. {\bf 150} (1910), 456--458. 


\bibitem[16]{16} J. Chazy, 
{\em Sur les \'equations diff\'erentielles du trousi\'eme ordre et d'ordre sup\'erieur dont l'int\'egrale a ses points critiques fixes}, 
Acta Math. {\bf 34} (1911), 317--385. 


\bibitem[17]{17} E. L. Ince, 
{\em Ordinary differential equations}, Dover Publications, New York, (1956). 

\bibitem[18]{Jimbo} M. Jimbo and T. Miwa, 
{\em Monodromy preserving deformation of linear ordinary differential equations with rational coefficients II}, Physica D 2, 407--448, (1981).

\bibitem[19]{T1} T. Shioda and K. Takano, {\em On some Hamiltonian structures of Painlev\'e systems I}, Funkcial. Ekvac. {\bf 40} (1997), 271--291.

\end{thebibliography}
\end{document}